\documentclass{elsarticle}
\journal{}
\usepackage{color}
\usepackage{amsfonts,enumerate,amsmath,amssymb}
\usepackage{mathrsfs}
\usepackage{geometry}
\def\qed{\strut\hfill $\Box$}
\newtheorem{thm}{Theorem}[section]

\newtheorem{lem}[thm]{Lemma}

\newtheorem{rem}[thm]{Remark}
\newtheorem{defn}[thm]{Definition}

\newcommand{\lemref}[1]{Lemma~{\rm \ref{#1}}}

\def\para#1{\vskip .4\baselineskip\noindent{\bf #1}}
\numberwithin{equation}{section}
\begin{document}
\begin{frontmatter}	
\title{Pathwise unique solutions and stochastic averaging for mixed stochastic partial differential equations driven by fractional Brownian motion and Brownian motion}

\author[mymainaddress,mythirdaddress]{Bin \textsc{Pei}}
\ead{binpei@nwpu.edu.cn}

\author[mysecondaryaddress]{Yuzuru \textsc{Inahama}}
\ead{inahama@math.kyushu-u.ac.jp}

\author[mymainaddress]{Yong \textsc{Xu}}
\ead{hsux3@nwpu.edu.cn}

\address[mymainaddress]{School of Mathematics and Statistics, Northwestern Polytechnical University, Xi'an, 710072, China}
\address[mysecondaryaddress]{Faculty of Mathematics,
	Kyushu University, Fukuoka, 819-0395, Japan}
\address[mythirdaddress]{
Innovation Center, NPU-Chongqing, Chongqing, 401135, China}

\begin{abstract}
This paper is devoted to a system of stochastic partial differential equations (SPDEs) that have
a slow component driven by fractional Brownian motion (fBm)
with the Hurst parameter $H >1/2$
and a fast component
driven by fast-varying diffusion. It improves previous work in two aspects: Firstly, using a stopping time technique and an approximation of the fBm, we prove an existence and uniqueness theorem for a class of mixed SPDEs driven by both fBm and Brownian motion;  Secondly, an averaging principle in the mean square sense for SPDEs driven by fBm subject to an additional fast-varying diffusion process is established. To carry out these improvements, we combine
the pathwise approach based on the generalized Stieltjes integration theory with the It\^o stochastic calculus.  Then, we obtain a desired limit process of the slow component which strongly relies on an invariant measure of the fast-varying diffusion process.
\end{abstract}		
\end{frontmatter}

\section{Introduction}\label{se-1}
It is widely known that there are many phenomena the well-studied
theory of semimartingales cannot describe.
For example, telecommunication connections, climate, weather derivatives and other objects have long memory \cite{decreusefond1998fractional,Mandelbrot1968,Taqqu2003fractional}.
Brownian motion (Bm) with independent increments which has no memory turns out to be insufficient to describe this effect. Another example could be that the concept of turbulence in
hydrodynamics can be described with the help of stationary (dependent)
increments \cite{Monin1967Statistical}. Thus, the long-range dependence properties of fracional Brownian motion (fBm) make this process a suitable candidate to describe this kind of phenomena.

For $H \in  (0, 1)$, a continuous centered Gaussian process
$\beta^H =(\beta^H (t))_{t \ge 0}$ with the covariance function
$$\mathbb{E} [\beta^{H}(t) \beta^{H}(s)]=\frac{1}{2}\left(|t|^{2 H}+|s|^{2 H}-|t-s|^{2 H}\right), \qquad t, s \ge 0,$$
is called one-dimensional fBm \cite{Mandelbrot1968} with
the Hurst parameter $H$.
An fBm differs significantly from an Bm
and semimartingales. It is characterized by the stationarity of its  (dependent) increments
and long-memory property only for $H\in(\frac{1}{2},1)$. In the case $H\in(0,\frac{1}{2})$ it is a process with short memory.
Note that if $H\neq \frac{1}{2}$,
an fBm is not a semimartingale nor a Markov process.

Now, we recall the definition of an infinite-dimensional fBm
following \cite{GarridoAtienza2010dcds,maslowski2003evolution}.
Let $(V,  \langle \cdot,\cdot \rangle)$ be a separable Hilbert space.
Its norm is denoted by $ |\cdot|$.
For a sequences $\{\lambda_{i}\}_{i \in \mathbb{N}}$
of positive real numbers with $\sum_{i=1}^{\infty} {\lambda_{i}}<\infty$
and an orthonormal basis
$\{e_{i}\}_{i \in \mathbb{N}}$ of $V$,
a $V$-valued fBm $B^{H}$ is defined by
$$
B^{H}_t=\sum_{i=1}^{\infty} \sqrt{\lambda_{i}} e_{i} \beta_{i}^{H}(t), \qquad t \ge 0,
$$
where $\left\{\beta_{i}^{H}\right\}_{i \in \mathbb{N}}$ is a sequence of  independent one-dimensional fBm's.
It is known that the right hand side is convergent in $L^{2}$
for every $t$ and has a continuous modification in $t$.

Let
$(\Omega, \mathscr{F}, \{\mathscr{F}_t\}_{t\ge 0}, \mathbb{P})$
be a suitable probability space with a filtration
satisfying the usual condition.
We assume that
a $V$-valued $\{\mathscr{F}_t\}$-Brownian motion $W$
and
$\{\mathscr{F}_t\}$-adapted one-dimensional fBm's
$\beta^H_i, i \in {\mathbb N},$
are defined on this probability space.
We further assume that $W$ and $\beta^H_i, i \in {\mathbb N},$
are all independent.
(For the definition of
a $V$-valued Brownian motion, see \cite[Proposition 4.3]{da2014stochastic} where is called $Q$-Wiener process.)

From now on we assume $H \in (\frac12, 1)$
and work on the time interval $[0,T]$,
where $T >0$ is arbitrary but fixed.
Let $A$ be the infinitesimal generator of
an analytic semigroup $S$ on $V$ and assume that
$-A$ has discrete
spectra $0 < \bar{\lambda}_1 < \bar \lambda_{2} < \cdots <\bar \lambda_{k}  < \cdots $ and
$\lim_{k\rightarrow \infty} \bar \lambda_k =\infty$.

This paper firstly will prove an existence and uniqueness theorem for a class of mixed stochastic partial differential equations (SPDEs)
driven by both fBm and Bm with a given initial value $u_0$,
which is given by
\begin{eqnarray}\label{spde-uniq}
\mathrm{d} u_t =(A u_t+f(u_t))\mathrm{d}t+\sigma(u_t)\mathrm{d}W_t+g(u_t) \mathrm{d}B^H_t,
\qquad
0 \le  t \le T.
\end{eqnarray}
Precise conditions on the
nonlinear coefficients $f,\sigma, g$ will be given in Section \ref{se-3}.

The idea of this part is based on a pathwise approach developed by Z\"ahle \cite{Zahle1998}, who defined the stochastic integral with respect to fBm based on a sort of generalized integration by parts formula with respect to fractional derivatives. Nualart and  R\u{a}\c{s}canu \cite{Rascanu2002} and Garrido-Atienza et al. \cite{Random2010garr} investigated stochastic differential equations (SDEs) in finite dimension. Infinite-dimensional equations were
treated with the same success as finite-dimensional ones, e.g. Tindel et al. \cite{tindel2003stochastic} and Garrido-Atienza et al. \cite{GarridoAtienza2010dcds,Garrido-Atienzalu2010}. Pathwise solutions of this kind of equation without Bm term ($\sigma(\cdot)=0$ in Eq. (\ref{spde-uniq})) were studied in
Maslowski and Nualart \cite{maslowski2003evolution} and Garrido-Atienza, Lu, and Schmalfuss \cite{Garrido-Atienzalu2010}, and recently by Chen, Gao,
Garrido-Atienza, and Schmalfuss \cite{chen2014pathwise} when the stochastic evolution equations are driven by a H\"older continuous
function with H\"older exponent in $(\frac{1}{2},1)$ and with nontrivial multiplicative
noise. Guerra and Nualart \cite{Guerra2008} proved an existence and uniqueness
theorem of solutions to multidimensional SDEs driven by fBm with Hurst parameter $H>\frac{1}{2}$ and Bm. Using the theory of Wiener integral, Caraballo, Garrido-Atienza and Taniguchi \cite{caraballo2011existence} investigated the existence and exponential behavior of solutions to stochastic delay evolution equations with an additive fractional noise.

However, the method proposed in \cite{Guerra2008} fails for
the infinite-dimensional case and the method in \cite{caraballo2011existence} fails
for the multiplicative fractional noise case.
The main difference (and, of course, difficulty) is that we
cannot apply directly the existence and uniqueness results in \cite{Garrido-Atienzalu2010,Guerra2008} and \cite{caraballo2011existence,chen2014pathwise}. Thus, to close this gap, as one of  two main results of our paper, we  obtain pathwise unique solutions to Eq. (\ref{spde-uniq}) relying on a pathwise approach, a stopping time technique and an approximation for the fractional noise (See Theorem \ref{mixuniq}).

Then, as the second main result,
this paper will establish an averaging principle in the mean square sense for a class of SPDEs driven by fBm subject to an additional fast-varying diffusion process, which is given by
\begin{eqnarray}\label{spde-fs}
\begin{cases}
\mathrm{d} X^{\varepsilon}_t&=(AX^{\varepsilon}_t+b( X^{\varepsilon}_t,Y^{\varepsilon}_t))\mathrm{d}t+g(X^{\varepsilon}_t)\mathrm{d}B^H_t,\\ \mathrm{d} Y^{\varepsilon}_t& =\frac{1}{\varepsilon}(AY^{\varepsilon}_t+F(X^{\varepsilon}_t,Y^{\varepsilon}_t))\mathrm{d}t+\frac{1}
{\sqrt{\varepsilon}}G(X^{\varepsilon}_t,Y^{\varepsilon}_t)\mathrm{d}W_t,
\end{cases}
\end{eqnarray}
where $b,g,F,G$ are nonlinear coefficients and $X^{\varepsilon}_0=X_{0},Y^{\varepsilon}_0=Y_{0}$ are initial values.
The parameter $0<\varepsilon \ll 1$ represents the ratio between the natural time scale of the $X_{t}^{\varepsilon}$ and $Y_{t}^{\varepsilon}$ variables. For more precise setting and assumptions, see Section 4.

The theory of stochastic averaging principles has been studied extensively (see for instance the
paper by Khasminskii \cite{khas1966limit} and
Freidlin
and Wentzell \cite{freidlin2012random},  or the recent paper by Xu et al. \cite{xu2011averaging,xu2014,xu2015approximation,Xu2015Stochastic,Xu2017Stochastic}, Liu \cite{liu2010strong}, Liu, R\"{o}ckner, Sun and Xie \cite{Liu2019} and Thompson, Kuske, and Monahan \cite{thompson2015stochastic}) and is used in many applications.
Cerrai and Freidlin \cite{cerrai2009averaging} developed  stochastic averaging
principles for two-time-scale stochastic reaction diffusion equations whose additive noise is included in the fast motion.
In this infinite-dimensional setting, there are also interesting papers such as Br\'{e}hier \cite{brehier2012strong}, Xu and Miao \cite{xu2015strong}, Fu et al. \cite{fu2011strong,fu2015strong}, Bao, Yin, and Yuan \cite{Bao2016Two} and Sun and Zhai \cite{sun2020baveraging}.  However, the
literature concerning fast-slow mixed SPDEs driven by multiplicative fractional noise is still, to some extent, in its infancy.  Pei, Xu and Yin  \cite{peistochastic2017} established an averaging principle for a system of SPDEs that have
a slow component driven by an additive fractional noise and a fast component
driven by fast-varying diffusion. Pei et al. \cite{pei2018sd,Averaging2018pei} examined
averaging principles for
SPDEs
driven by an additive fractional noise
with two-time-scale Markovian switching processes. But, till now, in multiplicative fractional noise case, only the averaging results for SDE cases were obtained. Hairer and Li \cite{Hairer2019averaging} considered slow-fast systems where the slow system is driven by fBm and proved the convergence to the averaged solution took place in probability. Very recently, Pei, Inahama and Xu \cite{pei2020averaging} answered affirmatively that an averaging principle still holds for fast-slow mixed SDEs if
disturbances involve both Bm and long-range dependence modeled by
fBm $H\in(\frac{1}{2},1)$ in the mean square sense.

To the best of our
knowledge, the second part of our paper is the first attempt to study stochastic averaging for fast-slow mixed SPDEs driven by multiplicative fractional noise. The main goal of this part is to generalize the results in  \cite{pei2020averaging,pei2020convergence,peistochastic2017} by using directly a pathwise approach to deal
with multiplicative fractional noise term. In order to reach this
objective, we shall borrow the construction of stochastic integral with respect to infinite-dimensional fBm given in \cite{Garrido-Atienzalu2010,maslowski2003evolution} and the stopping time technique to control the fBm term given in \cite{Mishura2012mixed} which will be recalled in Section \ref{se-2} and Section \ref{se-3}, respectively (See Theorem \ref{mainthm}).

The paper is organized as follows. Section \ref{se-2} presents some necessary
notations and assumptions. Section \ref{se-3} presents pathwise unique solutions for a class of mixed SPDEs driven by fBm and Bm.
An averaging principle for fast-slow mixed SPDEs driven by fBm subject to an additional fast-varying diffusion process is then established in Section \ref{se-4}. Appendix A provides the arguments of the ergodicity for the fast component in which the slow component is kept
frozen. Some technical complements are included in
Appendix B.

\section{Preliminaries}\label{se-2}
Although the results on fractional calculus and stochastic integrals with respect to the one-dimensional fBm $\beta^H$ and $V$-valued infinite-dimensional fBm $B^H$ have
already been done in the recent paper \cite{Biagini2008stochastic,GarridoAtienza2010dcds,Garrido-Atienzalu2010,maslowski2003evolution,mishura2008stochastic}, we present them here for the sake of completeness.

For $T>0$ and $0<\alpha<\frac{1}{2}$, we denote $W^{\alpha,1}(0,T;V)$, the space of measurable functions $h:[0,T]\to V$ such that
\begin{eqnarray*}
	|h|_{\alpha,1}:=\int_0^T\Big(\frac{|{h(s)}|}{s^\alpha}+\int_0^s\frac{|{h(s)-h(r)}|}{(s-r)^{\alpha+1}}\mathrm{d}r\Big)\mathrm{d}s<\infty.
\end{eqnarray*}	

For $a<t<c$ the Weyl derivatives are given by
\begin{eqnarray*}	
	D_{a+}^\alpha h(t)&=&\frac{1}{\Gamma(1-\alpha)}\Big(\frac{h(t)}{(t-a)^\alpha}+\alpha\int_a^t\frac{h(t)-h(\zeta)}{(t-\zeta)^{\alpha+1}}\mathrm{d}\zeta\Big),\cr
	D_{c-}^{1-\alpha} l_{c-} (t)&=&\frac{(-1)^{1-\alpha}}{\Gamma(\alpha)}\Big(\frac{l(t)-l(c)}{(c-t)^{1-\alpha}}+(1-\alpha) \int_t^c\frac{l(t)-l(\zeta)}{(\zeta-t)^{2-\alpha}}\mathrm{d}\zeta\Big),
\end{eqnarray*}	
where, for $0\le a < c\le T, l_{c-}(r):=l(r)-l(c),$ and $\Gamma$ denotes the Gamma function. Then, according to Z\"{a}hle \cite{Zahle1998}, for $h\in W^{\alpha,1}(0,T;V), 0 \le s<t \le T$, the generalized Stieltjes integrals
\begin{align}\label{g-s}
\int_0^T h(r) \mathrm{d}l(r)=&(-1)^\alpha\int_0^TD_{0+}^{\alpha}h(r)D_{T-}^{1-\alpha}l_{T-}(r)\mathrm{d}r,\\
\int_s^t h(r)\mathrm{d}l(r)=&\int_0^T h(r) \mathbf{1}_{(s,t)}\mathrm{d}l(r),
\end{align}
are defined. In addition, the integral (\ref{g-s}) exists and has the following estimatate
\begin{eqnarray*}
	\Big|\int_0^T h(t) \mathrm{d}l(t)\Big| \leq \frac{\|l\|_{\alpha,0,T}}{\Gamma(1-\alpha)\Gamma(\alpha)} |h|_{\alpha,1},
\end{eqnarray*}	
where
$$\|l\|_{\alpha,0,T}:=\sup_{0\le s<t\le T}\Big(\frac{|l(t)-l(s)|}{(t-s)^{1-\alpha}}+\int_s^t\frac{|l(\zeta)-l(s)|}{(\zeta-s)^{2-\alpha}}\mathrm{d}\zeta\Big)<\infty.$$
For the sake of shortness, we denote $\Lambda_\alpha^{0,T}(l):=\frac{\|l\|_{\alpha,0,T}}{\Gamma(1-\alpha)\Gamma(\alpha)}.$

Let $H\in(\frac{1}{2},1)$, take a parameter $\alpha \in( 1-H,\frac{1}{2})$ which will be fixed througout this paper.
For $h\in W^{\alpha,1}(0,T;V)$ the integral
$$\int_{0}^{T}h(s)\mathrm{d}\beta^H(s)$$
will be understood in the sense of definition (\ref{g-s}) pathwise, which makes sense because $\Lambda_\alpha^{0,T}(\beta^H)<\infty$ a.s. (cf. \cite{Rascanu2002}).

Let $L(V)$ denote the space of linear bounded operators on $V$ and let $G:\Omega \times[0,T]\to L(V)$ be an operator valued map such that $G(\omega,\cdot) e_i\in W^{\alpha,1}(0,T;V)$ for each $i\in \mathbb{N}$ and almost $\omega\in\Omega$. We define
\begin{eqnarray}\label{infbm}
\int_0^TG(s)\mathrm{d} B_s^H:=\sum_{i=1}^\infty \sqrt{\lambda_i} \int_0^TG(s)e_i \mathrm{d} \beta_i^H(s),
\end{eqnarray}		
where the convergence of the sums in (\ref{infbm}) is understood as $\mathbb{P}$-a.s. convergence in $V$.

From now on, to make the pathwise integral (\ref{infbm}) well-defined,  we assume $\sum_{i=1}^\infty\sqrt{\lambda_i}<\infty$.
The following result can be found in \cite[Proposition 2.1]{maslowski2003evolution}.
\begin{rem}\label{itofbm}{\rm
		Assume that $\sum_{i=1}^\infty\sqrt{\lambda_i}<\infty.$ Then there exists $\Omega_1 \subset \Omega$, $\mathbb{P}(\Omega_1)=1$, such that the pathwise integral Eq. (\ref{infbm}) is well-defined on $\Omega_1$ for each $G:\Omega_1\times [0,T]\to L(V)$ satisfying $G(\omega,\cdot)e_i\in W^{\alpha,1}(0,T;V),$ for $   \omega\in \Omega_1,$ such that $\sup_{i\in\mathbb{N}}|G(\omega,\cdot)e_i|_{\alpha,1}<\infty.$ In addition
		\begin{eqnarray*}
			\Big|\int_0^TG(s) \mathrm{d} B^H_s\Big| \le \Lambda_{\alpha,B^H}^{0,T}\sup_{i\in\mathbb{N}} |G(\cdot)e_i|_{\alpha,1}, \qquad
			\omega \in \Omega_1,
		\end{eqnarray*}
		where $\Lambda_{\alpha,B^H}^{0,T}:=\sum_{i=1}^\infty\sqrt{\lambda_i} \Lambda_\alpha^{0,T}(\beta_i^H)$. Note that $\Lambda_{\alpha,B^H}^{0,T}$ is finite a.s. }
\end{rem}

We recall the following two auxiliary technical lemmas from \cite{GarridoAtienza2010dcds}.
\begin{lem}\label{inq-s}
	For any positive constants $a, d$, if $a+d-1>0$ and $a<1 $, one has
	\begin{align*}
	\int_{0}^{r}(r-s)^{-a}(t-s)^{-d} d s &\leq (t-r)^{1-a-d} B(1-a, d+a-1),\\
	\int_{r}^{t}(s-r)^{-a}(t-s)^{-d} d s &\leq (t-r)^{1-a-d}B(1-a, d+a-1),
	\end{align*}
	where $r \in(0, t)$ and $B$ is the Beta Function.
\end{lem}
\begin{lem} \label{inq-rho}
	For any non-negative a and $d$ such that $a+d<1,$ and for any  $\rho \geq 1$, there exists a positive constant $C$ such that
	\begin{eqnarray*}
		\int_{0}^{t} e^{-\rho(t-r)}(t-r)^{-a} r^{-d} d r \leq C \rho^{a+d-1}.
	\end{eqnarray*}
	In addition, for $d \leq 0$ and $0 \leq a<1,$ and for any $\rho \geq 1$, we have
	\begin{eqnarray*}
		\int_{0}^{t} e^{-\rho(t-r)}(t-r)^{-a} r^{-d} d r \leq \Gamma(1-a) t^{-d} \rho^{a-1}.
	\end{eqnarray*}
\end{lem}

Please note that $C$ and $C_\ast$ denote certain positive constants that may depend on the parameters $\alpha,\beta,T$ and the initial values and vary from line to line. $C_\ast$
is used to emphasize that the constant depends on the corresponding parameter $\ast$ which is one or more than one parameter.

\section{Mixed SPDEs driven by fBm and Bm}\label{se-3}
This section will prove an existence and uniqueness theorem  for the mixed SPDEs driven by both fBm and Bm (\ref{spde-uniq}). Let $V_\beta, \beta\geq 0$, denote the domain of the fractional power $(-A)^\beta$ equipped with the graph norm $|x|_{V_\beta}:=|(-A)^{\beta}x|, x\in V_\beta$. For shortness, denote, $|x|_\beta:=|x|_{V_\beta}.$ We recall here some properties of the analytic semigroup, which will be used later in our analysis. For $0 \leq   \gamma \leq \varsigma \leq 1$ and $\upsilon\in[0,1), \mu\in(0,1-\upsilon)$, there exists a constant $C>0$, such that for $0\leq s <t \leq T$, we have
\begin{align}\label{semi1}
|S_t|_{L(V_\gamma, V_{\varsigma})} &\leq C t^{-\varsigma+\gamma}e^{-\bar{\lambda}_1 t},\\
\label{semi2}
|S_{t-s}-\mathrm{id}|_{L(V_{\upsilon+\mu}, V_{\upsilon})} &\leq C (t-s)^{\mu}.
\end{align}
We also note that, for $\varrho,\nu \in(0,1]$ and $0\leq \nu<\gamma+\varrho$, there exists $C>0$ such that
for $0\leq q \leq r \leq s \leq t$, we derive
\begin{eqnarray}\label{semi3}
|S_{t-r}-S_{t-q}|_{L\left(V_{\nu}, V_{\gamma}\right)} \leq C(r-q)^{\varrho}(t-r)^{-\varrho-\gamma+\nu},
\end{eqnarray}
and
\begin{eqnarray}\label{semi4}
|S_{t-r}-S_{s-r}-S_{t-q}+S_{s-q}|_{L(V,V)} \leq C(t-s)^{\varrho}(r-q)^{\nu}(s-r)^{-(\varrho+\nu)}.
\end{eqnarray}
From now on, we use the symbol $\|\cdot\|$ to denote $|\cdot|_{L(V,V)}$ for shortness.

Let $f: V \rightarrow V$ and $\sigma: V \rightarrow L_2(V)$ be measurable and satisfy Lipschitz and linear growth conditions with constants $L_{f}$ and $L_{\sigma}$ respectively, where $L_{2}(V)$ is the family of Hilbert-Schmidt operators from $V$ to itself. We assume that 
$g: V \rightarrow L(V)$ is Fr\'echet $C^1$ and that $g$ and $g^{\prime} : V \rightarrow L(V, L(V))$ are Lipschitz continuous with constants $L_{g},M_{g}$ in the following senses:
\begin{eqnarray}\label{g1}
\sup _{i \in \mathbb{N}}\left|g\left(v_{1}\right) e_{i}-g\left(v_{2}\right) e_{i}\right| \leq L_{g}\left|v_{1}-v_{2}\right|,
\end{eqnarray}
\begin{eqnarray}\label{g2}
\sup _{i \in \mathbb{N}}\left|g^{\prime}\left(v_{1}\right) e_{i}-g^{\prime}\left(v_{2}\right) e_{i}\right|_{L(V)} \leq M_{g}\left|v_{1}-v_{2}\right|,
\end{eqnarray}
where $\left\{e_{i}\right\}_{i \in \mathbb{N}}$ is the complete orthonormal basis in $V$.
\begin{rem}{\rm
		For $v_1, v_2, u_1, u_2\in V$, there exist $c_1,c_2>0$, such that
		\begin{eqnarray}\label{growth}
		\sup _{i \in \mathbb{N}}|g(v_{1}) e_{i}| \leq \sup _{i \in \mathbb{N}}|g(v_{1}) e_{i}-g(0) e_{i}| +\sup _{i \in \mathbb{N}}|g(0) e_{i}|
		\leq  c_1 (1+|v_{1}|),
		\end{eqnarray}
		here $c_1:=\max\{L_{g},\sup _{i \in \mathbb{N}}|g(0) e_{i}|\}$ and by \cite[Lemma 7.1]{Rascanu2002},
		\begin{eqnarray}\label{linear}
		\sup _{i \in \mathbb{N}}|g(v_{1}) e_{i}-g(v_{2})e_{i}-g(u_1) e_{i}+g(u_2) e_{i}| &\leq& c_2 |v_{1}-v_{2}-u_1+u_2|\cr
		&&+c_2 |v_1-v_2|(|v_1-u_1|+|v_2-u_2|),
		\end{eqnarray}
		holds.}
\end{rem}

Taking a parameter $ 1-H<\alpha<\frac{1}{2}$, for
the measurable functions $h:[0,T] \rightarrow V$, let
\begin{eqnarray}\label{alphanorm}
\|h(t)\|_{\alpha}=|h(t)|+\int_{0}^{t} \frac{|h(t)-h(s)|}{(t-s)^{\alpha+1}} \mathrm{d}s.
\end{eqnarray}

Denote by $\mathcal{B}^{\alpha,2}(0, T ; V)$ the space of measurable functions $h:[0, T] \rightarrow V$ endowed with the norm $\|\cdot\|_{\alpha,T}$ defined by $$\|h\|_{\alpha,T}^2:=\sup_{t\in [0,T]}|h(t)|^2+\int_{0}^{T}\Big(\int_{0}^{t} \frac{|h(t)-h(s)|}{(t-s)^{\alpha+1}} \mathrm{d} s \Big)^2\mathrm{d}t<\infty.$$
\begin{rem}
{\rm Note that $\mathcal{B}^{\alpha,2}(0, T ; V)$ is continuously embedded in $W^{\alpha,1}(0,T:V)$.}
\end{rem}
\begin{defn}{\rm
		For $\alpha \in (1-H,\frac{1}{2})$, $V$-valued process $(u_t)_{ t\in[0,T]}$, is a solution of Eq. (\ref{spde-uniq}) in the mild sense if the following two conditions are satisfied:
		\begin{enumerate}
			\item $u \in \mathcal{B}^{\alpha,2}(0, T ; V)$ a.s.;
			\item $\{u_t\}$ is adapted to $\{\mathscr{F}_t\}$
			and satisfies the following integral equation:
			\begin{eqnarray}\label{mild1}
			u_t=S_t u_{0}+\int_{0}^{t} S_{t-s} f(u_s) \mathrm{d}s+\int_{0}^{t} S_{t-s} \sigma(u_s) \mathrm{d} W_s+\int_{0}^{t} S_{t-s} g(u_s) \mathrm{d} B^H_s.
			\end{eqnarray}
	\end{enumerate}}
\end{defn}

The following Lemma \ref{itofbm2} which will be proved in Appendix B
provides the basic estimates needed to prove the pathwise unique solutions of Eq. (\ref{spde-uniq}).
\begin{lem}\label{itofbm2}
	Taking $\alpha<\alpha'<1-\beta$ and $\alpha<\beta$, there exists a constant $C>0$, such that for any $0\leq s<t \leq T$, $u,v\in V$,
	\begin{eqnarray*}
		\mathcal{K}_1 (s, t) &:=&\Big|\int_{s}^{t} S_{t-r} g(u_r) \mathrm{d} B^H_r\Big|\cr
		&\leq& C \Lambda_{\alpha,B^H}^{0,t} \int_{s}^{t}\Big([(r-s)^{-\alpha}+(t-r)^{-\alpha}](1+|u_r|)+\int_{s}^{r}\frac{|u_r-u_q|}{(r-q)^{1+\alpha}}\mathrm{d}q\Big)\mathrm{d}r,
		\cr
		\mathcal{K}_2(0,s) &:=&\Big|\int_{0}^{s} (S_{t-r}-S_{s-r})g(u_r) \mathrm{d} B^H_r\Big| \cr
		&\leq&   C\Lambda_{\alpha,B^H}^{0,t} (t-s)^\beta\int_{0}^s  [(s-r)^{-\beta} r^{-\alpha}+(s-r)^{-\alpha-\beta}](1+|u_r|)\mathrm{d}r\cr
		&&+C\Lambda_{\alpha,B^H}^{0,t} (t-s)^\beta\int_{0}^s  (s-r)^{-\beta}\Big(\int_{0}^{r}\frac{|u_r-u_q|}{(r-q)^{1+\alpha}}\mathrm{d}q\Big)\mathrm{d}r,\cr
		\mathcal{K}_3(s,t) &:=&\Big|\int_{s}^{t} S_{t-r} (g(u_{r})-g(v_{r})) \mathrm{d} B^H_r\Big| \cr
		&\leq& C \Lambda_{\alpha,B^H}^{0,t}
		\int_{s}^t[(r-s)^{-\alpha}+(t-r)^{-\alpha}]|u_r-v_r|
		\mathrm{d}r \cr
		&&+C  \Lambda_{\alpha,B^H}^{0,t}
		\int_{s}^t|u_r-v_r|\Big(\int_{s}^{r}\frac{|u_r-u_q|+|v_r-v_q|}{(r-q)^{1+\alpha}}\mathrm{d}q\Big)\mathrm{d}r
		\cr
		&&+C  \Lambda_{\alpha,B^H}^{0,t}
		\int_{s}^t\Big(\int_{s}^{r}\frac{|u_r-v_r-u_q+v_q|}{(r-q)^{1+\alpha}}\mathrm{d}q\Big)\mathrm{d}r,
		\cr
		\mathcal{K}_4(0,s) &:=& \Big|\int_{0}^{s} (S_{t-r}-S_{s-r}) (g(u_{r})-g(v_{r})) \mathrm{d} B^H_r\Big|\cr
		&\leq& C  \Lambda_{\alpha,B^H}^{0,t} (t-s)^\beta \int_{0}^s [(s-r)^{-\beta} r^{-\alpha}+(s-r)^{-\alpha-\beta}]|u_r-v_r|\mathrm{d}r\cr
		&&+C   \Lambda_{\alpha,B^H}^{0,t} (t-s)^\beta \int_{0}^s (s-r)^{-\beta} |u_r-v_r| \Big(\int_{0}^{r}\frac{|u_r-u_q|+|v_r-v_q|}{(r-q)^{1+\alpha}}\mathrm{d}q\Big)\mathrm{d}r\cr
		&&+C   \Lambda_{\alpha,B^H}^{0,t}  (t-s)^\beta \int_{0}^s (s-r)^{-\beta}\Big(\int_{0}^{r}\frac{|u_r-v_r-u_q+v_q|}{(r-q)^{1+\alpha}}\mathrm{d}q\Big)\mathrm{d}r.
	\end{eqnarray*}
\end{lem}

\begin{thm}\label{mixuniq}
	Assume that $f,\sigma$ satisfy
	Lipschitz and linear growth conditions, and $g$ and $g^{\prime}$ satisfy (\ref{g1}) and (\ref{g2}), respectively. Then, for any initial value $u_{0} \in V_\beta,$ $\beta >\alpha$,  there exists a unique mild pathwise  solution to Eq. (\ref{spde-uniq}).
\end{thm}

Note that the unique solution for given $u_0$ is independent of $\alpha.$ The proof of Theorem \ref{mixuniq} will be divided into several logical steps.
\para{Step1: Construction of approximations.}
We recall the following auxiliary technical lemma from \cite{Mishura2012mixed}. The proof can be obtained by \cite[Lemma 2.1, Proposition 2.1]{Mishura2012mixed}, thus, we omit it.
\begin{lem}\label{smooth}
	Let $\varpi\in(0,1]$ and $\hbar:[0,T]\rightarrow \mathbb{R}$ be a $\varpi$-H\"older continuous function. Define for $\epsilon>0,$ $$\hbar^\epsilon(t)=\epsilon^{-1}\int_{0  \vee (t-\epsilon)}^{t}\hbar(s)ds.$$ Then, for $\alpha\in(1-\varpi,1)$, there exists a constant $C>0$ such that
	$$\|\hbar-\hbar^\epsilon\|_{\alpha,0,T}\leq C K_{\varpi}(\hbar)\epsilon^{\varpi+\alpha-1},$$
	where $K_{\varpi}(\hbar)=\sup_{0\leq s <t \leq T}\frac{|\hbar(t)-\hbar(s)|}{(t-s)^{\varpi}}$ is the $\varpi$-H\"older seminorm of $\hbar$.
\end{lem}

Fix $N \geq 1$, we define the following stopping time $\tau_N$,
\begin{eqnarray}\label{stoptime}
\tau_N :=\inf \{t\geq 0:\Lambda^{0,t}_{\alpha,B^H} \geq N \} \wedge T,
\end{eqnarray}
where $\Lambda_{\alpha,B^H}^{0,t}:=\sum_{i=1}^\infty\sqrt{\lambda_i} \Lambda_\alpha^{0,t}(\beta_i^H)$.

Put $\beta_i^{H, N}(t)=\beta^H_{i}(t\wedge \tau_N), t\in [0,T],i \in \mathbb{N}$ and taking $\epsilon=\frac{1}{n}$ in Lemma \ref{smooth},  then, for each $n,i \in \mathbb{N}$,  define an approximation of $\beta_i^{H, N, n}$ by
\begin{eqnarray}\label{betaHn}
\beta_i^{H, N, n}(t)=n\int_{(t-\frac{1}{n})\vee 0}^{t}\beta_i^{H,N}(s)\mathrm{d}s.
\end{eqnarray}

Similarly, denote $\Lambda_{\alpha,B^{H,N, n}}^{0,T} :=\sum_{i=1}^\infty\sqrt{\lambda_i} \Lambda_\alpha^{0,T}(\beta_i^{H,N,n})$ which will be used in next step.
\begin{lem}\label{betaHn1}
	For any $ i\in \mathbb{N}$ and $N$, we have
	\begin{eqnarray*}
		\|\beta_{i}^{H,N}\|_{\alpha,0,T} \leq C \|\beta^H_i\|_{\alpha,0,\tau_N},
	\end{eqnarray*}
	almost surely, where $C$ is a constant which is independent of $N$ and $i$.
\end{lem}

Note that the above lemma will be proven in Appendix B.

To proceed, by Lemma \ref{smooth} and Lemma \ref{betaHn1}, for any $i\in \mathbb{N}$ and $\varpi \in(1-\alpha, H)$, we have
\begin{eqnarray}\label{betaNn-N}
\|\beta_{i}^{H,N,n}-\beta_{i}^{H,N}\|_{\alpha,0,T} &\leq& CK_{\varpi}(\beta_{i}^{H,N})(\frac{1}{n})^{\varpi+\alpha-1} \cr
&\leq& C\|\beta_{i}^{H,N}\|_{1-\varpi,0,T}(\frac{1}{n})^{\varpi+\alpha-1}\cr
&\leq& C\|\beta_{i}^{H}\|_{1-\varpi,0,\tau_N}(\frac{1}{n})^{\varpi+\alpha-1},
\end{eqnarray}
almost surely. Moreover, since $\varpi \in(1-\alpha, H)$,  by Lemma \ref{betaHn1} and (\ref{betaNn-N}), we have
\begin{eqnarray}\label{betaNn-Nn}
\Lambda_{\alpha,B^{H,N, n}}^{0,T}
&\leq& \frac{ \sum_{i=1}^\infty\sqrt{\lambda_i}\big(\|\beta_{i}^{H,N,n}-\beta_{i}^{H,N}\|_{\alpha,0,T}+\|\beta_{i}^{H,N}\|_{\alpha,0,T}\big)}{\Gamma(1-\alpha)\Gamma(\alpha)}\cr
&\leq& \frac{C  \sum_{i=1}^\infty\sqrt{\lambda_i}(\|\beta_{i}^{H}\|_{\alpha,0,\tau_N}+\|\beta_{i}^{H}\|_{1-\varpi,0,\tau_N})}{\Gamma(1-\alpha)\Gamma(\alpha)}\cr
&\leq& CN.
\end{eqnarray}

Now,
let us consider
\begin{eqnarray}\label{map1}
u_t^{N,n}&=& S_t u_{0}+\int_{0}^{t} S_{t-s} f(u^{N,n}_s) \mathrm{d} s+\int_{0}^{t} S_{t-s} \sigma(u^{N,n}_s) \mathrm{d} W_s\cr
&&+\sum_{i=1}^{\infty} \sqrt{\lambda_{i}}\int_{0}^{t} S_{t-s} g(u^{N,n}_s)e_i  \mathrm{d} \beta_{i}^{H,N,n}(s),
\end{eqnarray}
or equivalently
\begin{eqnarray}\label{map2}
u_t^{N,n}=S_t u_{0}+\int_{0}^{t} S_{t-s} f^{N,n}(u^{N,n}_s) \mathrm{d} s+\int_{0}^{t} S_{t-s} \sigma(u^{N,n}_s) \mathrm{d}W_s,
\end{eqnarray}
where $f^{N,n}(u):=f(u)+g(u) \sum_{i=1}^{\infty} \sqrt{\lambda_{i}}e_i  \frac{\mathrm{d}}{\mathrm{d}s} \beta_{i}^{H,N,n}(s)$ is a random drift.

Such equations were studied in \cite[Section 7]{da2014stochastic}. To proceed, it is easy to obtain
\begin{eqnarray}\label{lambeta}
\Big|\sum_{i=1}^{\infty} \sqrt{\lambda_{i}}\frac{\mathrm{d}}{\mathrm{d}s} \beta_{i}^{H,N,n}(s)\Big|&\leq& \sum_{i=1}^{\infty} \sqrt{\lambda_{i}} \Big|\frac{\mathrm{d}}{\mathrm{d}s} \beta_{i}^{H,N,n}(s)\Big|\cr
&\leq & n \sum_{i=1}^{\infty} \sqrt{\lambda_{i}} \Big| \beta_{i}^{H,N}(s)-\beta_{i}^{H,N}((s-\frac{1}{n})\vee 0)\Big|\cr
&\leq & n^\alpha \sum_{i=1}^{\infty} \sqrt{\lambda_{i}} \Lambda^{0,s}_{\alpha} (\beta_{i}^{H,N})\cr
&\leq & C_{n,N},
\end{eqnarray}
where $C_{n,N}$ is a constant dependent on $n$ and $N$. Thus, by (\ref{lambeta}), we obtain the function $f^{N,n}$ satisfies Lipschitz and growth conditions. Then, by \cite[Theorem 7.4]{da2014stochastic}, there exists a unique mild solution $u^{N,n}$ to Eq.(\ref{map2}).

\para{Step 2: Convergence of approximations.}
To obtain the convergence of approximations, we give the following two key lemmas which will be proved in Appendix B.
\begin{lem}\label{lemboundN}
	Under the assumptions of Theorem \ref{mixuniq}, there exists a constant $C_{N}$, such that
	\begin{eqnarray*}
		\mathbb{E}[\|u^{N,n}\|^2_{\alpha, T}] \leq C_{N}.
	\end{eqnarray*}
\end{lem}

\begin{lem}\label{lemboundNm}
	Under the assumptions of Theorem \ref{mixuniq}. Then,
	there exists a constant $C_{N,R}$, such that
	\begin{eqnarray}\label{bounduNm}
	\mathbb{E}\Big[\|u^{N,n}-u^{N,m}\|^2_{\alpha,T}\mathbf{1}_{D_{T}^{N,R}}\Big]  \leq C_{N,R}\mathbb{E}\Big[\Big(\sum_{i=1}^{\infty}\sqrt{\lambda_{i}}\|\beta_{i}^{H,N,n}-\beta_{i}^{H,N,m }\|_{\alpha,0,T}\Big)^2\Big],
	\end{eqnarray}
	where $D_{T}^{N,R}:=\{\|u^{N,n}\|_{\alpha,T} \leq R, \|u^{N,m}\|_{\alpha,T} \leq R\}$.
\end{lem}

For fixed $N\geq 1$, we will show that the sequence $\{u^{N,n}, n\geq 1\}$ is Cauchy sequence in the norm $\|\cdot\|_{\alpha,T}$. For all $\varepsilon>0$, $R\geq 1$, we have
\begin{eqnarray*}
\mathbb{P}\big(\|u^{N,n}-u^{N,m}\|_{\alpha,T}>\varepsilon\big) &\leq&
\mathbb{P}\big( \|u^{N,n}\|_{\alpha,T} > R~{\rm or}~ \|u^{N,m}\|_{\alpha,T} > R \big) \cr
&&+\mathbb{P}\big(\|u^{N,n}-u^{N,m}\|_{\alpha,T}>\varepsilon, \|u^{N,n}\|_{\alpha,T} \leq R, \|u^{N,m}\|_{\alpha,T} \leq R \big).
\end{eqnarray*}

Since $\varpi \in(1-\alpha, H)$, by (\ref{betaNn-N}), we have
\begin{eqnarray}\label{boundbetan}
\sup_{i\in \mathbb{N}}\mathbb{E}\big[\|\beta_{i}^{H,N,n}-\beta_{i}^{H,N}\|^2_{\alpha,0,T}\big] &\leq& C(\frac{1}{n})^{2(\varpi+\alpha-1)} \sup_{i\in \mathbb{N}} \mathbb{E}\big[\|\beta_{i}^{H}\|^2_{1-\varpi,0,T}\big]\cr
&\leq& C(\frac{1}{n})^{2(\varpi+\alpha-1)}.
\end{eqnarray}
where $C$ is a constant which is independent of $i$ and $N$.

Then, by Cauchy-Schwarz's inequality, we have
\begin{eqnarray*}
	\mathbb{E}\Big[\Big(\sum_{i=1}^{\infty}\sqrt{\lambda_{i}}\|\beta_{i}^{H,N,n}-\beta_{i}^{H,N}\|_{\alpha,0,T}\Big)^2\Big]
	&=&\mathbb{E}\Big[\Big(\sum_{i=1}^{\infty}\sqrt[4]{\lambda_{i}}\sqrt[4]{\lambda_{i}}\|\beta_{i}^{H,N,n}-\beta_{i}^{H,N}\|_{\alpha,0,T}\Big)^2\Big]\cr
	&\leq& \sum_{i=1}^{\infty}\sqrt{\lambda_{i}}\Big(\sum_{i=1}^{\infty}\sqrt{\lambda_{i}}\mathbb{E}\big[\|\beta_{i}^{H,N,n}-\beta_{i}^{H,N}\|^2_{\alpha,0,T}\big]\Big).
\end{eqnarray*}

Thus, the condition $\sum_{i=1}^{n}\sqrt{\lambda_{i}}<\infty$ and (\ref{boundbetan}) yield that
\begin{eqnarray}
\mathbb{E}\Big[\Big(\sum_{i=1}^{\infty}\sqrt{\lambda_{i}}\|\beta_{i}^{H,N,n}-\beta_{i}^{H,N}\|_{\alpha,0,T}\Big)^2\Big] \leq  C(\frac{1}{n})^{2(\varpi+\alpha-1)}\Big(\sum_{i=1}^{\infty}\sqrt{\lambda_{i}}\Big)^2\rightarrow 0 \qquad ~{\rm as}~ n \rightarrow \infty.
\end{eqnarray}

Due to Lemma \ref{lemboundNm} and Markov's inequality, we see that for all $\varepsilon>0, R \geq 1$, $$\mathbb{P}\big(\|u^{N,n}-u^{N,m}\|_{\alpha,T}>\varepsilon, \|u^{N,n}\|_{\alpha,T} \leq R, \|u^{N,m}\|_{\alpha,T} \leq R \big) \rightarrow 0\qquad ~{\rm as}~n \rightarrow \infty,$$
and $$\limsup_{n,m \rightarrow \infty}\mathbb{P}\big( \|u^{N,n}\|_{\alpha,T} > R~{\rm or}~\|u^{N,m}\|_{\alpha,T} > R \big) \leq 2 \sup_{n\in\mathbb{N}} \mathbb{P}\big( \|u^{N,n}\|_{\alpha,T} > R\big).$$
Next, Lemma \ref{lemboundN} and Markov's inequality imply that $$\sup_{n\in\mathbb{N}} \mathbb{P}\big( \|u^{N,n}\|_{\alpha,T} > R\big)\rightarrow 0 \qquad ~{\rm as}~ R \rightarrow \infty   .$$
Thus, $$\|u^{N,n}-u^{N,m}\|_{\alpha,T}\rightarrow0\qquad ~{\rm as}~n,m \rightarrow \infty,$$ in probability. Then there exists a random process $u^N$ such that $$\|u^{N,n}-u^{N}\|_{\alpha,T}\rightarrow0,\qquad ~{\rm as}~ n \rightarrow \infty,$$ in probability. Denoting an almost surely convergent subsequence by the same symbol, we have
$$\|u^{N,n}-u^{N}\|_{\alpha,T}\rightarrow 0\qquad ~{\rm as}~n \rightarrow \infty,$$ almost surely.

\para{Step 3: The limit provides a solution.}
Since $\|u^{N,n}-u^{N}\|_{\alpha,T}\rightarrow0, ~{\rm as}~n \rightarrow \infty,$ a.s., we easily obtain
$$\int_{0}^{t} S_{t-s} f(u^{N,n}_{s})\mathrm{d}s \rightarrow \int_{0}^{t} S_{t-s} f(u^{N}_{s}) \mathrm{d}s\qquad ~{\rm as}~n \rightarrow \infty,$$
almost surely.

Similar to the cases of the proof in Lemma \ref{lemboundNm}, denoting $\mathbf{1}_{t}:=\mathbf{1}_{\{\|u^{N,n}\|_{\alpha,t}\leq R, \|u^{N}\|_{\alpha,t} \leq R\}}$, we have
\begin{eqnarray*}
	&& \Big|\sum_{i=1}^{\infty}\sqrt{\lambda_{i}}\int_{0}^{t} S_{t-s}  (g(u^{N,n}_{s})-g(u^{N}_{s}))e_i \mathrm{d} \beta_{i}^{H,N,n}(s)\mathbf{1}_{t}\Big|^2\cr
	&\leq& C \Big(\sum_{i=1}^{\infty}\sqrt{\lambda_{i}}\|\beta_{i}^{H,N,n}\|_{\alpha,0,T}
	\int_{0}^t[r^{-\alpha}+(t-r)^{-\alpha}]|u^{N,n}_r-u^{N}_r|
	\mathrm{d}r\mathbf{1}_{t}\Big)^2 \cr
	&&+C \Big(\sum_{i=1}^{\infty}\sqrt{\lambda_{i}}\|\beta_{i}^{H,N,n}\|_{\alpha,0,T}
	\int_{0}^t|u^{N,n}_r-u^{N}_r|\big(\int_{0}^{r}\frac{|u^{N,n}_r-u^{N,n}_q|}{(r-q)^{1+\alpha}}\mathrm{d}q\big)\mathrm{d}r\mathbf{1}_{t}\Big)^2
	\cr
	&&+C \Big(\sum_{i=1}^{\infty}\sqrt{\lambda_{i}}\|\beta_{i}^{H,N,n}\|_{\alpha,0,T}
	\int_{0}^t|u^{N,n}_r-u^{N}_r|\Big(\int_{0}^{r}\frac{|u^{N}_r-u^{N}_q|}{(r-q)^{1+\alpha}}\mathrm{d}q\Big)\mathrm{d}r\mathbf{1}_{t}\Big)^2
	\cr
	&&+C \Big(\sum_{i=1}^{\infty}\sqrt{\lambda_{i}}\|\beta_{i}^{H,N,n}\|_{\alpha,0,T}
	\int_{0}^t\Big(\int_{0}^{r}\frac{|u^{N,n}_r-u^{N}_r-u^{N,n}_q+u^{N}_q|}{(r-q)^{1+\alpha}}\mathrm{d}q\Big)\mathrm{d}r\mathbf{1}_{t}\Big)^2 \cr
	&\leq& C_{N,R}
	\int_{0}^t\|u^{N,n}-u^{N}\|^2_{\alpha,r}\mathbf{1}_{r}\mathrm{d}r+C_{N} \|u^{N,n}-u^{N}\|^2_{\alpha,t}\mathbf{1}_{t}.
\end{eqnarray*}
Thus, it is easy to obtain
$$\sum_{i=1}^{\infty}\sqrt{\lambda_{i}}\int_{0}^{t} S_{t-s} g(u^{N,n}_{s}) e_i \mathrm{d} \beta_{i}^{H,N,n}(s) \rightarrow \sum_{i=1}^{\infty}\sqrt{\lambda_{i}}\int_{0}^{t} S_{t-s}  g(u^{N}_{s}) e_i \mathrm{d} \beta_{i}^{H,N,n}(s)\qquad ~{\rm as}~n \rightarrow \infty,$$
almost surely. Finally,  we will obtain
\begin{eqnarray*}
	&& \Big|\sum_{i=1}^{\infty}\sqrt{\lambda_{i}}\int_{0}^{t} S_{t-s}  g(u^{N}_{s})e_i d (\beta_{i}^{H,N,n}(s)-\beta_{i}^{H,N}(s))\mathbf{1}_{t}\Big|\cr
	&\leq& C \Big(\sum_{i=1}^{\infty}\sqrt{\lambda_{i}}\|\beta_{i}^{H,N,n}-\beta_{i}^{H,N}\|_{\alpha,0,T}\int_{0}^{t}(r^{-\alpha}+(t-r)^{-\alpha})(1+\|u^{N}\|_{\alpha,r})\mathrm{d}r\Big)\mathbf{1}_{t}\cr
	&&+C \Big(\sum_{i=1}^{\infty}\sqrt{\lambda_{i}}\|\beta_{i}^{H,N,n}-\beta_{i}^{H,N}\|_{\alpha,0,T}\int_{0}^t\Big( \int_{0}^{r}\frac{|u^{N}_r-u^{N}_q|}{(r-q)^{1+\alpha}}\mathrm{d}q\Big)\mathrm{d}r\Big)\mathbf{1}_{t}\cr
	&\leq& C_{R} \sum_{i=1}^{\infty}\sqrt{\lambda_{i}}\|\beta_{i}^{H,N,n}-\beta_{i}^{H,N}\|_{\alpha,0,T}.
\end{eqnarray*}

Since $\|\beta_{i}^{H,N,n}-\beta_{i}^{H,N}\|_{\alpha,0,T} \rightarrow 0, n \rightarrow \infty$, a.s., then, we have
$$\Big|\sum_{i=1}^{\infty}\sqrt{\lambda_{i}}\int_{0}^{t} S_{t-s}  g(u^{N}_{s})e_i \mathrm{d} (\beta_{i}^{H,N,n}(s)-\beta_{i}^{H,N}(s))\Big|\rightarrow 0\qquad ~{\rm as}~n \rightarrow \infty,$$
almost surely.

Next, we have
\begin{eqnarray*}
	\mathbb{E} \Big[ \Big|\int_{0}^{t} S_{t-s} (\sigma(u^{N,n}_{s})-\sigma(u^{N}_{s})) \mathrm{d}W_s\Big|^2\mathbf{1}_{t}\Big]
	&\leq&  \int_{0}^{t} \mathbb{E} \big[|\sigma(u^{N,n}_{s})-\sigma(u^{N}_{s}))|_{L_2(V)}^2\mathbf{1}_{s}\big] \mathrm{d}s  \cr
	&\leq&  \int_{0}^{t} \mathbb{E} \big[\|u^{N,n}-u^{N}\|_{\alpha,s}^2\mathbf{1}_{s}\big] \mathrm{d}s.
\end{eqnarray*}

Thus, by Lemma \ref{lemboundNm}, it is easy to obtain
$$\mathbb{E} \Big[ \Big|\int_{0}^{t} S_{t-s} (\sigma(u^{N,n}_{s})-\sigma(u^{N}_{s})) \mathrm{d}W_s\Big|^2\mathbf{1}_{t}\Big] \rightarrow 0\qquad ~{\rm as}~ n\rightarrow \infty, $$
and, consequently
$$\int_{0}^{t} S_{t-s} \sigma(u^{N,n}_{s})\mathrm{d}W_s-\int_{0}^{t} S_{t-s}\sigma(u^{N}_{s})\mathrm{d}W_s \rightarrow 0\qquad ~{\rm as}~ n\rightarrow \infty,$$
in probability. Since $\|u^{N,n}-u^{N}\|_{\alpha,T}\rightarrow0, n \rightarrow \infty,$ a.s., we have the convergence of the integrals in probability on $\{\|u^{N}\|_{\alpha,T}\leq R\}$, where $R\geq 1$ is arbitrary, therefore the convergence holds on $\Omega$. This means that $u^N$ is a solution to
\begin{eqnarray}\label{map3}
u_t^{N}= S_t u_{0}+\int_{0}^{t} S_{t-s} f(u^{N}_s) d s+\int_{0}^{t} S_{t-s} \sigma(u^{N}_s) \mathrm{d}W_s+\int_{0}^{t} S_{t-s} g(u^{N}_s)\mathrm{d} B_s^{H,N}.
\end{eqnarray}
\para{Step 4: Letting $N\rightarrow \infty$ and uniqueness.}
It follows from Lemma \ref{lemboundNm} that the processes $u^N$ and $u^M$ with $M\geq N$ coincide almost surely on the set $A_{N,T}=\{\Lambda^{0,T}_{\alpha,B^H} \leq N\}$. Hence, there exists a process $u$ such that $u^N=u$, a.s. on $A_{N,T}$ for each $N\geq 1$, hence, almost surely. Finally, the pathwise uniqueness follows in similar way. Thus, the proof is finished.
\qed

\section{Fast-slow SPDEs Driven by fBm and Bm}\label{se-4}
Throughout this section, we assume that the following conditions are fulfilled. We assume that
\begin{enumerate}
	\item[(A1)] The coefficients $b(x,y): V\times V \rightarrow V,F(x,y): V \times V \rightarrow V,G(x,y): V \times V \rightarrow L_2(V)$ of Eq. (\ref{spde-fs}) are globally Lipschitz continuous in $x,y$, i.e., there exist two positive constants $C_{1},C_{2}$, such that
	\begin{eqnarray*}
		|b(x_{1},y_{1})-b(x_{2},y_{2})|^{2}
		&\leq& C_{1}(|x_{1}-x_{2}|^{2}+|y_{1}-y_{2}|^{2}),\cr
		|F(x_{1},y_{1})-F(x_{2},y_{2})|^{2}+|G(x_{1},y_{1})-G(x_{2},y_{2})|_{L_2(V)}^{2} &\leq& C_{2}(|x_{1}-x_{2}|^{2}+|y_{1}-y_{2}|^{2}),
	\end{eqnarray*}
	for all $x_{1},x_{2},y_{1},y_{2}\in V$.
	\item [(A2)]
	The coefficients $b(x,y),F(x,y),G(x,y)$ of Eq. (\ref{spde-fs}) satisfy linear growth conditions,  i.e., there exist two positive constants $C_{3}, C_{4}$ such that
	\begin{eqnarray*}
		|F(x,y)|^{2}+|G(x,y)|_{L_2(V)}^{2} &\leq& C_{3}(1+|x|^{2}+|y|^{2})\cr
		|b(x,y)|^{2} &\leq& C_{4}(1+|x|^{2}+|y|^{2}),
	\end{eqnarray*}
	for all $x,y\in V$.
	\item[(A3)] $g: V \rightarrow L(V)$ and $g^{\prime} : V \rightarrow L(V, L(V))$ are Lipschitz continuous in the senses of Eq. (\ref{g1}) and Eq. (\ref{g2}).
	\item[(A4)] There exist constants $\beta_{1},C_5>0$ and $\beta_{2},\beta_{3}\in \mathbb{R}$ which are independent of $(x,y_{1},y_{2})$, such that
	\begin{eqnarray*}
		\langle y_{1}, F(x,y_{1})\rangle & \leq& -\beta_{1} |y_{1}|^{2}+\beta_{2},\cr
		\langle y_{1}-y_{2}, F(x_1,y_{1})-F(x_2,y_{2}) \rangle &\leq & \beta_{3} |y_{1}-y_{2}|^{2}+C_5|x_{1}-x_{2}|^{2},
	\end{eqnarray*}
	for all $x,y_{1},y_{2}\in V$.
	\item[(A5)] $\eta:=2\bar{\lambda}_1-2\beta_{3}-C_{2}>0,\kappa:=2\bar{\lambda}_1+2\beta_{1}-C_{3}>0, \theta:=\bar{\lambda}_1^2-C_3>0$, where $\bar{\lambda}_1$ is the first eigenvalue of $-A$, $C_i, i=2,3$ and $\beta_i, i=1,2,3$ were given in (A1), (A2) and (A4).
\end{enumerate}

\begin{rem}\label{remdissipative}
	{\rm Assumptions (A4) and (A5) are known as the strong dissipative conditions that imply the existence of a unique invariant measure and moreover, it has
		exponentially mixing property for the Markov semigroup associated to the fast variable.
	}
\end{rem}

Through Theorem \ref{mixuniq} and a similar argument as in the proof of \cite[Theorem 2.2]{Liu2019}, it is easy to prove that Eq. (\ref{spde-fs}) has a unique mild pathwise solution. Here, we omit the proof.
\begin{lem}
	Suppose that conditions {\rm (A1)-(A3)} are satisfied. Then, for any initial values $X_{0},Y_0\in V_\beta,\beta>\alpha$, Eq. (\ref{spde-fs}) has a unique mild pathwise  solution $(X^\varepsilon_t,Y^\varepsilon_t)$, i.e.,
	\begin{align}\label{mild2}
	X^{\varepsilon}_t
	&=S_t X^{\varepsilon}_{0}+\int_0^tS_{t-s}b(X^{\varepsilon}_s,Y^{\varepsilon}_s)\mathrm{d}s+\int_0^tS_{t-s}g(X^{\varepsilon}_s)\mathrm{d}B^H_s, \qquad X^{\varepsilon}_0= X_0, \\
	\label{mild3} Y^{\varepsilon}_t& =S_\frac{t}{\varepsilon}Y^{\varepsilon}_{0}+\frac{1}{\varepsilon}\int_0^tS_{\frac{t-s}{\varepsilon}}F(X^{\varepsilon}_s,Y^{\varepsilon}_s)\mathrm{d}s+\frac{1}{\sqrt{\varepsilon}}\int_0^tS_{\frac{t-s}{\varepsilon}}G(X^{\varepsilon}_s,Y^{\varepsilon}_s)\mathrm{d}W_s, \qquad  Y^{\varepsilon}_0= Y_0.
	\end{align}
\end{lem}

Denote by $\bar{X}_t, t\in[0,T]$, the solution of the following SPDEs driven by fBm,
\begin{eqnarray}\label{ave-eq}
\begin{split}
\mathrm{d} \bar{X}_t=(A\bar{X}_t+\bar{b}(\bar{X}_t))\mathrm{d}t+g(\bar{X}_t)\mathrm{d}B^H_t,  \qquad \bar{X}_0= X_0,
\end{split}
\end{eqnarray}
where
\begin{eqnarray}\label{invar-mea}
\bar{b}(x)=\int_{V}b(x,z)\mu^x(\mathrm{d} z),\qquad x \in V,
\end{eqnarray}
and $\mu^x$ is the unique invariant measure on $V$ of the transition
semigroups for the following frozen equation:
\begin{eqnarray}\label{invariant}
\mathrm{d} Y^{x}_t =(AY^{x}_t+F(x,Y^{x}_t))\mathrm{d}t+G(x,Y^{x}_t)\mathrm{d}{W}_t, \qquad Y^{x}_0=y \in V.
\end{eqnarray}

According to the definition of $\bar b$ (\ref{invar-mea}) and conditions {\rm (A1)-(A4)}, it is easy to prove $\bar b$ also satisfies the Lipschitz and growth conditions. Then, we have the following lemma.
\begin{lem}
	Suppose that conditions {\rm (A1)-(A5)} are satisfied. Then, for any initial value $X_{0} \in V_\beta,\beta>\alpha$, Eq. (\ref{ave-eq}) has a unique mild pathwise  solution.
\end{lem}
\para{Proof:} For any $x_1, x_2, x \in V$ and any initial value $y\in V$, by (\ref{yxz}) and (\ref{measure2}) in Appendix A, we have
\begin{eqnarray*}
|\bar{b}(x_1)-\bar{b}(x_2)|
&\leq&\bigg|\int_{V} b(x_1, z)\mu^{x_1}(\mathrm{d} z)-\mathbb{E}[b(x_1,Y^{x_1,y}_s)]\bigg|^2\cr
&&+\bigg|\int_{V} b(x_2,z) \mu^{x_2}(\mathrm{d} z)-\mathbb{E}[b(x_2,Y^{x_2,y}_s)]\bigg|^2\cr
&&+\big|\mathbb{E}[b(x_1,Y^{x_1,y}_s)-b(x_2,Y^{x_2,y}_s)]\big|^2\cr
&\leq&Ce^{-\eta s}(1+|x_1|^2+|x_2|^2+|y|^2)+C|x_1-x_2|^2,
\end{eqnarray*}
Let $s\rightarrow \infty$, then we obtain that $\bar{b}_1$ is Lipschitz continuous in $x$, and
\begin{eqnarray*}
|\bar{b}(x)|^2\leq \bigg(\int_{V}|b(x, z)|\mu^{x}(d z)\bigg)^2\leq C(1+|x|^2).
\end{eqnarray*}
So, $\bar{b}_1$ satisfies the growth condition. Thus, according to Theorem \ref{mixuniq},  (\ref{ave-eq}) has a unique strong solution. \qed

From now on, we assume $\beta\in(\frac{1}{2},1-\alpha)$ and to present our main averaging result, we need to impose another condition.
 \begin{enumerate}
	\item[(B1)] $\sup_{x,y \in V}(|b(x,y)|+|G(x,y)|) <\infty$.
\end{enumerate}

\begin{thm}\label{mainthm}
	Suppose that conditions {\rm (A1)-(A5)} and {\rm (B1)} hold. Then, we have
	$$\lim\limits_{\varepsilon \rightarrow 0}\mathbb{E}\big[\|X^{\varepsilon}-\bar{X}\|^2_{\alpha,T}\big]=0.$$
\end{thm}

\begin{rem}
	{\rm To obtain the strong convergence, it is known that the diffusion coefficient $g$ in (\ref{spde-fs}) should not depend on the fast variable $Y^\varepsilon$ (see e.g. \cite{Givon2007}). }
\end{rem}

To prove Theorem \ref{mainthm},
firstly, following the discretization techniques inspired by Khasminskii in \cite{khas1966limit}, we introduce an auxiliary process $(\hat{X}_{t}^{\varepsilon},\hat{Y}_{t}^{\varepsilon})$ and divide $[0,T]$ into intervals of size $\delta$, where $\delta \in(0,1)$ is a fixed number depending on $\varepsilon$ and $\delta > \varepsilon $, which
will be chosen later.  Then,
we construct auxiliary processes $\hat{Y}^{\varepsilon}$ and $\hat{X}^{\varepsilon}$, by
\begin{align}
\label{yj} \hat{Y}^{\varepsilon}_t&=S_{\frac{t}{\varepsilon}}Y_{0}+\frac{1}{\varepsilon}\int_0^tS_{\frac{t-s}{\varepsilon}}F(X^{\varepsilon}_{s(\delta)},\hat{Y}^{\varepsilon}_s)\mathrm{d}s+\frac{1}{\sqrt{\varepsilon}}\int_0^tS_{\frac{t-s}{\varepsilon}}G(X^{\varepsilon}_{s(\delta)},\hat{Y}^{\varepsilon}_s)\mathrm{d}W_s,\\
\label{xj} \hat{X}^{\varepsilon}_t&=S_{t}X_{0}+\int_{0}^{t}S_{t-s}b(X^{\varepsilon}_{s(\delta)},
\hat{Y}^{\varepsilon}_s)\mathrm{d}s+\int_{0}^{t}S_{t-s}g(X^{\varepsilon}_s)\mathrm{d}B^H_s,
\end{align}
where $s(\delta)=\lfloor s / \delta\rfloor \delta$  is the nearest breakpoint preceding $s$.  For $t\in [k\delta, \min\{(k+1)\delta,T\}]$, we assume the fast component $\hat{Y}^{\varepsilon}_{k\delta}$
is reset to equal ${Y}^{\varepsilon}_{k\delta}$ at each breakpoint $k\delta$.
To proceed, we can derive uniform bounds $\|X^{\varepsilon}-\hat{X}^{\varepsilon}\|_{\alpha,T}^2 $.
Next, based on the ergodic property of the frozen equation, we obtain
appropriate control of $\|\hat{X}^{\varepsilon}-\bar{X}\|_{\alpha,T}^2 $.
Finally, we can estimate $\|{X}^{\varepsilon}-\bar{X}\|_{\alpha,T}^2 $.

\subsection{A Priori Estimate} Estimates of the auxiliary provess $(X^{\varepsilon}_t,Y^{\varepsilon}_t)$ wil be given in this subsection.
\begin{lem}\label{xbound}
	Suppose that conditions {\rm (A1)-(A3)} and {\rm (B1)} are satisfied. Then, for any $p\geq 2$, there exists a constant $C_p $ which is independent of $\varepsilon$ such that
	$$\mathbb{E}\Big[\sup_{t\in[0,T]}\| X^{\varepsilon}_t\|_{\alpha}^p\Big] \leq C_p,$$
	where $\|\cdot\|_{\alpha}$ was defined in (\ref{alphanorm}).
\end{lem}
\para{Proof:} For shortness, denote, $\Lambda:=\Lambda_{\alpha,B^H}^{0,T} \vee 1,$ and for $\rho \geq 1,$ let
\begin{eqnarray*}
	\|f\|_{\rho,T} &:=&\sup _{t \in[0, T]}e^{-\rho t}|f(t)|,\cr
	\|f\|_{1,\rho,T} &:=&\sup _{t \in[0, T]}e^{-\rho t} \int_{0}^{t} \frac{|f(t)-f(r)|}{(t-r)^{\alpha+1}} \mathrm{d}r.\end{eqnarray*}

Using techniques similar to those used in \cite[Lemma 4.1 ]{Shevchenko2014mixed}, we start by estimating $\|X^\varepsilon\|_{\rho,T}$. By the similar step as for the terms $\mathcal{M}_{21}$, $\mathcal{M}_{42}$ in Appendix B and using Lemma \ref{inq-rho}, Lemma \ref{itofbm2}, (A1)-(A3) and (B1), we have
\begin{eqnarray}\label{in1}
\|X^\varepsilon\|_{\rho,T}&=&\sup _{t \in[0, T]}e^{-\rho t}\Big|S_t X_0+\int_{0}^{t} S_{t-r}b(X_{r}^{\varepsilon}, Y_{r}^{\varepsilon}) \mathrm{d} r+ \int_{0}^{t} S_{t-r} g(X^\varepsilon_r) \mathrm{d} B^{H}_r\Big|\cr
&\leq& C \Big(1+|X_0|_{{\beta}}+ \Lambda_{\alpha,B^H}^{0,T}\sup _{t \in[0, T]}\int_{0}^{t} e^{-\rho (t-r)}[(r^{-\alpha}+(t-r)^{-\alpha}) \|X^\varepsilon\|_{\rho,T} +\|X^\varepsilon\|_{1,\rho,T}]\mathrm{d}r \Big)\cr
&\leq& K \Lambda\left(1+\rho^{\alpha-1}\|X^\varepsilon\|_{\rho, T}+\rho^{-1}\|X^\varepsilon\|_{1, \rho, T}\right),
\end{eqnarray}
with some constant $K$ (which is dependent on $|X_0|_{{\beta}}$ and can be assumed to be greater than 1 without loss of generality).

For $\|X^{\varepsilon}\|_{1,\rho,T}$, by the similar step as for the terms $\mathcal{M}_{41}$ in Appendix B and using Fubini's theorem, Lemma \ref{inq-rho}, Lemma \ref{itofbm2}, (A1)-(A3) and (B1), we have
\begin{eqnarray}\label{in2}
\|X^{\varepsilon}\|_{1,\rho,T} &\leq& C \sup_{t \in[0, T]} \int_{0}^t e^{-\rho t} \frac{|(S_t-S_s)X_{0}|}{(t-s)^{1+\alpha}}\mathrm{d}s\cr
&&
+C\sup _{t \in[0, T]} \int_{0}^t e^{-\rho t} \frac{\big|\int_{s}^{t}S_{t-r} b(X^{\varepsilon}_r, \hat Y^{\varepsilon}_r) \mathrm{d}r\big|}{(t-s)^{1+\alpha}}\mathrm{d}s\cr
&&+C\sup _{t \in[0, T]} \int_{0}^t e^{-\rho t} \frac{\big|\int_{0}^{s}(S_{t-r} -S_{s-r}) b(X^{\varepsilon}_r, \hat Y^{\varepsilon}_r)\mathrm{d} r\big|}{(t-s)^{1+\alpha}}\mathrm{d}s\cr
&&+C \sup _{t \in[0, T]}  \int_{0}^t e^{-\rho t} \frac{\big|\int_{s}^{t} S_{t-r} g(X^{\varepsilon}_r) \mathrm{d} B^H_r  \big|}{(t-s)^{1+\alpha}} \mathrm{d}s \cr
&&+C\sup _{t \in[0, T]}  \int_{0}^t e^{-\rho t} \frac{\big|\int_{0}^{s}(S_{t-r}-S_{s-r}) g(X^{\varepsilon}_r) \mathrm{d}B^H_r  \big|}{(t-s)^{1+\alpha}} \mathrm{d}s\cr
&\leq&  C(1+|X_0|_{{\beta}})\cr
&&+C\Lambda_{\alpha,B^H}^{0,T} \sup _{t \in[0, T]}  \int_{0}^t e^{-\rho t} \int_{s}^{t} \frac{(r-s)^{-\alpha}+(t-r)^{-\alpha}}{(t-s)^{1+\alpha}}(1+|X^{\varepsilon}_r|)\mathrm{d} r \mathrm{d}s\cr
&&+C\Lambda_{\alpha,B^H}^{0,T} \sup _{t \in[0, T]}  \int_{0}^t e^{-\rho t} \int_{s}^{t}(t-s)^{-1-\alpha}\Big(\int_{s}^{r} \frac{|X^{\varepsilon}_r-X^{\varepsilon}_q|}{(r-q)^{1+\alpha}} \mathrm{d} q\Big)\mathrm{d} r \mathrm{d}s\cr
&&+C\Lambda_{\alpha,B^H}^{0,T} \sup _{t \in[0,T]}  \int_{0}^t e^{-\rho t} \int_{0}^{s}\frac{(s-r)^{-\beta} r^{-\alpha}+(s-r)^{-\alpha-\beta}}{(t-s)^{1+\alpha-\beta}}(1+|X^{\varepsilon}_r|)\mathrm{d}r \mathrm{d}s\cr
&&+C\Lambda_{\alpha,B^H}^{0,T} \sup _{t \in[0,T]}  \int_{0}^t e^{-\rho t} \int_{0}^{s}(t-s)^{-1-\alpha+\beta}(s-r)^{-\beta}\Big(\int_{0}^r\frac{|X^{\varepsilon}_r-X^{\varepsilon}_q|}{(r-q)^{1+\alpha}}\mathrm{d}q\Big)\mathrm{d}r \mathrm{d}s\cr
&\leq&  C(1+|X_0|_{{\beta}})\cr
&&+C\Lambda_{\alpha,B^H}^{0,T}\sup _{t \in[0,T]}  \int_{0}^t e^{-\rho (t-r)}\frac{r^{-\alpha}+(t-r)^{-\alpha}}{(t-r)^{\alpha}}e^{-\rho r}(1+|X^{\varepsilon}_r|) \mathrm{d}r\cr
&&+C\Lambda_{\alpha,B^H}^{0,T} \sup _{t \in[0,T]}  \int_{0}^te^{-\rho (t-r)}(t-r)^{-\alpha}e^{-\rho r}\Big( \int_{0}^{r}\frac{|X^{\varepsilon}_r-X^{\varepsilon}_q|}{(r-q)^{1+\alpha}}\mathrm{d}q\Big)\mathrm{d}r \cr
&\leq&  C(1+|X_0|_{{\beta}})\cr
&&+C\Lambda_{\alpha,B^H}^{0,T}\Big(1+\rho^{2\alpha-1}\sup _{t \in[0, T]}e^{-\rho t}|X^{\varepsilon}_t|+\rho^{\alpha-1}\sup _{t \in[0, T]}e^{-\rho t} \int_{0}^t \frac{|X^{\varepsilon}_t-X^{\varepsilon}_q|}{(t-q)^{1+\alpha}}\mathrm{d}q\Big)\cr
&\leq& K \Lambda (1+\rho^{2 \alpha-1}\|X^{\varepsilon}\|_{\rho, T}+\rho^{\alpha-1}\|X^{\varepsilon}\|_{1, \rho, T}).
\end{eqnarray}

Putting $\rho=(4 K \Lambda)^{\frac{1}{1-\alpha}}$, we get from the inequality (\ref{in1}) that
\begin{eqnarray}\label{in4}
\|X^\varepsilon\|_{\rho, T} \leq \frac{4}{3} K \Lambda\left(1+\rho^{-1}\|X^\varepsilon\|_{1, \rho, T}\right).
\end{eqnarray}
Plugging this into the inequality (\ref{in2}) and making simple transformations, we arrive at
\begin{eqnarray*}
	\|X^\varepsilon\|_{1, \rho, T} \leq \frac{3}{2} K \Lambda+2 (K \Lambda)^{\frac{1}{1-\alpha}}\leq C \Lambda^{\frac{1}{1-\alpha}}.
\end{eqnarray*}
Substituting this into (\ref{in4}), we get $$\|X^\varepsilon\|_{\rho, T} \leq C  \Lambda^{\frac{1}{1-\alpha}}.$$
Thus, we have
\begin{eqnarray*}
	\sup_{t\in[0,T]} \| X^{\varepsilon}_t\|_{\alpha}&\leq& e^{\rho T}(\|X^\varepsilon\|_{\rho, T}+\|X^\varepsilon\|_{1,\rho, T}) \cr
	&\leq& C  \exp(C  \Lambda^{\frac{1}{1-\alpha}}) \Lambda^{\frac{1}{1-\alpha}}\cr
	&\leq& C  \exp\big((\Lambda_{\alpha,B^H}^{0,T})^{\frac{1}{1-\alpha}}\big).
\end{eqnarray*}

Since $0<\frac{1}{1-\alpha}<2$, by the classical Fernique's theorem, we have $$\mathbb{E}\big[\exp\big((\Lambda_{\alpha,B^H}^{0,T})^{\frac{1}{1-\alpha}}\big)\big]<\infty.$$
Then, the statement follows. \qed

Using similar techniques as for the Lemma \ref{xbound}, we have the following remark.
\begin{rem}\label{hatxbarx}
	Suppose that conditions {\rm (A1)-(A5)} and {\rm (B1)} are satisfied. Then, for any $p\geq 2$, we have
	\begin{eqnarray*}
		\sup_{t\in[0,T]} \big(\|\hat{X}^\varepsilon_t\|_{\alpha}+ \|\bar{X}_t\|_{\alpha}\big) &\leq& C  \exp\big((\Lambda_{\alpha,B^H}^{0,T})^{\frac{1}{1-\alpha}}\big),\cr
		\mathbb{E} \Big[\sup_{t\in[0,T]} \big(\|\hat{X}^\varepsilon_t\|_{\alpha}^p+ \|\bar{X}_t\|_{\alpha}^p\big)\Big]  &\leq& C_p.
	\end{eqnarray*}
\end{rem}

\begin{lem}\label{xh}
	Suppose that conditions {\rm (A1)-(A3)} and {\rm (B1)} are satisfied. Then, we have
	$$\mathbb{E}\big[|X^{\varepsilon}_t-X^{\varepsilon}_s |^2\big]\leq C|t-s|^{2\beta},$$
	where $C$ is independent of $\varepsilon, t, s$.
\end{lem}
\para{Proof:} From (\ref{mild2}), we have
\begin{eqnarray*}
	|X^{\varepsilon}_t-X^{\varepsilon}_s| &\leq& |(S_t-S_s)X_{0}|+\Big|\int_{s}^{t} S_{t-r} b(X^{\varepsilon}_r,Y^{\varepsilon}_r)\mathrm{d}r\Big| \cr
	&&+\Big|\int_{0}^{s}(S_{t-r}-S_{s-r}) b(X^{\varepsilon}_r,Y^{\varepsilon}_r)\mathrm{d}r\Big|\cr
	&&+\Big|\int_{s}^{t} S_{t-r} g(X^{\varepsilon}_r) \mathrm{d} B^H_r\Big| \cr
	&&+ \Big|\int_{0}^{s}(S_{t-r}-S_{s-r}) g(X^{\varepsilon}_r)\mathrm{d}B^H_r\Big|\cr
	&=:& \mathcal{V}_1+ \mathcal{V}_2+ \mathcal{V}_3+ \mathcal{V}_4+ \mathcal{V}_5.
\end{eqnarray*}

Since $X_{0}\in V_\beta$, by (\ref{semi2}) and (B1), we have
\begin{eqnarray*}
	\mathcal{V}_1+ \mathcal{V}_2 + \mathcal{V}_3 &\leq&|(S_{t-s}-\mathrm{id})S_s X_{0}|+\int_{s}^{t} |S_{t-r}b(X^{\varepsilon}_r,Y^{\varepsilon}_r)|\mathrm{d}r\cr
	&&+\int_{0}^{s}|(S_{t-r}-S_{s-r})b(X^{\varepsilon}_r,Y^{\varepsilon}_r)|\mathrm{d}r\cr
	&\leq&  C|X_{0}|_{{\beta}} (t-s)^\beta +C (t-s)+C (t-s)^\beta.
\end{eqnarray*}

Next, for $ \mathcal{V}_4$ and  $\mathcal{V}_5$, by Remark \ref{itofbm} and Lemma \ref{itofbm2}, we have
\begin{eqnarray*}
	\mathcal{V}_4+\mathcal{V}_5 &\leq& C \Lambda_{\alpha,B^H}^{0,T}\int_{s}^{t}[(r-s)^{-\alpha}+(t-r)^{-\alpha}](1+\|X^\varepsilon_r\|_{\alpha})\mathrm{d}r\cr
	&&+C(t-s)^{\beta}\Lambda_{\alpha,B^H}^{0,T}\int_{0}^s  [(s-r)^{-\beta} r^{-\alpha}+(s-r)^{-\alpha-\beta}](1+\|X^\varepsilon_r\|_{\alpha})\mathrm{d}r\cr
	&\leq& C[(t-s)^{1-\alpha}+(t-s)^{\beta}]\Lambda_{\alpha,B^H}^{0,T}\sup_{t\in[0,T]} \| X^{\varepsilon}_t\|_{\alpha}.
\end{eqnarray*}

Finally, by Lemma \ref{xbound}, we have
$$\mathbb{E}\big[|X^{\varepsilon}_t-X^{\varepsilon}_s |^2\big]\leq C|t-s|^{2\beta}.$$ Then,
the statement follows.
\qed

\begin{lem}\label{ybound}
	Suppose that conditions {\rm (A1)-(A5) and (B1)} are satisfied. Then, we have
	$$\sup_{t\in[0,T]}\mathbb{E} \big[| Y^{\varepsilon}_t |^{2}\big] \leq C,$$
	where $C$ is a positive constant which is independent of $\varepsilon$.
\end{lem}

\para{Proof:} Note that
\begin{eqnarray*}
	\sup_{t\in[0,T]}\mathbb{E}\big[|Y^{\varepsilon}_t|^{2}\big]
	& \leq& \sup_{t\in[0,T]} \mathbb{E}\big[|S_{\frac{t}{\varepsilon}}Y_{0}|^{2}\big]+\sup_{t\in[0,T]} \mathbb{E}\Big[\Big|\frac{1}{\varepsilon}\int_0^tS_{\frac{t-s}{\varepsilon}}F(X^{\varepsilon}_s,Y^{\varepsilon}_s)\mathrm{d}s\Big|^2\Big]\cr
	&&+\sup_{t\in[0,T]}\mathbb{E} \Big[\Big|\frac{1}{\sqrt{\varepsilon}} \int_{0}^{t} S_{\frac{t-s}{\varepsilon}}G(X^{\varepsilon}_s,Y^{\varepsilon}_s)\mathrm{d} W_s\Big|^{2}\Big]\cr
	&=:&\mathcal{I}_1+\mathcal{I}_2+\mathcal{I}_3.
\end{eqnarray*}

It is clear that $\mathcal{I}_1\leq e^{-2 \bar{\lambda}_1 \frac{t}{\varepsilon}}|Y_{0}  |^{2}< \infty.$

Next, we have
\begin{eqnarray*}
	\mathcal{I}_2&\leq&\sup_{t\in[0,T]} \Big(\frac{1}{\varepsilon}\int_0^te^{-\bar{\lambda}_1\frac{t-s}{\varepsilon}}\mathrm{d}s\times  \frac{1}{\varepsilon}\int_0^te^{-\bar{\lambda}_1\frac{t-s}{\varepsilon}}\mathbb{E}\big[|F(X^{\varepsilon}_s,Y^{\varepsilon}_s)  |^{2}\big]\mathrm{d}s\Big)\cr
	&\leq& \sup_{t\in[0,T]}\Big(\big(1+\mathbb{E}\big[|X^{\varepsilon}_t|^{2}\big]+\mathbb{E}\big[|Y^{\varepsilon}_t|^{2}\big]\big)\frac{C_3}{\bar{\lambda}_1 \varepsilon}\int_0^te^{-\bar{\lambda}_1\frac{t-s}{\varepsilon}}\mathrm{d}s\Big)\cr
	&\leq&  \frac{C_3}{\bar{\lambda}_1^2}\sup_{0\leq t \leq T}\mathbb{E}\big[|Y^{\varepsilon}_t|^{2}\big]+C.
\end{eqnarray*}

Then, (B1), It\^{o} isometry and \lemref{xbound} yield
 \begin{eqnarray*}
	\mathcal{I}_3
	\leq \sup_{t\in[0,T]}\mathbb{E}\big[|G(X^{\varepsilon}_t,Y^{\varepsilon}_t)|^{2}\big] \int_0^t\frac{C_3}{{\varepsilon}}e^{-2\bar{\lambda}_1\frac{t-s}{\varepsilon}}\mathrm{d}s
	\leq C.
\end{eqnarray*}

Therefore, due to (A5), $C_3< 2\bar{\lambda}_1^2$, we have
$$\sup_{t\in[0,T]}\mathbb{E} \big[| Y^{\varepsilon}_t |^{2}\big] \leq C.$$Then,
the statement follows.  \qed

Using similar techniques, under conditions {\rm (A1)-(A5) and (B1)}, we can prove $$\sup_{t\in[0,T]}\mathbb{E} \big[| \hat{Y}^{\varepsilon}_t |^{2}\big] \leq C.$$
Here, $C$ is also a positive constant which is independent of $\varepsilon$.

Now, using the definitions of $Y^{\varepsilon}_t$ (\ref{mild3}) and $\hat{Y}^{\varepsilon}_t$ (\ref{yj}),
we proceed to estimate  $\mathbb{E}\big[|Y^{\varepsilon}_t-\hat{Y}^{\varepsilon}_t|^{2}\big]$.
\begin{lem}\label{y-yhat}
	Suppose that conditions {\rm (A1)-(A5)} are satisfied. Then, we have
	\begin{eqnarray}
	\int_{k\delta}^{\min\{(k+1)\delta,T\}}
	\mathbb{E}[|Y^{\varepsilon}_t-\hat{Y}^{\varepsilon}_t|^{2}]\mathrm{d}t
	\leq C \varepsilon \delta^{2 \beta}e^{C\frac{\delta^2 }{\varepsilon^2}},
	\end{eqnarray}
	where $C$ is a constant which is indenpendent of $\varepsilon,\delta,k$.
\end{lem}
\para{Proof:} The resetting of the auxiliary process at the breakpoints $k \delta$ implies that
$Y^{\varepsilon}_{k\delta}=\hat{Y}^{\varepsilon}_{k\delta},$ for all $k$. Then, for $t\in [k\delta, \min\{(k+1)\delta,T\}]$, we start with
\begin{eqnarray}
e^{2 \bar{\lambda}_1 \frac{t}{\varepsilon}}\mathbb{E}[|Y^{\varepsilon}_t-\hat{Y}^{\varepsilon}_t|^{2}]&\leq& e^{2 \bar{\lambda}_1 \frac{t}{\varepsilon}}\mathbb{E}\Big[\Big|\frac{1}{\varepsilon}\int_{k\delta}^tS_{\frac{t-s}{\varepsilon}}(F(X^{\varepsilon}_s,Y^{\varepsilon}_s)-F(X^{\varepsilon}_{k\delta},\hat{Y}^{\varepsilon}_s))\mathrm{d}s\Big|^{2}\Big]\cr
&&+e^{2\bar{\lambda}_1 \frac{t}{\varepsilon}}\mathbb{E}\Big[\Big|\frac{1}{\sqrt{\varepsilon}}\int_{k\delta}^tS_{\frac{t-s}{\varepsilon}}(G(X^{\varepsilon}_s,Y^{\varepsilon}_s)-G(X^{\varepsilon}_{k\delta},\hat{Y}^{\varepsilon}_s))dW(s)\Big|^{2}\Big]\cr
&\leq& C (\frac{\delta }{\varepsilon^2}+\frac{1}{\varepsilon})\int_{k\delta}^te^{2 \bar{\lambda}_1 \frac{s}{\varepsilon}}\mathbb{E}[(|X^{\varepsilon}_s-X^{\varepsilon}_{k\delta}|^{2}+|Y^{\varepsilon}_s-\hat{Y}^{\varepsilon}_s|^{2})]\mathrm{d}s\cr
&\leq& C  (\frac{ \delta }{\varepsilon^2}+\frac{1}{\varepsilon})\Big(\delta^{2 \beta}\int_{k\delta}^te^{2 \bar{\lambda}_1 \frac{s}{\varepsilon}}\mathrm{d}s+\int_{k\delta}^t e^{2 \bar{\lambda}_1 \frac{s}{\varepsilon}}\mathbb{E} [|Y^{\varepsilon}_s-\hat{Y}^{\varepsilon}_s|^{2}] \mathrm{d}s\Big).
\end{eqnarray}

By Gr\"onwall's inequality, we have
\begin{eqnarray}
e^{2 \bar{\lambda}_1 \frac{t}{\varepsilon}}\mathbb{E}[|Y^{\varepsilon}_t-\hat{Y}^{\varepsilon}_t|^{2}]&\leq& C \delta^{2 \beta} \frac{ \delta+\varepsilon }{\varepsilon^2}\Big(\int_{k\delta}^te^{2 \bar{\lambda}_1 \frac{s}{\varepsilon}}\mathrm{d}s\Big) e^{C (\frac{ \delta }{\varepsilon^2}+\frac{1}{\varepsilon})(t-k\delta)}.\end{eqnarray}
It is clear that
\begin{eqnarray}\label{y-y}
\mathbb{E}[|Y^{\varepsilon}_t-\hat{Y}^{\varepsilon}_t|^{2}]
&\leq& C \delta^{2 \beta} \frac{ \delta+\varepsilon }{\varepsilon}e^{C (\frac{ \delta }{\varepsilon^2}+\frac{1}{\varepsilon})(t-k\delta)}.
\end{eqnarray}

Integrate (\ref{y-y}) from $k\delta$ to $\min\{(k+1)\delta,T\}$, we have
\begin{eqnarray*}
	\int_{k\delta}^{\min\{(k+1)\delta,T\}} \mathbb{E}[|Y^{\varepsilon}_t-\hat{Y}^{\varepsilon}_t|^{2}] \mathrm{d}t
	&\leq& C \delta^{2 \beta} \frac{ \delta+\varepsilon }{\varepsilon} \int_{k\delta}^{\min\{(k+1)\delta,T\}}  e^{C(\frac{ \delta }{\varepsilon^2}+\frac{1}{\varepsilon})(t-k\delta)}\mathrm{d}t\cr
	&\leq&C \delta^{2 \beta} \frac{ \delta+\varepsilon }{\varepsilon} \frac{1}{(\frac{ \delta }{\varepsilon^2}+\frac{1}{\varepsilon})}(e^{C(\frac{ \delta }{\varepsilon^2}+\frac{1}{\varepsilon})\delta}-1)\cr
	&\leq& C \varepsilon \delta^{2 \beta} e^{C\frac{\delta^2 }{\varepsilon^2}}.
\end{eqnarray*}
This completes  the proof of \lemref{y-yhat}. \qed

\subsection{The Proof of Theorem \ref{mainthm}}

We divide the  proof into three steps.
\para{Step 1: This step will estimate $\|\hat{X}^{\varepsilon}-X^{\varepsilon}\|_{\alpha,T}$.}
By (\ref{mild2}) and (\ref{xj}), we have	
\begin{eqnarray*}
	\mathbb{E}\big[\|{X}^{\varepsilon}-\hat{X}^{\varepsilon}\|_{\alpha,T}^2\big] &\leq&C \mathbb{E}\Big[\Big\|\int_{0}^{\cdot}S_{\cdot-s}(b(X^{\varepsilon}_s, \hat Y^{\varepsilon}_s)-b(X^{\varepsilon}_{s(\delta)},\hat{Y}^{\varepsilon}_s))\mathrm{d}s\Big\|^2_{\alpha,T}\Big]\cr
	&&+C \mathbb{E}\Big[\Big\|\int_{0}^{\cdot}S_{\cdot-s}(b(X^{\varepsilon}_s, Y^{\varepsilon}_s)-b(X^{\varepsilon}_s,\hat{Y}^{\varepsilon}_s))\mathrm{d}s\Big\|^2_{\alpha,T}\Big]\cr
	&=:& I_1+I_2.
\end{eqnarray*}

Following the similar steps as for the terms $\mathcal{M}_{21}$ and $\mathcal{M}_{22}$ in Appendix B, we obtain
\begin{eqnarray*}
	I_1 &\leq&C \int_{0}^{T}\mathbb{E}\big[|X^{\varepsilon}_r-X^{\varepsilon}_{r(\delta)}|^2\big]\mathrm{d}r\cr
	&&+C\sup _{t \in[0, T]} \mathbb{E}\Big[\Big(\int_{0}^{t}\int_0^r (t-s)^{-1-\alpha}\mathrm{d}s |b(X^{\varepsilon}_r, \hat Y^{\varepsilon}_r)-b(X^{\varepsilon}_{r(\delta)},\hat{Y}^{\varepsilon}_r)|\mathrm{d}r\Big)^2\Big]\cr
	&&+C\sup _{t \in[0, T]} \mathbb{E}\Big[\Big(\int_{0}^t \int_{r}^{t}(t-s)^{-1-\alpha+\beta}(s-r)^{-\beta}\mathrm{d}s |b(X^{\varepsilon}_r, \hat Y^{\varepsilon}_r)-b(X^{\varepsilon}_{r(\delta)},\hat{Y}^{\varepsilon}_r)|\mathrm{d}r\Big)^2\Big]\cr
	&\leq&C\int_{0}^{T}\mathbb{E}\big[|X^{\varepsilon}_r-X^{\varepsilon}_{r(\delta)}|^2\big]\mathrm{d}r\cr
	&\leq&C \delta^{2\beta}.
\end{eqnarray*}

Next, for $I_2$, we have
\begin{eqnarray*}
	I_2 &\leq&C\mathbb{E}\Big[\sup _{t \in[0, T]}\Big|\int_{0}^{t}S_{t-r}(b(X^{\varepsilon}_r, Y^{\varepsilon}_r)-b(X^{\varepsilon}_r,\hat{Y}^{\varepsilon}_r))\mathrm{d}r\Big|^2\Big]\cr
	&&+C\sup _{t \in[0, T]} \mathbb{E}\Big[\Big(\int_{0}^t (t-s)^{-1-\alpha}\Big|\int_{s}^{t}S_{t-r} (b(X^{\varepsilon}_r, Y^{\varepsilon}_r)-b(X^{\varepsilon}_r,\hat{Y}^{\varepsilon}_r)) \mathrm{d}r\Big|\mathrm{d}s\Big)^2\Big]\cr
	&&+C\sup _{t \in[0, T]} \mathbb{E}\Big[\Big(\int_{0}^t (t-s)^{-1-\alpha}\Big|\int_{0}^{s}(S_{t-r} -S_{s-r}) (b(X^{\varepsilon}_r, Y^{\varepsilon}_r)-b(X^{\varepsilon}_r,\hat{Y}^{\varepsilon}_r))d r\Big|\mathrm{d}s\Big)^2\Big]\cr
	&=:& I_{21}+I_{22}+I_{23}.
\end{eqnarray*}

If we set $\ell:=\{t\geq (\lfloor \frac{s}{\delta} \rfloor +2)\delta\}$, $\ell^c:=\{t < (\lfloor \frac{s}{\delta} \rfloor +2)\delta\}$, $\jmath:=\{\lfloor \frac{t}{\delta} \rfloor \leq 1\}$ and $\jmath^c:=\{\lfloor \frac{t}{\delta} \rfloor > 1\}$ and by H\"older's inequality and the fact that $(\lfloor \frac{t}{\delta}\rfloor-\lfloor \frac{s}{\delta} \rfloor -1)\leq \frac{t-s}{\delta}$, for any $0\leq s \leq t$, we obtain
\begin{eqnarray*}
	&&I_{21}+I_{22}\cr
	&\leq&C\mathbb{E}\Big[\sup _{t \in[0, T]}\Big|\int_{t(\delta)}^{t}S_{t-r} (b(X^{\varepsilon}_r, Y^{\varepsilon}_r)-b(X^{\varepsilon}_r,\hat{Y}^{\varepsilon}_r))\mathrm{d}r\Big|^2\Big]\cr
	&&+C \delta^{-1} \sum_{k=0}^{\lfloor \frac{T}{\delta} \rfloor-1} \mathbb{E}\Big[\Big| \int_{k\delta}^{(k+1)\delta} S_{t-r} (b(X^{\varepsilon}_r, Y^{\varepsilon}_r)-b(X^{\varepsilon}_r,\hat{Y}^{\varepsilon}_r))\mathrm{d}r\Big|^2\Big]\cr
	&&+C\sup _{t \in[0, T]}\mathbb{E}\Big[\int_{0}^t (t-s)^{-\frac{3}{2}-\alpha}\Big|\int_{s}^{t}S_{t-r} (b(X^{\varepsilon}_r, Y^{\varepsilon}_r)-b(X^{\varepsilon}_r,\hat{Y}^{\varepsilon}_r)) \mathrm{d}r\Big|^2\mathrm{d}s\Big]\cr
	&\leq&C\delta^2+C\delta^{-2}\max_{0\leq k \leq \lfloor \frac{T}{\delta} \rfloor-1} \mathbb{E}\Big[\Big| \int_{k\delta}^{(k+1)\delta} S_{t-r} (b(X^{\varepsilon}_r, Y^{\varepsilon}_r)-b(X^{\varepsilon}_r,\hat{Y}^{\varepsilon}_r))\mathrm{d}r\Big|^2\Big]\cr
	&&+C\sup _{t \in[0, T]} \Big\{\int_{0}^{t}(t-s)^{-\frac{1}{2}-\alpha} \int_{s}^{t} \mathbb{E}\big[|S_{t-r} (b(X^{\varepsilon}_r, Y^{\varepsilon}_r)-b(X^{\varepsilon}_r,\hat{Y}^{\varepsilon}_r))|^2\big]\mathrm{d}r\mathrm{d}s\Big\}\mathbf{1}_{\ell^c\bigcap \jmath}\cr
	&&+C\sup _{t \in[0, T]} \Big\{\int_{0}^{(\lfloor \frac{t}{\delta} \rfloor -1) \delta}(t-s)^{-\frac{3}{2}-\alpha} \mathbb{E}\Big[\Big|\int_{s}^{t} S_{t-r} (b(X^{\varepsilon}_r, Y^{\varepsilon}_r)-b(X^{\varepsilon}_r,\hat{Y}^{\varepsilon}_r))\mathrm{d}r\Big|^2\Big]\mathbf{1}_{\ell^c\bigcap \jmath^c}\mathrm{d}s\Big\}\cr
	&&+C\sup _{t \in[0, T]} \Big\{\int_{(\lfloor \frac{t}{\delta} \rfloor -1) \delta }^{t}(t-s)^{-\frac{1}{2}-\alpha} \int_{s}^{t} \mathbb{E}\big[|S_{t-r} (b(X^{\varepsilon}_r, Y^{\varepsilon}_r)-b(X^{\varepsilon}_r,\hat{Y}^{\varepsilon}_r))|^2\big]\mathrm{d}r \mathbf{1}_{\ell^c\bigcap \jmath^c} ds\Big\}\cr
	&&+C\sup _{t \in[0, T]} \Big\{\int_{0}^{t}(t-s)^{-\frac{1}{2}-\alpha} \int_{s}^{(\lfloor \frac{s}{\delta} \rfloor +1)\delta} \mathbb{E}\big[|S_{t-r} (b(X^{\varepsilon}_r, Y^{\varepsilon}_r)-b(X^{\varepsilon}_r,\hat{Y}^{\varepsilon}_r))|^2\big]\mathrm{d}r \mathbf{1}_{\ell} ds\Big\}\cr
	&&+C\sup _{t \in[0, T]} \Big\{\int_{0}^{t}(t-s)^{-\frac{1}{2}-\alpha} \int^{t}_{(\lfloor \frac{t}{\delta} \rfloor)\delta} \mathbb{E}\big[|S_{t-r} (b(X^{\varepsilon}_r, Y^{\varepsilon}_r)-b(X^{\varepsilon}_r,\hat{Y}^{\varepsilon}_r))|^2\big]\mathrm{d}r \mathbf{1}_{\ell} ds\Big\}\cr
	&&+C\sup _{t \in[0, T]} \Big\{\int_{0}^{t}\frac{\lfloor \frac{t}{\delta} \rfloor-\lfloor \frac{s}{\delta} \rfloor -1}{(t-s)^{\frac{3}{2}+\alpha}} \sum_{k=\lfloor \frac{s}{\delta} \rfloor +1}^{\lfloor \frac{t}{\delta} \rfloor-1} \mathbb{E}\Big[\Big| \int_{k\delta}^{(k+1)\delta } S_{t-r} (b(X^{\varepsilon}_r, Y^{\varepsilon}_r)-b(X^{\varepsilon}_r,\hat{Y}^{\varepsilon}_r))\mathrm{d}r\Big|^2\Big] \mathbf{1}_{\ell} \mathrm{d}s\Big\}\cr
	&\leq&C\delta^2+C\delta^{-2}\max_{0\leq k \leq \lfloor \frac{T}{\delta} \rfloor-1} \mathbb{E}\Big[\Big| \int_{k\delta}^{(k+1)\delta} S_{t-r} (b(X^{\varepsilon}_r, Y^{\varepsilon}_r)-b(X^{\varepsilon}_r,\hat{Y}^{\varepsilon}_r))\mathrm{d}r\Big|^2\Big]\cr
	&&+C\delta^{\frac{3}{2}-\alpha}+C\sup _{t \in[0, T]} \Big\{\delta^2 \int_{0}^{(\lfloor \frac{t}{\delta} \rfloor -1) \delta}(t-s)^{-\frac{3}{2}-\alpha}\mathbf{1}_{\ell^c\bigcap \jmath^c} \mathrm{d}s\Big\}+C\delta\cr
	&&+C \delta^{-1} \sup _{t \in[0, T]} \Big\{\int_{0}^{t}\frac{\sum_{k=\lfloor \frac{s}{\delta} \rfloor +1}^{\lfloor \frac{t}{\delta} \rfloor-1} \mathbb{E}\big[\big| \int_{k\delta}^{(k+1)\delta} S_{t-r} (b(X^{\varepsilon}_r, Y^{\varepsilon}_r)-b(X^{\varepsilon}_r,\hat{Y}^{\varepsilon}_r))\mathrm{d}r\big|^2\big]}{(t-s)^{\alpha+\frac{1}{2}}} \mathbf{1}_{\ell} \mathrm{d}s\Big\}\cr
	&\leq&C\delta+C\delta^{-2} \max_{0\leq k \leq \lfloor \frac{T}{\delta} \rfloor-1} \mathbb{E}\Big[\Big| \int_{k\delta}^{(k+1)\delta} S_{t-r} (b(X^{\varepsilon}_r, Y^{\varepsilon}_r)-b(X^{\varepsilon}_r,\hat{Y}^{\varepsilon}_r))\mathrm{d}r\Big|^2\Big].
\end{eqnarray*}

Once again, by H\"older's inequality and  inequality (\ref{semi2}) and taking $\beta'\in(\frac{1+2\alpha}{4},\frac{1}{2})$, we have
\begin{eqnarray*}
	I_{23}&\leq&C\sup _{t \in[0, T]}\Big\{\int_{0}^{t}\frac{\mathbb{E}\big[\big|\int_{0}^{s}(S_{t-r} -S_{s-r}) (b(X^{\varepsilon}_r, Y^{\varepsilon}_r)-b(X^{\varepsilon}_r,\hat{Y}^{\varepsilon}_r)) \mathrm{d} r\big|^2\big]}{(t-s)^{\frac{3}{2}+\alpha}}\mathrm{d}s\Big\}\cr
	&\leq & C \sup _{t \in[0, T]}\Big\{\int_{0}^{t}\frac{\mathbb{E}\big[\big(\int_{0}^s \|(-A)^{\beta'}  (S_{s-r}-S_{s-r(\delta)})\||(b(X^{\varepsilon}_r, Y^{\varepsilon}_r)-b(X^{\varepsilon}_r,\hat{Y}^{\varepsilon}_r))| \mathrm{d}r\big)^2\big]}{(t-s)^{\frac{3}{2}+\alpha-2\beta' }}\mathrm{d}s\Big\}\cr
	&&+C\sup _{t \in[0, T]}\Big\{\int_{0}^{t}\frac{\mathbb{E}\big[\big|\int_{0}^s (-A)^{\beta'} S_{s-r(\delta)}(b(X^{\varepsilon}_r, Y^{\varepsilon}_r)-b(X^{\varepsilon}_r,\hat{Y}^{\varepsilon}_r)) \mathrm{d}r\big|^2\big]}{(t-s)^{\frac{3}{2}+\alpha-2\beta' }}\mathrm{d}s\Big\}\cr
	&\leq & C \delta^{2\beta'} \sup _{t \in[0, T]}\Big\{\int_{0}^{t}\frac{\mathbb{E}\big[\big(\int_{0}^s\| (-A)^{2\beta'}S_{s-r}\| |(b(X^{\varepsilon}_r, Y^{\varepsilon}_r)-b(X^{\varepsilon}_r,\hat{Y}^{\varepsilon}_r))| \mathrm{d}r\big)^2\big]}{(t-s)^{\frac{3}{2}+\alpha-2\beta' }}\mathrm{d}s\Big\}\cr
	&&+ C \sup _{t \in[0, T]}\Big\{\int_{0}^{t}\frac{\mathbb{E}\big[\big|\int_{\lfloor \frac{s}{\delta} \rfloor \delta}^s (-A)^{\beta'} S_{s-\lfloor \frac{r}{\delta} \rfloor \delta} (b(X^{\varepsilon}_r, Y^{\varepsilon}_r)-b(X^{\varepsilon}_r,\hat{Y}^{\varepsilon}_r)) \mathrm{d}r\big|^2\big]}{(t-s)^{\frac{3}{2}+\alpha-2\beta' }}\mathrm{d}s\Big\}\cr
	&&+ C \sup _{t \in[0, T]}\Big\{\int_{0}^{t}\frac{\mathbb{E}\big[\big|\sum_{k=0}^{\lfloor \frac{s}{\delta} \rfloor -1}\int_{k\delta}^{(k+1)\delta} (-A)^{\beta'} S_{s-k\delta} (b(X^{\varepsilon}_r, Y^{\varepsilon}_r)-b(X^{\varepsilon}_r,\hat{Y}^{\varepsilon}_r)) \mathrm{d}r\big|^2\big]}{(t-s)^{\frac{3}{2}+\alpha-2\beta' }}\mathrm{d}s\Big\}\cr
	&\leq & C \delta^{2\beta'} \sup _{t \in[0, T]}\Big\{\int_{0}^{t}\frac{\big(\int_{0}^s(s-r)^{-2\beta'}\mathrm{d}r\big)^2}{(t-s)^{\frac{3}{2}+\alpha-2\beta' }}\mathrm{d}s\Big\}+ C\sup _{t \in[0, T]}\Big\{\int_{0}^{t}\frac{\big(\int_{\lfloor \frac{s}{\delta} \rfloor \delta}^s (s-\lfloor \frac{s}{\delta} \rfloor \delta)^{-\beta'}\mathrm{d}r\big)^2}{(t-s)^{\frac{3}{2}+\alpha-2\beta' }}\mathrm{d}s\Big\}\cr
	&&+ C \sup _{t \in[0, T]}\Big\{\int_{0}^{t}\frac{\mathbb{E}\big[\big|\sum_{k=0}^{\lfloor \frac{s}{\delta} \rfloor -1}\int_{k\delta}^{(k+1)\delta} (-A)^{\beta'} S_{s-k\delta} (b(X^{\varepsilon}_r, Y^{\varepsilon}_r)-b(X^{\varepsilon}_r,\hat{Y}^{\varepsilon}_r)) \mathrm{d}r\big|^2\big]}{(t-s)^{\frac{3}{2}+\alpha-2\beta' }}\mathrm{d}s\Big\}\cr
	&\leq& C\delta^{2\beta'}+C \delta^{2-2\beta'}\cr
	&&+C \delta^{-1}\sup _{t \in[0, T]} \Big\{\int_{0}^{t}\frac{\sum_{k=0}^{\lfloor \frac{s}{\delta} \rfloor -1}\big[(s-k\delta)^{-2\beta'}\mathbb{E}\big[\big|\int_{k\delta}^{(k+1)\delta} (b(X^{\varepsilon}_r, Y^{\varepsilon}_r)-b(X^{\varepsilon}_r,\hat{Y}^{\varepsilon}_r)) \mathrm{d}r\big|^2\big]}{(t-s)^{\frac{3}{2}+\alpha-2\beta' }}\mathrm{d}s\Big\}.
\end{eqnarray*}

As a consequence, we have
\begin{eqnarray}\label{I2}
I_2 &\leq& C \delta^{-1}\sup _{t \in[0, T]} \Big\{\int_{0}^{t}\frac{\sum_{k=0}^{\lfloor \frac{s}{\delta} \rfloor -1}\big((s-k\delta)^{-2\beta'}\mathbb{E}\big[\big|\int_{k\delta}^{(k+1)\delta} (b(X^{\varepsilon}_r, Y^{\varepsilon}_r)-b(X^{\varepsilon}_r,\hat{Y}^{\varepsilon}_r)) \mathrm{d}r\big|^2\big]\big)}{(t-s)^{\frac{3}{2}+\alpha-2\beta' }}\mathrm{d}s\Big\}\cr
&&+C \delta^{2\beta'}+C\delta^{-2} \max_{0\leq k \leq \lfloor \frac{T}{\delta} \rfloor-1} \mathbb{E}\Big[\Big| \int_{k\delta}^{(k+1)\delta} S_{t-r} (b(X^{\varepsilon}_r, Y^{\varepsilon}_r)-b(X^{\varepsilon}_r,\hat{Y}^{\varepsilon}_r))\mathrm{d}r\Big|^2\Big].
\end{eqnarray}

Note that for $2\beta'<1$,
\begin{eqnarray*}
	\sum_{k=0}^{\lfloor \frac{s}{\delta} \rfloor -1}  (s-k\delta)^{-2\beta'} &\leq& \delta^{-2\beta'}\sum_{k=1}^{\lfloor \frac{s}{\delta} \rfloor}  k^{-2\beta'}
	=\delta^{-2\beta'}\sum_{k=1}^{\lfloor \frac{s}{\delta} \rfloor}  \int_{k-1}^k k^{-2\beta'}dv \cr
	&\leq& \delta^{-2\beta'}\int_{0}^{\lfloor \frac{s}{\delta} \rfloor }v^{-2\beta'}dv
	\leq \delta^{-2\beta'}\frac{1}{1-2\beta'} \lfloor \frac{s}{\delta} \rfloor ^{1-2\beta'}\cr
	&\leq& C \delta^{-1},
\end{eqnarray*}
holds, then by Lemma \ref{y-yhat} and applying H\"older's inequality again
\begin{eqnarray*}
	I_2 &\leq&C \delta^{2\beta'}+C \delta^{-1}\max_{0\leq k \leq \lfloor \frac{T}{\delta} \rfloor-1} \int_{k\delta}^{(k+1)\delta} \mathbb{E}\big[|Y^{\varepsilon}_r-\hat{Y}^{\varepsilon}_r|^2\big]\mathrm{d}r\cr
	&&+C \sup _{t \in[0, T]} \Big\{\int_{0}^{t}\frac{\sum_{k=0}^{\lfloor \frac{s}{\delta} \rfloor -1}(s-k\delta)^{-2\beta'}\int_{k\delta}^{(k+1)\delta}\mathbb{E}\big[|Y^{\varepsilon}_r-\hat{Y}^{\varepsilon}_r|^2\big]\mathrm{d}r}{(t-s)^{\frac{3}{2}+\alpha-2\beta'}}\mathrm{d}s\Big\}\cr
	&\leq&C \delta^{2\beta'}+C \varepsilon \delta^{2 \beta}e^{C\frac{\delta^2}{\varepsilon^2}}\Big(\delta^{-1}+\sup _{t \in[0, T]}\Big\{ \int_{0}^{t}\frac{\sum_{k=0}^{\lfloor \frac{s}{\delta} \rfloor -1}(s-k\delta)^{-2\beta'}}{(t-s)^{\frac{3}{2}+\alpha-2\beta' }}\mathrm{d}s\Big\}\Big)\cr
	&\leq&C \delta^{2\beta'}+C \varepsilon\delta^{-1} \delta^{2 \beta}e^{C\frac{\delta^2}{\varepsilon^2}}.
\end{eqnarray*}

Therefore, we have
\begin{eqnarray}\label{x-hatx}
\mathbb{E}[\|{X}^{\varepsilon}-\hat{X}^{\varepsilon}\|_{\alpha,T}^2]\leq C \big(\varepsilon \delta^{-1} \delta^{2 \beta}e^{C\frac{\delta^2 }{\varepsilon^2}}+\delta^{2\beta'}\big).
\end{eqnarray}

This completes the proof of Step 1. \qed
\para{Step 2:} This step will estimate $\|\hat{X}^{\varepsilon}-\bar{X}\|_{\alpha,T}$.
\begin{lem}\label{ptau}
	The following inequality hold:
	$$\mathbb{P}(\tau_N < T) \leq N^{-1}\mathbb{E}\big[\Lambda_{\alpha,B^H}^{0,T}\big],$$
	and it tends to $0$ when $N \rightarrow \infty$.
\end{lem}
\para{Proof:} By Chebyshev's inequality, we have
\begin{eqnarray*}
	\mathbb{P}(\tau_N < T)&\leq& \mathbb{P}\Big(\sum_{i=1}^\infty\sqrt{\lambda_i} \Lambda_\alpha^{0,T}(\beta_i^H) \geq N\Big) \leq N^{-1}{\sum_{i=1}^\infty\sqrt{\lambda_i} \mathbb{E} \big[\Lambda_\alpha^{0,T}(\beta_i^H)}\big].
\end{eqnarray*}

Because $\Lambda_\alpha^{0,T}(\beta_i^H)$ has moments of all order, see \cite[Lemma 7.5]{Rascanu2002}, thus we have
$$\lim_{N \rightarrow \infty}N^{-1}\mathbb{E}\big[\Lambda_{\alpha,B^H}^{0,T}\big]=0.$$
This completes the proof of Lemma \ref{ptau}. \qed	

Then, by (\ref{mild2}) and (\ref{ave-eq}), we get
\begin{eqnarray*}
	\mathbb{E}\big[\|\hat{X}^{\varepsilon}-\bar{X}\|_{\alpha,T}^2\big]
	\leq \mathbb{E}\big[\|\hat{X}^{\varepsilon}-\bar{X}\|^2_{\alpha,T}\mathbf{1}_{\{\tau_N < T\}}\big]+\mathbb{E}\big[\|\hat{X}^{\varepsilon}-\bar{X}\|^2_{\alpha,T}\mathbf{1}_{\{\tau_N \geq T\}}\big].
\end{eqnarray*}

For the first term on the right-hand side of above inequality, by Chebyshev's inequality, we have
\begin{eqnarray}
\mathbb{E}\big[\|\hat{X}^{\varepsilon}-\bar{X}\|^2_{\alpha,T}\mathbf{1}_{\{\tau_N < T\}}\big] \leq \sqrt{\mathbb{E}[\|\hat{X}^{\varepsilon}-\bar{X}\|^4_{\alpha,T}]}\cdot \sqrt{\mathbb{P}(\tau_N < T)}.
\end{eqnarray}

It follows from Lemma \ref{ptau} that $\mathbb{P}(\tau_N < T) \leq N^{-1}\mathbb{E}\big[\Lambda_{\alpha,B^H}^{0,T}\big]$. Then, by Lemma \ref{xbound}, summing up all bounds we obtain
\begin{eqnarray}\label{j0}
\mathbb{E}\big[\|\hat{X}^{\varepsilon}-\bar{X}\|^2_{\alpha,T}\mathbf{1}_{\{\tau_N < T\}}\big]\leq C \sqrt{N^{-1}\mathbb{E}[\Lambda_{\alpha,B^H}^{0,T}]}.
\end{eqnarray}

For the second term, set $A_{N,T}=\{\Lambda_{\alpha,B^H}^{0,T}\leq N\}$, we have
\begin{eqnarray*}
	\mathbb{E}\big[\|\hat{X}^{\varepsilon}-\bar{X}\|^2_{\alpha,T}\mathbf{1}_{A_{N,T}}\big]
	&\leq& C \mathbb{E}\Big[\Big\|\int_{0}^{\cdot}S_{\cdot-s}(b(X^{\varepsilon}_{s(\delta)}, \hat Y^{\varepsilon}_s)-\bar{b}(X^{\varepsilon}_{s(\delta)}))\mathrm{d}s\Big\|^2_{\alpha,T}\mathbf{1}_{A_{N,T}}\Big]\cr
	&&+ C \mathbb{E}\Big[\Big\|\int_{0}^{\cdot}S_{\cdot-s}(\bar{b}(X^{\varepsilon}_{s(\delta)})-\bar{b}(X^{\varepsilon}_s))\mathrm{d}s\Big\|^2_{\alpha,T}\mathbf{1}_{A_{N,T}}\Big]\cr
	&&+ C \mathbb{E}\Big[\Big\|\int_{0}^{\cdot}S_{\cdot-s}(\bar{b}(X^{\varepsilon}_s)-\bar{b}(\hat{X}^{\varepsilon}_s))\mathrm{d}s\Big\|^2_{\alpha,T}\mathbf{1}_{A_{N,T}}\Big]\cr
	&&+C \mathbb{E}\Big[\Big\|\int_{0}^{\cdot}S_{\cdot-s}(\bar{b}(\hat{X}^{\varepsilon}_s)-\bar{b}(\bar{X}_s))\mathrm{d}s\Big\|^2_{\alpha,T}\mathbf{1}_{A_{N,T}}\Big]\cr
	&&+C \mathbb{E}\Big[\Big\|\int_{0}^{\cdot}S_{\cdot-s}(g(X^{\varepsilon}_s)-g(\hat{X}^{\varepsilon}_{s}))\mathrm{d}B^H_s\Big\|^2_{\alpha,T}\mathbf{1}_{A_{N,T}}\Big]\cr
	&&+C \mathbb{E}\Big[\Big\|\int_{0}^{\cdot}S_{\cdot-s}(g(\hat{X}^{\varepsilon}_s)-g(\bar{X}_{s}))\mathrm{d}B^H_s\Big\|^2_{\alpha,T}\mathbf{1}_{A_{N,T}}\Big]\cr
	&=:& \sum_{i=1}^{6}J_i.
\end{eqnarray*}

For $J_1$, we only need to replace the estimate for $(b(X^{\varepsilon}_r, Y^{\varepsilon}_r)-b(X^{\varepsilon}_r,\hat{Y}^{\varepsilon}_r))$ which appeared in the estimation of $I_2$ in (\ref{I2}) by the corresponding estimate of $(b(X^{\varepsilon}_{r(\delta)}, \hat Y^{\varepsilon}_r)-\bar{b}(X^{\varepsilon}_{r(\delta)}))$.
\begin{eqnarray}\label{J1}
J_1 &\leq&C \delta^{-1}\sup _{t \in[0, T]} \int_{0}^{t}\frac{\sum_{k=0}^{\lfloor \frac{s}{\delta} \rfloor -1}\big((s-k\delta)^{-2\beta'}\mathbb{E}\big[\big|\int_{k\delta}^{(k+1)\delta} (b(X^{\varepsilon}_{k\delta}, \hat Y^{\varepsilon}_r)-\bar{b}(X^{\varepsilon}_{k\delta}))\mathrm{d}r\big|^2\big]\big)}{(t-s)^{\frac{3}{2}+\alpha-2\beta' }}\mathrm{d}s\cr
&&+C \delta^{2\beta'}+C \delta^{-2} \max_{0\leq k \leq \lfloor \frac{T}{\delta} \rfloor-1} \mathbb{E}\Big[\Big| \int_{k\delta}^{(k+1)\delta} S_{t-r} (b(X^{\varepsilon}_{k\delta}, \hat Y^{\varepsilon}_r)-\bar{b}(X^{\varepsilon}_{k\delta}))\mathrm{d}r\Big|^2\Big] .\end{eqnarray}

By a time shift transformation, we note that it follows from the definition of $\hat{Y}^{\varepsilon}_s$ that for $s\in[0,\delta]$, we have
\begin{eqnarray}\label{construct1}
\qquad \hat{Y}^{\varepsilon}_{k\delta+s} =S_{\frac{s}{\varepsilon}}\hat{Y}^{\varepsilon}_{k\delta}+\frac{1}{\varepsilon}\int_{0}^{s}S_{\frac{s-r}{\varepsilon}}F(X^{\varepsilon}_{k\delta},\hat{Y}^{\varepsilon}_{r+k\delta})\mathrm{d}r+\frac{1}{\sqrt{\varepsilon}}\int_{0}^{s}S_{\frac{s-r}{\varepsilon}}G(X^{\varepsilon}_{k\delta},\hat{Y}^{\varepsilon}_{r+k\delta})\mathrm{d}W^{*}_r,
\end{eqnarray}
where $W^{*}_r:=W_{r+k\delta}-W_{k\delta}$ is the shift of $W_r$, both of which
have the same distribution. Let $\bar{W}$ be a $V$-valued Brownian motion defined on the same stochastic basis and independent of $(B^{H},W)$. Construct a process  $Y^{X^{\varepsilon}_{k\delta},\hat{Y}^{\varepsilon}_{k\delta}}_s$   by means of	
\begin{eqnarray}\label{construct2}
Y^{X^{\varepsilon}_{k\delta},\hat{Y}^{\varepsilon}_{k\delta}}_{\frac{s}{\varepsilon}}&=&S_{\frac{s}{\varepsilon}}\hat{Y}^{\varepsilon}_{k\delta}+
\int_{0}^{\frac{s}{\varepsilon}}S_{\frac{s}{\varepsilon}-r}F(X^{\varepsilon}_{k\delta},Y^{X^{\varepsilon}_{k\delta},\hat{Y}^{\varepsilon}_{k\delta}}_r)\mathrm{d}r +\int_{0}^{\frac{s}{\varepsilon}}S_{\frac{s}{\varepsilon}-r}G(X^{\varepsilon}_{k\delta},Y^{X^{\varepsilon}_{k\delta},
	\hat{Y}^{\varepsilon}_{k\delta}}_r) \mathrm{d} \bar{W}_r\cr
&=&S_{\frac{s}{\varepsilon}}\hat{Y}^{\varepsilon}_{k\delta}+\frac{1}{\varepsilon}\int_{0}^{s}
S_{\frac{s-r}{\varepsilon}}F(X^{\varepsilon}_{k\delta},Y^{X^{\varepsilon}_{k\delta},\hat{Y}^{\varepsilon}_{k\delta}} _{\frac{r}{\varepsilon}})\mathrm{d}r\cr
&&+\frac{1}{\sqrt{\varepsilon}}\int_{0}^{s}S_{\frac{s-r}{\varepsilon}}G(X^{\varepsilon}_{k\delta},
Y^{X^{\varepsilon}_{k\delta},\hat{Y}^{\varepsilon}_{k\delta}}_{\frac{r}{\varepsilon}})\mathrm{d}\bar{W}^{*}_r,
\end{eqnarray}	
where $\bar{W}^{*}_r:=\sqrt{\varepsilon}\bar{W}_{{r}/{\varepsilon}}$ is the shift of $\bar{W}_r$
with the same distribution. Because both $W^{*}$ and $\bar{W}^{*}$ are independent of $(X_{k\delta}^\varepsilon,\hat{Y}_{k\delta}^\varepsilon)$, 	
comparison of (\ref{construct1}) and (\ref{construct2})
yields
$$(X^{\varepsilon}_{k\delta},\{\hat{Y}^{\varepsilon}_{r+k\delta}\}_{r\in[0,\delta)})\sim(X^{\varepsilon}_{k\delta},\{Y
^{X^{\varepsilon}_{k\delta},\hat{Y}^{\varepsilon}_{k\delta}}_{\frac{r}{\varepsilon}}\}_{r\in[0,\delta)}),$$
where $\sim$ denotes a coincidence in distribution sense.	

To proceed, for $t\in[0,T]$, we have
\begin{eqnarray*}
	\mathbb{E}\Big[\Big|\int_{k\delta}^{(k+1)\delta }S_{t-r} (b(X^{\varepsilon}_{k \delta},\hat{Y}^{\varepsilon}_r)-\bar{b}(X^{\varepsilon}_{k\delta}))\mathrm{d}r\Big|^2\Big]
	\leq  C \varepsilon^2 \int_{0}^{\frac{\delta}{\varepsilon}}\int_{\tau}^{\frac{\delta}{\varepsilon}}\mathcal{J}_{k}(s,\tau)\mathrm{d}s \mathrm{d}\tau,
\end{eqnarray*}
where
\begin{eqnarray*}
	\mathcal{J}_{k}(s,\tau)&=&\mathbb{E}[\langle S_{t-k\delta-s\varepsilon}(b(X^{\varepsilon}_{k\delta},\hat{Y}^{\varepsilon}_{s\varepsilon+k\delta})
	-\bar{b}(X^{\varepsilon}_{k\delta})),S_{t-k\delta-\tau\varepsilon}(b(X^{\varepsilon}_{k\delta},\hat{Y}^{\varepsilon}_{\tau \varepsilon+k\delta})
	-\bar{b}(X^{\varepsilon}_{k\delta}))\rangle]\cr
	&=&\mathbb{E}[\langle S_{t-k\delta-s\varepsilon}(b(X^{\varepsilon}_{k\delta},Y^{X^{\varepsilon}_{k\delta},\hat{Y}^{\varepsilon}_{k\delta}}_{s})
	-\bar{b}(X^{\varepsilon}_{k\delta})), S_{t-k\delta-\tau\varepsilon}(b(X^{\varepsilon}_{k\delta},Y^{X^{\varepsilon}_{k\delta},\hat{Y}^{\varepsilon}_{k\delta}}_{\tau})
	-\bar{b}(X^{\varepsilon}_{k\delta}))\rangle].
\end{eqnarray*}

We now present a key lemma for an estimate of $\mathcal{J}_{k}(s,\tau)$ which will be proved in Appendix B.
\begin{lem}\label{claim}
	For any $k$, we have
	$$\mathcal{J}_{k}(s,\tau)\leq C e^{-\frac{\eta}{2}(s-\tau)}\mathbb{E}[(1+|X^{\varepsilon}_{k\delta}|^{2}+|\hat{Y}^{\varepsilon}_{k\delta}|^{2})],$$
	where $\eta$ is defined in condition {\rm (A5)} and $C>0$ is a constant independent of $\varepsilon,\delta,k,s,\tau$.
\end{lem}

According to Lemma \ref{claim} and by choosing $\delta=\delta(\varepsilon)$ such that $\frac{\delta}{\varepsilon}$ is sufficiently large, we have
\begin{eqnarray*}
	\mathbb{E}\Big[\Big|\int_{k\delta}^{(k+1)\delta }S_{t-r} (b(X^{\varepsilon}_{k \delta},\hat{Y}^{\varepsilon}_r)-\bar{b}(X^{\varepsilon}_{k\delta}))\mathrm{d}r\Big|^2\Big]
	&\leq&  C \varepsilon^2 \int_{0}^{\frac{\delta}{\varepsilon} }\int_{\tau}^{\frac{\delta}{\varepsilon} }\mathcal{J}_{k}(s,\tau)\mathrm{d}s\mathrm{d}\tau\cr
	&\leq& C\varepsilon^2 \int_{0}^{\frac{\delta}{\varepsilon}}\int_{\tau}^{\frac{\delta}{\varepsilon}}e^{-\frac{\eta}{2}(s-\tau)}\mathrm{d}s\mathrm{d}\tau\cr
	&\leq& C \varepsilon^2 (\frac{2}{\eta}\frac{\delta}{\varepsilon} -\frac{4}{\eta^{2}}+e^{\frac{-\eta}{2}\frac{\delta}{\varepsilon}}).
\end{eqnarray*}

Hence, we get
\begin{eqnarray}\label{gamma}
J_1  \leq C \delta^{2\beta'}+C  \delta^{-2}\varepsilon^2 (\frac{2}{\eta}\frac{\delta}{\varepsilon} -\frac{4}{\eta^{2}}+e^{\frac{-\eta}{2}\frac{\delta}{\varepsilon}})\leq C (\varepsilon \delta^{-1}+\delta^{2\beta'}).
\end{eqnarray}

By similar caculations as for the terms $\mathcal{M}_{21},\mathcal{M}_{22}$ in Appendix B, we have
\begin{eqnarray*}
	J_2+J_3+J_4
	&\leq&C\int_{0}^{T}\mathbb{E}[|X^{\varepsilon}_r-X^{\varepsilon}_{r(\delta)}|^2\mathbf{1}_{A_{N,r}}]\mathrm{d}r\cr
	&&+C \int_{0}^{T}\mathbb{E}[|X^{\varepsilon}_r-\hat{X}^{\varepsilon}_r|^2\mathbf{1}_{A_{N,r}}]\mathrm{d}r\cr
	&&+C \int_{0}^{T}\mathbb{E}[|\hat{X}^{\varepsilon}_r-\bar{X}_r|^2\mathbf{1}_{A_{N,r}}]\mathrm{d}r\cr
	&\leq& C\int_{0}^{T}\mathbb{E}[\|\hat{X}^{\varepsilon}-\bar{X}\|_{\alpha,r}^2\mathbf{1}_{A_{N,r}}]\mathrm{d}r\cr
	&&+C  \big(\varepsilon \delta^{-1} \delta^{2 \beta}e^{C\frac{\delta^2 }{\varepsilon^2}}+\delta^{2\beta'}\big).
\end{eqnarray*}

Using similar caculations as for the terms $\mathcal{N}_{33},\mathcal{N}_{34}$ and $\mathcal{N}_{36}$ in Appendix B,  we have
\begin{eqnarray*}
	J_5 &\leq & C\mathbb{E}\Big[\sup _{t \in[0, T]}\Big|\int_{0}^{t}S_{t-r}(g(X^{\varepsilon}_r)-g(\hat X^{\varepsilon}_r))dB^H_r \Big|^2\mathbf{1}_{A_{N,T}}\Big]\cr
	&&+C\mathbb{E}\Big[\int_{0}^T\Big(\int_{0}^t (t-s)^{-1-\alpha}\Big|\int_{s}^{t} S_{t-r} (g(X^{\varepsilon}_r)-g(\hat X^{\varepsilon}_r)) d B^H_r  \Big|\mathrm{d}s \Big)^2\mathrm{d}t \mathbf{1}_{A_{N,T}}  \Big]\cr
	&&+ C\mathbb{E}\Big[\int_{0}^T\Big(\int_{0}^t (t-s)^{-1-\alpha}\Big|\int_{0}^{s}(S_{t-r}-S_{s-r}) (g(X^{\varepsilon}_r)-g(\hat X^{\varepsilon}_r)) d B^H_r \Big|\mathrm{d}s\Big)^2\mathrm{d}t \mathbf{1}_{A_{N,T}} \Big]
	\cr
	&=:& J_{51}+J_{52}+J_{53},
\end{eqnarray*}
where
\begin{eqnarray*}
	J_{51}
	&\leq& C_{N} \int_{0}^{T}\mathbb{E}[\|X^{\varepsilon}-\hat{X}^{\varepsilon}\|_{\alpha,r}^2\mathbf{1}_{A_{N,r}}]\mathrm{d}r\cr
	&&+C_{N} \mathbb{E}\Big[\int_{0}^T \Big(\int_{0}^{t}\frac{|X^{\varepsilon}_t-\hat{X}^{\varepsilon}_t-X^{\varepsilon}_s+\hat{X}^{\varepsilon}_s|}{(t-s)^{1+\alpha}}\mathrm{d}s\Big)^2 \mathrm{d}t \mathbf{1}_{A_{N,T}}\Big]\cr
	&\leq& C_{N} \int_{0}^{T}\mathbb{E}[\|X^{\varepsilon}-\hat{X}^{\varepsilon}\|_{\alpha,r}^2\mathbf{1}_{A_{N,r}}]\mathrm{d}r+C_{N}\mathbb{E}\Big[\|X^{\varepsilon}-\hat{X}^{\varepsilon}\|_{\alpha,T}^2 \mathbf{1}_{A_{N,T}}\Big],
\end{eqnarray*}
and
\begin{eqnarray*}
	J_{52}+J_{53}
	\leq C_{N} \int_{0}^{T}\mathbb{E}[\|X^{\varepsilon}-\hat{X}^{\varepsilon}\|_{\alpha,r}^2\mathbf{1}_{A_{N,r}}]\mathrm{d}r.
\end{eqnarray*}
This yields that
\begin{eqnarray*}
	J_{5} \leq C_N \big(\varepsilon \delta^{-1} \delta^{2 \beta}e^{C\frac{\delta^2 }{\varepsilon^2}}+\delta^{2\beta'}\big).
\end{eqnarray*}

For $J_6$, using the same step as in the proof of $J_{5}$, we have
\begin{eqnarray*}
	J_6&\leq & C\mathbb{E}\Big[\sup _{t \in[0, T]}\Big|\int_{0}^{t}S_{t-r}(g(\hat X^{\varepsilon}_r)-g(\bar X _r))dB^H_r \Big|^2\mathbf{1}_{A_{N,T}}\Big]\cr
	&&+C\mathbb{E}\Big[\int_{0}^T\Big(\int_{0}^t (t-s)^{-1-\alpha}\Big|\int_{s}^{t} S_{t-r} (g(\hat X^{\varepsilon}_r)-g(\bar X _r)) d B^H_r  \Big|\mathrm{d}s\Big)^2\mathrm{d}t \mathbf{1}_{A_{N,T}} \Big]\cr
	&&+ C\mathbb{E}\Big[\int_{0}^T\Big(\int_{0}^t (t-s)^{-1-\alpha}\Big|\int_{0}^{s}(S_{t-r}-S_{s-r}) (g(\hat X^{\varepsilon}_r)-g(\bar X _r)) d B^H_r \Big| \mathrm{d}s\Big)^2\mathrm{d}t \mathbf{1}_{A_{N,T}}  \Big]
	\cr
	&=:& J_{61}+J_{62}+J_{63},
\end{eqnarray*}
where
\begin{eqnarray*}
	J_{61}
	\leq C_{N}  \int_{0}^{T}\mathbb{E}[\|\hat{X}^{\varepsilon}-\bar{X}\|_{\alpha,r}^2\mathbf{1}_{A_{N,r}}]\mathrm{d}r+C_{N}\mathbb{E}\Big[\int_{0}^T \Big(\int_{0}^{t}\frac{|\hat{X}^{\varepsilon}_t-\bar{X}_t-\hat{X}^{\varepsilon}_s+\bar{X}_s|}{(t-s)^{1+\alpha}}\mathrm{d}s\Big)^2  \mathrm{d}t \mathbf{1}_{A_{N,T}} \Big],
\end{eqnarray*}
and
\begin{eqnarray*}
	J_{62}+J_{63}
	\leq C_{N} \int_{0}^{T}\mathbb{E}[\|\hat{X}^{\varepsilon}-\bar{X}\|_{\alpha,r}^2\mathbf{1}_{A_{N,r}}]\mathrm{d}r.
\end{eqnarray*}

For $J_{61}$, from  Eq. (\ref{ave-eq}) and Eq. (\ref{xj}), using the fact that
\begin{eqnarray*}
	\hat{X}^{\varepsilon}_t-\bar{X}_t=\int_{0}^{t}S_{t-r}(b(X^{\varepsilon}_{r(\delta)},
	\hat{Y}^{\varepsilon}_r)-\bar b(\bar X_r)) \mathrm{d}r+\int_{0}^{t}S_{t-r}(g(X^{\varepsilon}_r)-\bar g(\bar X_r))\mathrm{d}B^H_r,
\end{eqnarray*}	
we have
\begin{eqnarray*}
	J_{61}
	&\leq& C_{N}  \int_{0}^{T}\mathbb{E}[\|\hat{X}^{\varepsilon}-\bar{X}\|_{\alpha,r}^2\mathbf{1}_{A_{N,r}}]\mathrm{d}r\cr
	&&+
	C_{N} \mathbb{E}\Big[\int_{0}^T \Big(\int_{0}^{t}{(t-s)^{-1-\alpha}}\Big|\int_{0}^{t}S_{t-r}(b(X^{\varepsilon}_{r(\delta)},
	\hat{Y}^{\varepsilon}_r)-\bar b(\bar X_r)) \mathrm{d}r\cr
	&&\qquad \qquad -\int_{0}^{s}S_{s-r}(b(X^{\varepsilon}_{r(\delta)},
	\hat{Y}^{\varepsilon}_r)-\bar b(\bar X_r)) \mathrm{d}r\Big|\mathrm{d}s\Big)^2 \mathrm{d}t \mathbf{1}_{A_{N,T}}\Big]\cr
	&&+C_{N} \mathbb{E}\Big[\int_{0}^T \Big(\int_{0}^{t}(t-s)^{-1-\alpha}\Big|\int_{0}^{t}S_{t-r}(g(\hat X^{\varepsilon}_r)-\bar g (\bar X_r)) \mathrm{d}r\cr
	&&\qquad \qquad -\int_{0}^{s}S_{s-r}(g(\hat X^{\varepsilon}_r)-\bar g(\bar X_r)) \mathrm{d}r\Big|\mathrm{d}s\Big)^2 \mathrm{d}t \mathbf{1}_{A_{N,T}}  \Big]\cr
	&&+C_{N} \mathbb{E}\Big[\int_{0}^T \Big(\int_{0}^{t}(t-s)^{-1-\alpha}\Big|\int_{0}^{t}S_{t-r}(g(X^{\varepsilon}_r)-\bar g (\hat X^{\varepsilon}_r)) \mathrm{d}r\cr
	&&\qquad \qquad -\int_{0}^{s}S_{s-r}(g(X^{\varepsilon}_r)-\bar g(\hat X^{\varepsilon}_r) \mathrm{d}r\Big|\mathrm{d}s\Big)^2 \mathrm{d}t \mathbf{1}_{A_{N,T}}  \Big]\cr
	&\leq& C_{N}  \int_{0}^{T}\mathbb{E}[\|\hat{X}^{\varepsilon}-\bar{X}\|_{\alpha,r}^2\mathbf{1}_{A_{N,r}}]\mathrm{d}r+\sum_{i=1}^{4}J_{i}+J_{52}+J_{53}+J_{62}+J_{63}\cr
	&\leq& C_{N} \int_{0}^{T}\mathbb{E}[\|\hat{X}^{\varepsilon}-\bar{X}\|_{\alpha,r}^2\mathbf{1}_{A_{N,r}}]\mathrm{d}r+C_{N}  \int_{0}^{T}\mathbb{E}[\|{X}^{\varepsilon}-\hat{X}^{\varepsilon}\|_{\alpha,r}^2\mathbf{1}_{A_{N,r}}]\mathrm{d}r
	\cr&&+C_N  \varepsilon \delta^{-1} \delta^{2 \beta}e^{C\frac{\delta^2 }{\varepsilon^2}}+ C_N  (\delta^{2\beta'}+\delta^{2\beta}).
\end{eqnarray*}	

This yields that
\begin{eqnarray*}
	J_6
	\leq C_{N} \int_{0}^{T}\mathbb{E}[\|\hat{X}^{\varepsilon}-\bar{X}\|_{\alpha,r}^2\mathbf{1}_{A_{N,r}}]\mathrm{d}r+C_N \big(\varepsilon \delta^{-1} \delta^{2 \beta}e^{C\frac{\delta^2 }{\varepsilon^2}}+  \delta^{2\beta'}\big).
\end{eqnarray*}	

Putting above results together, by Gronwall's lemma, we have
\begin{eqnarray*}
	\mathbb{E}[\|\hat{X}^{\varepsilon}-\bar{X}\|_{\alpha,T}^2\mathbf{1}_{A_{N,T}}] \leq C_N \big(\varepsilon \delta^{-1} \delta^{2 \beta}e^{C\frac{\delta^2 }{\varepsilon^2}}+  \delta^{2\beta'}\big).
\end{eqnarray*}

Finally, we have
\begin{eqnarray}\label{xhat-barx}
\mathbb{E}[\|\hat{X}^{\varepsilon}-\bar{X}\|_{\alpha,T}^2] \leq C_N \big(\varepsilon \delta^{-1} \delta^{2 \beta}e^{C\frac{\delta^2 }{\varepsilon^2}}+  \delta^{2\beta'}\big)+C\sqrt{N^{-1}\mathbb{E}[\Lambda_{\alpha,B^H}^{0,T}]}.
\end{eqnarray}
\para{Step 3:} Putting (\ref{x-hatx}) and (\ref{xhat-barx}) together, then, taking $\delta= \varepsilon \sqrt{-\ln \varepsilon}$, we have
\begin{eqnarray*}
	\lim\limits_{\varepsilon \rightarrow 0} \mathbb{E}[\|X^{\varepsilon}-\bar{X}\|_{\alpha,T}^2]= 0.
\end{eqnarray*}
This completes the proof of Theorem \ref{mainthm}.
\qed

\appendix\section{}
To make our paper self-contained, we recall the ergodicity for fast motion which was
introduced by Fu and Liu \cite{fu2011strong}. Consider the problem associated to the fast motion with frozen show component Eq. (\ref{invariant}).
In this section, we repalce the initial value $y$ defined in Eq.(\ref{invariant}) by $z$. Then, for any fixed $x\in V$ and initial value $z\in V$, Eq. (\ref{invariant}) has a unique strong solution (also a mild solution) which will be denoted by $(Y^{x,z}_t)_{t\geq 0}$. By energy equality \cite{fu2011strong} and conditions (A4)-(A5) and Poincar\'e inequality, one gets
\begin{eqnarray*}
	\mathbb{E}[|Y^{x,z}_t|^{2}] \leq |z|^{2}-(2\bar{\lambda}_1+2\beta_1-C_3)\mathbb{E}\Big[\int_{0}^{t}|Y^{x,z}_s|^{2}\mathrm{d}s\Big]+C_2(1+|x|^{2})t.
\end{eqnarray*}
By Gronwall's inequality, we have
$$\mathbb{E}[|Y^{x,z}_t|^{2}] \leq C(1+|z|^{2}e^{-(2\bar{\lambda}_1+2\beta_1-C_3)t}),$$
where $2\bar{\lambda}_1+2\beta_1-C_3>0$ owing to condition (A5) and $C>0$ is a constant.

Next, let $(Y^{x,z'}_t)_{t\geq 0}$ be a solution of Eq. (\ref{invariant}) with the initial value $Y_{0}^{x}=z'$, by the Poincar\'{e} inequality, conditions (A4) and (A5), we obtain
\begin{eqnarray*}
	\mathbb{E}[|Y^{x,z}_t-Y^{x,z'}_t|^{2}] &=& |z-z'|^{2}+2 \mathbb{E}\Big[\int_{0}^{t}\langle A(Y^{x,z}_s-Y^{x,z'}_s),Y^{x,z}_s-Y^{x,z'}_s \rangle \mathrm{d}s\Big]\cr
	&& +2 \mathbb{E}\Big[\int_{0}^{t}\langle F(x,Y^{x,z}_s)-F(x,Y^{x,z'}_s),Y^{x,z}_s-Y^{x,z'}_s \rangle  \mathrm{d}s\Big]\cr
	&& + \int_{0}^{t}\mathbb{E}[|G(x,Y^{x,z}_s)-G(x,Y^{x,z'}_s)|^{2}_{L_2(V)}] \mathrm{d}s\cr
	&\leq& |z-z'|^{2}- (2\bar{\lambda}_1-2\beta_{3}-C_{2})\int_{0}^{t}\mathbb{E}[|Y^{x,z}_s-Y^{x,z'}_s|^{2}]\mathrm{d}s
\end{eqnarray*}

Therefore, by Gronwall's inequality \cite[pp. 584]{Givon2007}, we have
\begin{eqnarray}\label{measure}
\mathbb{E}[|Y^{x,z}_t-Y^{x,z'}_t|^{2}]
&\leq |z-z'|^{2}e^{ -\eta t},
\end{eqnarray}
where $\eta=2\bar{\lambda}_1-2\beta_{3}-C_{2}>0.$

For any $x\in V$, denote by $(P^{x} _t)_{t\geq 0}$ the Markov semigroup associated to Eq. (\ref{invariant}) defined by
\begin{eqnarray*}
	P^{x} _t\Psi (z)= \mathbb{E}[\Psi (Y^{x,z}_t)], \quad  t \geq 0,  \quad z\in V,
\end{eqnarray*}
for any $\Psi \in \mathcal{B}_{b}(V)$, the space of bounded functions on $V$.
We
recall that a probability $\mu^{x}$ on $V$  is called an invariant measure for $(P^{x}_t)_{t\geq 0}$ if
$$\int_{V}P^{x}_t\Psi \mathrm{d} \mu^{x}=\int_{V}\Psi \mathrm{d} \mu^{x}, \qquad t\geq 0,$$
for any bounded function $\Psi \in \mathcal{B}_{b}(V)$. As in \cite{cerrai2009averaging}, it is possible to show the existence of the unique invariant measure $\mu^{x}$ for the semigroup  $(P^{x}_t)_{t\geq 0}$ which satisfies
$\int_{V}|z|\mu^{x}(\mathrm{d} z)\leq (1+|x|).$

Furthermore, according to Lipschitz assumption on $b$ and Eq. (\ref{measure}), we have
\begin{eqnarray}\label{yxz}
\Big|\mathbb{E}[b(x,Y^{x,y}_t)] -\int_{V}b(x,z)\mu^x(\mathrm{d} z) \Big|  &=&  \Big|\int_{V}\mathbb{E}[b(x,Y^{x,y}_t)-b(x,Y^{x,z}_t)]\mu^x(\mathrm{d} z) \Big|\cr
&\leq& C \int_{V}\mathbb{E} [|Y^{x,y}_t-Y^{x,z}_t|]\mu^x(\mathrm{d} z)\cr
&\leq& C e^{-\frac{1}{2}\eta t}\int_{V}  |y-z |\mu^x(\mathrm{d} z)\cr
& \leq& C e^{-\frac{1}{2}\eta t}(1+|x|+|y|),
\end{eqnarray}
where $C>0$ is a constant.

Suppose that {\rm (A1)-(A5)} hold. For any given value $x_1,x_2\in V,z\in V$, we have
\begin{eqnarray*}
	\mathbb{E}[|Y^{x_1,z}_t-Y^{x_2,z}_t|^{2}] &=& 2 \mathbb{E}\Big[\int_{0}^{t}\langle A(Y^{x_1,z}_s-Y^{x_2,z}_s),Y^{x_1,z}_s-Y^{x_2,z}_s \rangle \mathrm{d}s\Big]\cr
	&& +2 \mathbb{E}\Big[\int_{0}^{t}\langle F(x_1,Y^{x_1,z}_s)-F(x_2,Y^{x_2,z}_s),Y^{x_1,z}_s-Y^{x_2,z}_s \rangle  \mathrm{d}s\Big]\cr
	&& + \int_{0}^{t}\mathbb{E}[|G(x_1,Y^{x_1,z}_s)-G(x_2,Y^{x_2,z}_s)|^{2}_{L_2(V)}] \mathrm{d}s\cr
	&\leq&- (2\bar{\lambda}_1-2\beta_{3}-C_{2})\int_{0}^{t}\mathbb{E}[|Y^{x_1,z}_s-Y^{x_2,z}_s|^{2}]\mathrm{d}s+ C|x_1-x_2|^{2}.
\end{eqnarray*}

Therefore, by Gronwall's inequality \cite[pp. 584]{Givon2007}, we have
\begin{eqnarray}\label{measure2}
\mathbb{E}[|Y^{x_1,z}_t-Y^{x_2,z}_t|^{2}]
&\leq C |x_1-x_2|^{2}e^{ -\eta t}.
\end{eqnarray}

\section{}
In this section, we present the proofs of Lemma \ref{itofbm2}, Lemma \ref{betaHn1}, Lemma \ref{lemboundN},  Lemma \ref{lemboundNm} and Lemma \ref{claim} that have been
deferred from Section 3 and Section 4.
\para{The Proof of Lemma \ref{itofbm2}:} Note that $0\leq s<t \leq T$, by Remark \ref{itofbm} and Lemma \ref{inq-s} and Fubini's theorem, it follows
\begin{eqnarray*}
	\mathcal{K}_1(s,t)&\leq&
	\Lambda_{\alpha,B^H}^{0,t}\sup _{i} \int_{s}^{t} \Big(\frac{\|S_{t-r}\||g(u_r)e_{i}|}{(r-s)^{\alpha}}+\int_{s}^{r} \frac{\|S_{t-r}-S_{t-q}\||g(u_q)e_{i}|}{(r-q)^{1+\alpha}} d q \cr
	&&+\int_{s}^{r} \frac{\|S_{t-r}\|(g(u_r)-g(u_q)) e_{i}|}{(r-q)^{1+\alpha}} d q \Big) d r \cr
	&\leq& C \Lambda_{\alpha,B^H}^{0,t}\int_{s}^{t} \Big(\frac{1+|u_r|}{(r-s)^{\alpha}}+\int_{s}^r \frac{1+|u_q|}{(t-r)^{\beta}(r-q)^{1+\alpha-\beta}}\mathrm{d}q+\int_{s}^r \frac{|u_r-u_q|}{(r-q)^{1+\alpha}}\mathrm{d}q\Big)\mathrm{d}r\cr
	&\leq& C \Lambda_{\alpha,B^H}^{0,t} \int_{s}^{t}\Big([(r-s)^{-\alpha}+(t-r)^{-\alpha}](1+|u_r|)+\int_{s}^{r}\frac{|u_r-u_q|}{(r-q)^{1+\alpha}}\mathrm{d}q\Big)\mathrm{d}r.
\end{eqnarray*}

For $\mathcal{K}_2(0,s)$, taking $\alpha<\alpha'<1-\beta$, we obtain
\begin{eqnarray*}
	\mathcal{K}_2(0,s)  &\leq& \Lambda_{\alpha,B^H}^{0,t}\sup_{i}\int_{0}^{s}\frac{\|(S_{t-r}-S_{s-r})\||g(u_r)  e_{i}|}{r^\alpha} \mathrm{d}r\cr
	&&+ \Lambda_{\alpha,B^H}^{0,t}\sup_{i}\int_{0}^{s}\int_{0}^r\frac{\|(S_{t-r}-S_{s-r}-S_{t-q}+S_{s-q})\||g(u_q)e_i|}{(r-q)^{1+\alpha}}\mathrm{d}q \mathrm{d}r\cr
	&&+ \Lambda_{\alpha,B^H}^{0,t}\sup_{i}\int_{0}^{s}\int_{0}^r\frac{\|(S_{t-r}-S_{s-r})\||(g(u_r)-g(u_q))e_i|}{(r-q)^{1+\alpha}}\mathrm{d}q \mathrm{d}r\cr
	&\leq&C\Lambda_{\alpha,B^H}^{0,t}(t-s)^{\beta}\int_{0}^{s}\frac{1+|u_r|}{(s-r)^\beta r^\alpha} \mathrm{d}r\cr
	&&+C \Lambda_{\alpha,B^H}^{0,t}(t-s)^{\beta}\int_{0}^{s}\Big(\int_{0}^r\frac{(r-q)^{\alpha'}(1+|u_q|)}{(s-r)^{\alpha'+\beta}(r-q)^{1+\alpha}}\mathrm{d}q\Big) \mathrm{d}r\cr
	&&+C\Lambda_{\alpha,B^H}^{0,t} (t-s)^\beta \int_{0}^{s}(s-r)^{-\beta}\Big(\int_{0}^r\frac{|u_r-u_q|}{(r-q)^{1+\alpha}}\mathrm{d}q\Big)\mathrm{d}r\cr
	&\leq&   C\Lambda_{\alpha,B^H}^{0,t} (t-s)^\beta  \int_{0}^s  [(s-r)^{-\beta} r^{-\alpha}+(s-r)^{-\alpha-\beta}](1+|u_r|)\mathrm{d}r\cr
	&&+C\Lambda_{\alpha,B^H}^{0,t}  (t-s)^\beta  \int_{0}^s  (s-r)^{-\beta}\Big(\int_{0}^{r}\frac{|u_r-u_q|}{(r-q)^{1+\alpha}}\mathrm{d}q\Big)\mathrm{d}r.
\end{eqnarray*}

For $\mathcal{K}_3(s,t)$ and $\mathcal{K}_4(0,s)$,
we only need to replace the estimate for $g(u_r)$ which appeared in the above proofs of $\mathcal{K}_1(s,t)$ and $\mathcal{K}_2(0,s)$ by the corresponding estimate of $(g(u_{r})-g(v_{r}))$. By (\ref{linear}), we have
\begin{eqnarray*}
	\mathcal{K}_3(s,t)
	&\leq&
	\Lambda_{\alpha,B^H}^{0,t}\sup _{i} \int_{s}^{t} \Big(\frac{\|S_{t-r}\||(g(u_{r})-g(v_{r}))e_{i}|}{(r-s)^{\alpha}}\cr
	&&+\int_{s}^{r} \frac{\|S_{t-r}-S_{t-q}\| |(g(u_{q})-g(v_{q}))e_{i}|}{(r-q)^{1+\alpha}} d q \cr
	&&+\int_{s}^{r} \frac{\|S_{t-r}\||(g(u_{r})-g(v_{r})-(g(u_{q})-g(v_{q}))) e_{i}|}{(r-q)^{1+\alpha}} d q \Big) d r \cr
	&\leq& C \Lambda_{\alpha,B^H}^{0,t}
	\int_{s}^t[(r-s)^{-\alpha}+(t-r)^{-\alpha}]|u_r-v_r|
	\mathrm{d}r \cr
	&&+C \Lambda_{\alpha,B^H}^{0,t}
	\int_{s}^t|u_r-v_r|\Big(\int_{s}^{r}\frac{|u_r-u_q|+|v_r-v_q|}{(r-q)^{1+\alpha}}\mathrm{d}q\Big)\mathrm{d}r
	\cr
	&&+C \Lambda_{\alpha,B^H}^{0,t}
	\int_{s}^t\Big(\int_{s}^{r}\frac{|u_r-v_r-u_q+v_q|}{(r-q)^{1+\alpha}}\mathrm{d}q\Big)\mathrm{d}r,
\end{eqnarray*}
and
\begin{eqnarray*}
	\mathcal{K}_4(0,s) &\leq&
	\Lambda_{\alpha,B^H}^{0,t}\sup_{i}\int_{0}^{s}\frac{\|(S_{t-r}-S_{s-r})\||g(u_{r})-g(v_{r})e_{i}|}{r^\alpha} \mathrm{d}r\cr
	&&+ \Lambda_{\alpha,B^H}^{0,t}\sup_{i}\int_{0}^{s}\Big(\int_{0}^r\frac{\|(S_{t-r}-S_{s-r}-S_{t-q}+S_{s-q})\||(g(u_{q})-g(v_{q}))e_i|}{(r-q)^{1+\alpha}}\mathrm{d}q\Big) \mathrm{d}r\cr
	&&+ \Lambda_{\alpha,B^H}^{0,t}\sup_{i}\int_{0}^{s}\Big(\int_{0}^r\frac{\|(S_{t-r}-S_{s-r})\||(g(u_{r})-g(v_{r})-g(u_{q})+g(v_{q}))e_i|}{(r-q)^{1+\alpha}}\mathrm{d}q\Big) \mathrm{d}r\cr
	&\leq& C \Lambda_{\alpha,B^H}^{0,t} (t-s)^\beta \int_{0}^s [(s-r)^{-\beta} r^{-\alpha}+(s-r)^{-\alpha-\beta}]|u_r-v_r|\mathrm{d}r\cr
	&&+C \Lambda_{\alpha,B^H}^{0,t} (t-s)^\beta  \int_{0}^s (s-r)^{-\beta} |u_r-v_r| \Big(\int_{0}^{r}\frac{|u_r-u_q|+|v_r-v_q|}{(r-q)^{1+\alpha}}\mathrm{d}q\Big)\mathrm{d}r\cr
	&&+C \Lambda_{\alpha,B^H}^{0,t} (t-s)^\beta  \int_{0}^s (s-r)^{-\beta}\Big(\int_{0}^{r}\frac{|u_r-v_r-u_q+v_q|}{(r-q)^{1+\alpha}}\mathrm{d}q\Big)\mathrm{d}r.
\end{eqnarray*}
This completes the proof of Lemma \ref{itofbm2}.\qed
\para{The Proof of Lemma \ref{betaHn1}:}
Denote
$$A_{s,t}(\beta^{H,N}_i):=\Big( \frac{|\beta_{i}^{H,N}(t)-\beta_{i}^{H,N}(s)|}{(t-s)^{1-\alpha}}+\int_s^t\frac{|\beta_{i}^{H,N}(r)-\beta_{i}^{H,N}(s)|}{(r-s)^{2-\alpha}}\mathrm{d}r\Big).$$

Then, we have
\begin{eqnarray}\label{Ast}
A_{s,t}(\beta^{H,N}_i)& \leq & \Big( \frac{|\beta_{i}^{H,N}(t)-\beta_{i}^{H,N}(s)|}{(t-s)^{1-\alpha}}+\int_s^t\frac{|\beta_{i}^{H,N}(r)-\beta_{i}^{H,N}(s)|}{(r-s)^{2-\alpha}}\mathrm{d}r\Big)\mathbf{1}_{\{s \leq t \leq \tau_N\}}\cr
&&+\Big( \frac{|\beta_{i}^{H,N}(t)-\beta_{i}^{H,N}(s)|}{(t-s)^{1-\alpha}}+\int_s^t\frac{|\beta_{i}^{H,N}(r)-\beta_{i}^{H,N}(s)|}{(r-s)^{2-\alpha}}\mathrm{d}r\Big)\mathbf{1}_{\{\tau_N \leq s \leq t\}}\cr
&&+\Big( \frac{|\beta_{i}^{H,N}(t)-\beta_{i}^{H,N}(s)|}{(t-s)^{1-\alpha}}+\int_s^t\frac{|\beta_{i}^{H,N}(r)-\beta_{i}^{H,N}(s)|}{(r-s)^{2-\alpha}}\mathrm{d}r\Big)\mathbf{1}_{\{s \leq \tau_N \leq t\}}\cr
& \leq &
\Big( \frac{|\beta_{i}^{H}(t)-\beta_{i}^{H}(s)|}{(t-s)^{1-\alpha}}+\int_s^t\frac{|\beta_{i}^{H}(r)-\beta_{i}^{H}(s)|}{(r-s)^{2-\alpha}}\mathrm{d}r\Big)\mathbf{1}_{\{s \leq t \leq \tau_N\}}\cr
&&+\Big( \frac{|\beta_{i}^{H}(\tau_N)-\beta_{i}^{H}(\tau_N)|}{(t-s)^{1-\alpha}}+\int_s^t\frac{|\beta_{i}^{H}(\tau_N)-\beta_{i}^{H}(\tau_N)|}{(r-s)^{2-\alpha}}\mathrm{d}r\Big)\mathbf{1}_{\{\tau_N \leq s \leq t\}}\cr
&&+  \Big(\frac{|\beta_{i}^{H}(\tau_N)-\beta_{i}^{H}(s)|}{(\tau_N-s)^{1-\alpha}}\frac{(\tau_N-s)^{1-\alpha}}{(t-s)^{1-\alpha}}\Big)\mathbf{1}_{\{s \leq \tau_N \leq t\}}\cr
&&+\Big(\int_s^t\frac{|\beta_{i}^{H}(r \wedge \tau_N)-\beta_{i}^{H}(s)|}{(r\wedge \tau_N-s)^{2-\alpha}} \frac{(r\wedge \tau_N  -s)^{2-\alpha}}{(r -s)^{2-\alpha}}\mathrm{d}r\Big)\mathbf{1}_{\{s \leq \tau_N \leq t\}} \cr
&\leq &
C \Big( \frac{|\beta_{i}^{H}(t)-\beta_{i}^{H}(s)|}{(t-s)^{1-\alpha}}+\int_s^t\frac{|\beta_{i}^{H}(r)-\beta_{i}^{H}(s)|}{(r-s)^{2-\alpha}}\mathrm{d}r\Big)\mathbf{1}_{\{s \leq t \leq \tau_N\}}\cr
&&+ C  	\Big(\frac{|\beta_{i}^{H}(\tau_N)-\beta_{i}^{H}(s)|}{(\tau_N-s)^{1-\alpha}}+\int_s^t\frac{|\beta_{i}^{H}(r \wedge \tau_N)-\beta_{i}^{H}(s)|}{(r\wedge \tau_N-s)^{2-\alpha}} \mathrm{d}r\Big)\mathbf{1}_{\{s \leq \tau_N \leq t\}},
\end{eqnarray}
where $C$ is a constant which is independent of $N$ and $i$.

Thus, by (\ref{Ast}), we have
\begin{eqnarray*}
	\|\beta_{i}^{H,N}\|_{\alpha,0,T}= \sup_{0\le s<t\le T} A_{s,t}(\beta^{H,N}_i) \leq C \|\beta^H_i\|_{\alpha,0,\tau_N},
\end{eqnarray*}
almost surely. This completes the proof of Lemma \ref{betaHn1}.\qed

\para{The Proof of Lemma \ref{lemboundN}:}
We start with
\begin{eqnarray*}
	\mathbb{E} [\|u^{N,n}\|^2_{\alpha, T}]& \leq &
	C \mathbb{E}[\|S_{\cdot} u_{0}\|^2_{\alpha, T}]+C \mathbb{E}\Big[\Big\|\int_{0}^{\cdot} S_{\cdot-s} f(u^{N,n}_s) \mathrm{d} s\Big\|^2_{\alpha, T}\Big]\cr
	&&+C\mathbb{E}\Big[\Big\|\int_{0}^{\cdot} S_{\cdot-s} \sigma(u^{N,n}_s) \mathrm{d} W_s\Big\|^2_{\alpha, T}\Big]\cr
	&&+C\mathbb{E}\Big[\Big\|\int_{0}^{\cdot} S_{\cdot-s} g(u^{N,n}_s) \mathrm{d} B^{H,N,n}_s\Big\|^2_{\alpha, T}\Big]\cr
	&=:&\mathcal{M}_1+\mathcal{M}_2+\mathcal{M}_3+\mathcal{M}_4,
\end{eqnarray*}
where $B^{H,N,n}:=\sum_{i=1}^{\infty} \sqrt{\lambda_{i}} e_{i} \beta_{i}^{H,N,n}$.

Since $\alpha<\beta$, $u_0\in V_\beta$, it is easy to obtain
\begin{eqnarray*}
	\mathcal{M}_{1}&\leq &C\mathbb{E}\big[\sup _{t \in[0, T]} |S_t u_{0}|^2\big]+ C\mathbb{E}\Big[\int_{0}^{T}\Big(\int_{0}^{t} \frac{|(S_t-S_s)u_{0}|}{(t-s)^{1+\alpha}} \mathrm{d} s\Big)^2\mathrm{d}t\Big] \cr
	& \leq& C+ C\mathbb{E}\Big[\int_{0}^{T}\Big(\int_{0}^{t} |u_{0}|_{{\beta}}(t-s)^{\beta-1-\alpha} \mathrm{d} s\Big)^2\mathrm{d}t\Big] \cr
	&\leq& C.
\end{eqnarray*}

For $\mathcal{M}_2$, we have
\begin{eqnarray*}
	\mathcal{M}_2
	&\leq&C\mathbb{E}\Big[\sup _{t \in[0, T]} \Big|\int_{0}^{t} S_{t-s} f(u^{N,n}_s) \mathrm{d}s\Big|^2\Big]\cr
	&&+C\mathbb{E}\Big[\int_{0}^{T}\Big(  \int_{0}^{t}\frac{|\int_{0}^{t}S_{t-r} f(u^{N,n}_r) \mathrm{d} r-\int_{0}^{s}S_{s-r} f(u^{N,n}_r)\mathrm{d}r|}{(t-s)^{1+\alpha}}\mathrm{d}s\Big)^2\mathrm{d}t\Big]\cr
	&=:& \mathcal{M}_{21}+\mathcal{M}_{22}.
\end{eqnarray*}
By H\"older's inequality and the growth condition, we get
\begin{eqnarray*}
	\mathcal{M}_{21}
	&\leq& C \mathbb{E} \Big[\int_{0}^{T} \|S_{t-s}\|^2|f(u^{N,n}_s)|^2\mathrm{d}s\Big] \cr
	&\leq& C \int_{0}^{T}(1+\mathbb{E}[|u^{N,n}_s|^2])\mathrm{d}s.
\end{eqnarray*}
Next, for $\mathcal{M}_{22}$, by Fubini's theorem, Lemma \ref{inq-s} and H\"older's inequality, we get
\begin{eqnarray*}
	\mathcal{M}_{22}&\leq&C\mathbb{E}\Big[\int_{0}^{T}\Big(\int_{0}^{t}\frac{|\int_{s}^{t}S_{t-r} f(u^{N,n}_r) \mathrm{d}r|}{(t-s)^{1+\alpha}}\mathrm{d}s\Big)^2\mathrm{d}t\Big]\cr
	&&+C \mathbb{E}\Big[\int_{0}^{T}\Big(\int_{0}^{t}\frac{|\int_{0}^{s}(S_{t-r}-S_{s-r})f(u^{N,n}_r)\mathrm{d}r|}{(t-s)^{1+\alpha}}\mathrm{d}s\Big)^2\mathrm{d}t\Big]\cr
	&\leq&C\mathbb{E}\Big[\int_{0}^{T}\Big(\int_{0}^{t} \int_{0}^{r}(t-s)^{-1-\alpha}\mathrm{d}s|f(u^{N,n}_r)|\mathrm{d} r\Big)^2\mathrm{d}t\Big]\cr
	&&+ C\mathbb{E}\Big[\int_{0}^{T}\Big(\int_{0}^{t} \int_{r}^t (t-s)^{-1-\alpha+{\beta}}(s-r)^{-{\beta}}\mathrm{d}s |f(u^{N,n}_r)|\mathrm{d}r\Big)^2\mathrm{d}t\Big]\cr
	&\leq&C \mathbb{E}\Big[\int_{0}^{T}\Big(\int_{0}^{t} (t-r)^{-\alpha} |f(u^{N,n}_r)|d r \Big)^2\mathrm{d}t\Big]\cr
	&\leq&  C \int_{0}^{T}(1+\mathbb{E}[|u^{N,n}_s|^2])\mathrm{d}s.
\end{eqnarray*}

For $\mathcal{M}_3$, we get
\begin{eqnarray*}
	\mathcal{M}_3
	&\leq&C\mathbb{E}\Big[\sup _{t \in[0, T]}\Big|\int_{0}^{t} S_{t-s} \sigma(u^{N,n}_s) \mathrm{d} W_s\Big|^2\Big]\cr
	&&+C\mathbb{E}\Big[ \int_{0}^{T}\Big(  \int_{0}^{t}\frac{|\int_{0}^{t}S_{t-r} \sigma(u^{N,n}_r)\mathrm{d} W_r-\int_{0}^{s}S_{s-r} \sigma(u^{N,n}_r)\mathrm{d} W_r|}{(t-s)^{1+\alpha}}\mathrm{d}s\Big)^2\mathrm{d}t\Big]\cr
	&=:& \mathcal{M}_{31}+\mathcal{M}_{32},
\end{eqnarray*}
and then, by Burkholder-Davis-Gundy inequality,
\begin{eqnarray*}
	\mathcal{M}_{31}
	&\leq& C \mathbb{E} \Big[\int_{0}^{T} \|S_{t-s}\|^2|\sigma(u^{N,n}_s)|^2_{L_2(V)}\mathrm{d}s\Big]\cr
	&\leq& C \int_{0}^{T}(1+\mathbb{E}[|u^{N,n}_s|^2])\mathrm{d}s.
\end{eqnarray*}
Applying again H\"older's inequality, Burkholder-Davis-Gundy inequality, Fubini's theorem, Lemma \ref{inq-s} and Lemma \ref{inq-s}, we have
\begin{eqnarray*}
	\mathcal{M}_{32}&\leq&C\mathbb{E}\Big[\int_{0}^{T}\Big(\int_{0}^{t}\frac{|\int_{s}^{t}S_{t-r} \sigma(u^{N,n}_r) \mathrm{d}W_r|}{(t-s)^{1+\alpha}}\mathrm{d}s\Big)^2\mathrm{d}t\Big]\cr
	&&+C\mathbb{E}\Big[\int_{0}^{T}\Big(\int_{0}^{t}\frac{|\int_{0}^{s}(S_{t-r}-S_{s-r})\sigma(u^{N,n}_r)\mathrm{d}W_r|}{(t-s)^{1+\alpha}}\mathrm{d}s\Big)^2\mathrm{d}t\Big]\cr
	&\leq&C\mathbb{E}\Big[\int_{0}^{T}\int_{0}^{t}(t-s)^{-\frac{3}{2}-\alpha}\Big|\int_{s}^{t}S_{t-r} \sigma(u^{N,n}_r) \mathrm{d}W_r\Big|^2\mathrm{d}s\mathrm{d}t\Big]\cr
	&&+C\mathbb{E}\Big[\int_{0}^{T}\int_{0}^{t}(t-s)^{-\frac{3}{2}-\alpha}\Big|\int_{0}^{s}(S_{t-r}-S_{s-r})\sigma(u^{N,n}_r)\mathrm{d}W_r\Big|^2\mathrm{d}s\mathrm{d}t\Big]\cr
	&\leq&C\int_{0}^{T}\int_{0}^{t}(t-s)^{-\frac{3}{2}-\alpha}\int_{s}^{t}  \|S_{t-r}\| \mathbb{E}[|\sigma(u^{N,n}_r)|^2_{L_{2(V)}}]d rds\mathrm{d}t\cr
	&&+C\int_{0}^{T}\int_{0}^{t} \int_{0}^{s}(t-s)^{-\frac{3}{2}-\alpha}\|S_{t-r}-S_{s-r}\|^2 \mathbb{E}[|\sigma(u^{N,n}_r)|^2_{L_{2(V)}}]\mathrm{d}r \mathrm{d}s \mathrm{d}t\cr
	&\leq&C\int_{0}^{T}\int_{0}^{t}\int_{0}^{r}(t-s)^{-\frac{3}{2}-\alpha} \mathrm{d}s (1+\mathbb{E}[|u^{N,n}_r|^2])\mathrm{d} r\mathrm{d}t\cr
	&&+C\int_{0}^{T} \int_{0}^{t}\int_{r}^{t}(t-s)^{-\frac{3}{2}-\alpha+2\beta'} (s-r)^{-2\beta'}\mathrm{d}s(1+\mathbb{E}[|u^{N,n}_r|^2])\mathrm{d}r\mathrm{d}t\cr
	&\leq& C \int_{0}^{T}(1+\mathbb{E}[\|u^{N,n}\|^2_{\alpha, t}])\mathrm{d}t.
\end{eqnarray*}

To proceed, for $\mathcal{M}_4$, we obtain
\begin{eqnarray*}
	\mathcal{M}_4
	&\leq& C \mathbb{E}\Big[\int_{0}^{T}\Big(  \int_{0}^{t}\frac{|\int_{0}^{t}S_{t-r} g(u^{N,n}_r)\mathrm{d}B^{H,N,n}_r-\int_{0}^{s}S_{s-r} g(u^{N,n}_r)\mathrm{d}B^{H,N,n}_r|}{(t-s)^{1+\alpha}}\mathrm{d}s\Big)^2\mathrm{d}t\Big]\cr
	&&+C\mathbb{E}\Big[\sup _{t \in[0, T]} \Big|\int_{0}^{t} S_{t-r} g(u^{N,n}_r) \mathrm{d}B^{H,N,n}_r\Big|^2\Big]\cr
	&=:& \mathcal{M}_{41}+\mathcal{M}_{42}.
\end{eqnarray*}

For $\mathcal{M}_{41}$, by (\ref{betaNn-Nn}), Lemma \ref{itofbm2} and H\"older's inequality, taking $\alpha<\alpha'<1-\beta$, we have
\begin{eqnarray*}
	\mathcal{M}_{41}&\leq& C \mathbb{E}\Big[\int_{0}^{T}\Big(\int_{0}^t (t-s)^{-1-\alpha}\Big|\int_{s}^{t} S_{t-r} g(u^{N,n}_r) \mathrm{d}B^{H,N,n}_r  \Big| \mathrm{d}s \Big)^2\mathrm{d}t\Big] \cr
	&&+C\mathbb{E}\Big[\int_{0}^{T}\Big(\int_{0}^t (t-s)^{-1-\alpha}\Big|\int_{0}^{s}(S_{t-r}-S_{s-r}) g(u^{N,n}_r) \mathrm{d}B^{H,N,n}_r  \Big| \mathrm{d}s\Big)^2\mathrm{d}t\Big] \cr
	&\leq&
	C\mathbb{E}\Big[\int_{0}^{T}\Big(\int_{0}^t \frac{\Lambda_{\alpha,B^{H,N,n}}^{0,t}\int_{s}^{t}[(r-s)^{-\alpha}+(t-r)^{-\alpha}](1+|u^{N,n}_r|)\mathrm{d}r}{(t-s)^{1+\alpha}}
	\mathrm{d}s\Big)^2\mathrm{d}t\Big]\cr
	&&+
	C\mathbb{E}\Big[\int_{0}^{T}\Big(\int_{0}^t \frac{\Lambda_{\alpha,B^{H,N,n}}^{0,t}\int_{s}^{t}\big(\int_{s}^{r}\frac{|u^{N,n}_r-u^{N,n}_q|}{(r-q)^{1+\alpha}}\mathrm{d}q\big)\mathrm{d}r}{(t-s)^{1+\alpha}}
	\mathrm{d}s\Big)^2\mathrm{d}t\Big]\cr
	&&+C\mathbb{E}\Big[\int_{0}^{T}\Big(\int_{0}^t\frac{\Lambda_{\alpha,B^{H,N,n}}^{0,t}\int_{0}^s  \frac{r^{-\alpha}+(s-r)^{-\alpha}}{(s-r)^{\beta}}(1+|u^{N,n}_r|)\mathrm{d}r}{(t-s)^{1+\alpha-\beta}}  \mathrm{d}s\Big)^2\mathrm{d}t\Big]\cr
	&&+C\mathbb{E}\Big[\int_{0}^{T}\Big(\int_{0}^t\frac{\Lambda_{\alpha,B^{H,N,n}}^{0,t}\int_{0}^s  (s-r)^{-\beta}(\int_{s}^{r}\frac{|u^{N,n}_r-u^{N,n}_q|}{(r-q)^{1+\alpha}}\mathrm{d}q)\mathrm{d}r}{(t-s)^{1+\alpha-\beta}}\mathrm{d}s\Big)^2\mathrm{d}t\Big]\cr
	&\leq&
	C_N\mathbb{E}\Big[\int_{0}^{T}\Big(\int_{0}^t\int_{0}^{r} \frac{(r-s)^{-\alpha}+(t-r)^{-\alpha}}{(t-s)^{1+\alpha}}\mathrm{d}s(1+|u^{N,n}_r|)
	\mathrm{d}r \Big)^2\mathrm{d}t\Big]\cr
	&&+
	C_N\mathbb{E}\Big[\int_{0}^{T}\Big(\int_{0}^t\int_{0}^{r} (t-s)^{-1-\alpha}\mathrm{d}s\Big(\int_{0}^{r}\frac{|u^{N,n}_r-u^{N,n}_q|}{(r-q)^{1+\alpha}}\mathrm{d}q\Big)\mathrm{d}r
	\Big)^2\mathrm{d}t\Big]\cr
	&&+C_N\mathbb{E}\Big[\int_{0}^{T}\Big(\int_{0}^t \int_{r}^t \frac{(s-r)^{-\beta} r^{-\alpha}+(s-r)^{-\alpha-\beta}}{(t-s)^{1+\alpha-\beta}}\mathrm{d}s(1+|u^{N,n}_r|)\mathrm{d}r\Big)^2\mathrm{d}t\Big]\cr
	&&+C_N\mathbb{E}\Big[\int_{0}^{T}\Big(\int_{0}^t \int_{r}^t \frac{(s-r)^{-\beta}}{(t-s)^{1+\alpha-\beta}}\mathrm{d}s\Big(\int_{0}^{r}\frac{|u^{N,n}_r-u^{N,n}_q|}{(r-q)^{1+\alpha}}\mathrm{d}q\Big)\mathrm{d}r\Big)^2\mathrm{d}t\Big]\cr
	&\leq& C_{N} \mathbb{E}\Big[\int_{0}^{T} \int_{0}^t (t-r)^{-\alpha}[ r^{-\alpha}+(t-r)^{-\alpha}](1+|u^{N,n}_r|^2)\mathrm{d}r\mathrm{d}t\Big]\cr
	&&+C_{N} \mathbb{E}\Big[\int_{0}^{T}\int_{0}^t \Big(\int_{0}^{r}\frac{|u^{N,n}_r-u^{N,n}_q|}{(r-q)^{1+\alpha}}\mathrm{d}q\Big)^2\mathrm{d}r\mathrm{d}t\Big]\cr
	&\leq& C_{N} \int_{0}^{T} (1+\mathbb{E}[\|u^{N,n}\|_{\alpha,t}^2])\mathrm{d}t.
\end{eqnarray*}

For $\mathcal{M}_{42}$, by Lemma \ref{itofbm2} and H\"older's inequality, we have
\begin{eqnarray*}
	\mathcal{M}_{42}
	&\leq&
	C\mathbb{E}\Big[\sup _{t \in[0, T]} \Big( \Lambda_{\alpha,B^{H,N,n}}^{0,t}\int_{0}^t(r^{-\alpha}+(t-r)^{-\alpha})(1+|u^{N,n}_r|)\mathrm{d}r
	\Big)^2\Big]\cr
	&&+
	C\mathbb{E}\Big[\sup _{t \in[0, T]} \Big( \Lambda_{\alpha,B^{H,N,n}}^{0,t}\int_{0}^{t}\int_{0}^{r}\frac{|u^{N,n}_r-u^{N,n}_q|}{(r-q)^{1+\alpha}}\mathrm{d}q\mathrm{d}r
	\Big)^2\Big]\cr
	&\leq & C_{N} \int_{0}^{T}(1+\mathbb{E}[\|u^{N,n}\|^2_{\alpha, t}])\mathrm{d}t\cr
	&&+ C_{N}  \mathbb{E}\Big[\int_{0}^{T} \Big(\int_{0}^{t}\frac{|u^{N,n}_t-u^{N,n}_s|}{(t-s)^{1+\alpha}}\mathrm{d}s\Big)^2\mathrm{d}t\Big].
\end{eqnarray*}

The obtained results of $\mathcal{M}_1,\mathcal{M}_{22},\mathcal{M}_{32}$ and $\mathcal{M}_{41}$ yield that
\begin{eqnarray*}
	\mathbb{E}\Big[\int_{0}^{T} \Big(\int_{0}^{t}\frac{|u^{N,n}_t-u^{N,n}_s|}{(t-s)^{1+\alpha}}\mathrm{d}s\Big)^2\mathrm{d}t\Big]\leq C+ C_{N} \int_{0}^{T} (1+\mathbb{E}[\|u^{N,n}\|_{\alpha,t}^2])\mathrm{d}t.
\end{eqnarray*}
Thus, we have $$\mathcal{M}_{42}\leq C+C_{N} \int_{0}^{T} (1+\mathbb{E}[\|u^{N,n}\|_{\alpha,t}^2])\mathrm{d}t.$$

Finally, by Gronwall's lemma, we have
\begin{eqnarray*}
	\mathbb{E}[\|u^{N,n}\|^2_{\alpha, T}] \leq C_{N}.
\end{eqnarray*}
This completes the proof of Lemma \ref{lemboundN}.\qed

\para{The Proof of Lemma \ref{lemboundNm}:}
From (\ref{map1}), we have
\begin{eqnarray*}
	&&\mathbb{E}\Big[\|u^{N,n}-u^{N,m}\|^2_{\alpha,T}\mathbf{1}_{D_{T}^{N,R}}\Big] \cr
	&\leq& C \mathbb{E}\Big[\Big\|\int_{0}^{\cdot} S_{\cdot-r} (f(u^{N,n}_r)- f(u^{N,m}_r)) \mathrm{d} r\Big\|^2_{\alpha, T}\mathbf{1}_{D_{T}^{N,R}}\Big] \cr
	&&+C\mathbb{E}\Big[\Big\|\int_{0}^{\cdot} S_{\cdot-r} (\sigma(u^{N,n}_r)-\sigma(u^{N,m}_r)) \mathrm{d}W_r\Big\|^2_{\alpha, T}\mathbf{1}_{D_{T}^{N,R}}\Big] \cr
	&&+C\mathbb{E}\Big[\Big\|\int_{0}^{\cdot} S_{\cdot-r} g(u^{N,n}_r) \mathrm{d}B^{H,N,n}_r-\int_{0}^{\cdot} S_{\cdot-r} g(u^{N,m}_r) \mathrm{d}B^{H,N,m}_r\Big\|^2_{\alpha, T}\mathbf{1}_{D_{T}^{N,R}}\Big] \cr
	&=:&\mathcal{N}_1+\mathcal{N}_2+\mathcal{N}_3.
\end{eqnarray*}

The terms $\mathcal{N}_1$ and $\mathcal{N}_2$ can be estimated in the same way as the terms $\mathcal{M}_2$ and $\mathcal{M}_3$ in the proof of Lemma \ref{lemboundN}, using Lipschitz condition instead of the growth condition. This leads to the bounds
\begin{eqnarray*}
	\mathcal{N}_1+\mathcal{N}_2 \leq  C\int_{0}^{T}\mathbb{E}\big[\|u^{N,n}-u^{N,m}\|_{\alpha,t}^2\mathbf{1}_{D_{t}^{N,R}}\big]\mathrm{d}t.
\end{eqnarray*}

For $\mathcal{N}_3$, we have
\begin{eqnarray*}
	\mathcal{N}_3
	&\leq&C \mathbb{E}\Big[\int_{0}^{T}\Big(\int_{0}^{t}\Big|\int_{0}^{t} S_{t-r} g(u^{N,n}_r) \mathrm{d}B^{H,N,n}_r-\int_{0}^{t} S_{t-r} g(u^{N,m}_r) \mathrm{d}B^{H,N,m}_r\cr
	&&-\int_{0}^{s} S_{s-r} g(u^{N,n}_r) \mathrm{d}B^{H,N,n}_r+\int_{0}^{s} S_{s-r} g(u^{N,m}_r) \mathrm{d}B^{H,N,m}_r\Big|(t-s)^{-1-\alpha}\mathrm{d}s\Big)^2\mathbf{1}_{D_{t}^{N,R}}\mathrm{d}t\Big]\cr
	&&+
	C\mathbb{E}\Big[\sup _{t \in[0, T]} \Big(\Big|\int_{0}^{t} S_{t-r} g(u^{N,n}_r) \mathrm{d}B^{H,N,n}_r-\int_{0}^{t} S_{t-r} g(u^{N,m}_r) \mathrm{d}B^{H,N,m}_r\Big|^2\mathbf{1}_{D_{t}^{N,R}}\Big)\Big]\cr
	&\leq&C \mathbb{E}\Big[\int_{0}^{T}\Big(\int_{0}^{t}\frac{\big|\int_{s}^{t} S_{t-r} g(u^{N,n}_r) \mathrm{d}(B^{H,N,n}_r-B^{H,N,m}_r)\big|}{(t-s)^{1+\alpha}}\mathrm{d}s\Big)^2\mathbf{1}_{D_{t}^{N,R}}\mathrm{d}t\Big]\cr
	&&+C \mathbb{E}\Big[\int_{0}^{T}\Big( \int_{0}^{t}\frac{\big|\int_{0}^{s}(S_{t-r}-S_{s-r}) g(u^{N,n}_r) \mathrm{d}(B^{H,N,n}_r-B^{H,N,m}_r)\big|}{(t-s)^{1+\alpha}}\mathrm{d}s\Big)^2\mathbf{1}_{D_{t}^{N,R}}\mathrm{d}t\Big]\cr
	&&+C \mathbb{E}\Big[\int_{0}^{T}\Big(\int_{0}^{t}\frac{\big|\int_{s}^{t} S_{t-r} (g(u^{N,n}_r)-g(u^{N,m}_r)) \mathrm{d}B^{H,N,m}_r\big|}{(t-s)^{1+\alpha}}\mathrm{d}s\Big)^2\mathbf{1}_{D_{t}^{N,R}}\mathrm{d}t\Big]\cr
	&&+C \mathbb{E}\Big[\int_{0}^{T}\Big(\int_{0}^{t}\frac{\big|\int_{0}^{s} (S_{t-r}-S_{s-r}) (g(u^{N,n}_r)-g(u^{N,m}_r)) \mathrm{d}B^{H,N,m}_r\big|}{(t-s)^{1+\alpha}}\mathrm{d}s\Big)^2\mathbf{1}_{D_{t}^{N,R}}\mathrm{d}t\Big]\cr
	&&+C\mathbb{E}\Big[\sup _{t \in[0, T]} \Big|\int_{0}^{t} S_{t-r} g(u^{N,n}_r) \mathrm{d}(B^{H,N,n}_r-B^{H,N,m}_r)\Big|^2\mathbf{1}_{D_{T}^{N,R}}\Big]\cr
	&&+C\mathbb{E}\Big[\sup _{t \in[0, T]} \Big|\int_{0}^{t} S_{t-r} (g(u^{N,n}_r)-g(u^{N,m}_r)) \mathrm{d}B^{H,N,m}_r\Big|^2\mathbf{1}_{D_{T}^{N,R}}\Big]\cr
	&=:& \sum_{i=1}^{6}\mathcal{N}_{3i}.
\end{eqnarray*}

In the same way as for the term $\mathcal{M}_4$ in the proof of Lemma \ref{lemboundN}, we have
\begin{eqnarray*}
	\mathcal{N}_{31}+\mathcal{N}_{32} &\leq&
	C\mathbb{E}\Big[\int_{0}^{T}\big(\Lambda^{0,T}_{\alpha,B^{H,N,n,m}}\big)^2\Big(\int_{0}^t\frac{ r^{-\alpha}+(t-r)^{-\alpha}}{(t-r)^{\alpha}}(1+|u^{N,n}_r|)\mathrm{d}r\mathbf{1}_{D_{t}^{N,R}}\Big)^2\mathrm{d}t\Big]\cr
	&&+C\mathbb{E}\Big[\int_{0}^{T}\big(\Lambda^{0,T}_{\alpha,B^{H,N,n,m}}\big)^2\Big(\int_{0}^t \frac{\int_{0}^{r}\frac{|u^{N,n}_r-u^{N,n}_q|}{(r-q)^{1+\alpha}}\mathrm{d}q}{(t-r)^{\alpha}}\mathrm{d}r\mathbf{1}_{D_{t}^{N,R}}\Big)^2\mathrm{d}t\Big]\cr
	&\leq&C\mathbb{E}\Big[\int_{0}^{T}\big(\Lambda^{0,T}_{\alpha,B^{H,N,n,m}}\big)^2(1+\|u^{N,n}\|_{\alpha,t}^2)\mathbf{1}_{D_{t}^{N,R}}\mathrm{d}t\Big]\cr
	&\leq& C_{N,R, T, \alpha,\beta}\mathbb{E}\big[\big(\Lambda^{0,T}_{\alpha,B^{H,N,n,m}}\big)^2\big],
\end{eqnarray*}
where $\Lambda^{0,T}_{\alpha,B^{H,N,n,m}}:=\sum_{i=1}^{\infty}\sqrt{\lambda_{i}}\|\beta_{i}^{H,N,n}-\beta_{i}^{H,N,m }\|_{\alpha,0,T}.$

Next, for $\mathcal{N}_{33},$ and $\mathcal{N}_{34}$, by Lemma \ref{itofbm2}, we have
\begin{eqnarray*}
	&&\mathcal{N}_{33}+\mathcal{N}_{34}\cr
	&\leq&C \mathbb{E}\Big[\int_{0}^{T}\Big(\int_{0}^{t}\frac{\big|\int_{s}^{t} S_{t-r} (g(u^{N,n}_r)-g(u^{N,m}_r)) \mathrm{d}B^{H,N,m}_r\big|}{(t-s)^{1+\alpha}}\mathrm{d}s\mathbf{1}_{D_{t}^{N,R}}\Big)^2\mathrm{d}t\Big]\cr
	&&+C \mathbb{E}\Big[\int_{0}^{T}\Big(\int_{0}^{t}\frac{\big| \int_{0}^{s} (S_{t-r}-S_{s-r}) (g(u^{N,n}_r)-g(u^{N,m}_r)) \mathrm{d}B^{H,N,m}_r\big|}{(t-s)^{1+\alpha}}\mathrm{d}s\mathbf{1}_{D_{t}^{N,R}}\Big)^2\mathrm{d}t\Big]\cr
	&\leq& C_{N}\mathbb{E}\Big[\int_{0}^{T} \int_{0}^t (t-r)^{-\alpha}[ r^{-\alpha}+(t-r)^{-\alpha}]|u^{N,n}_r-u^{N,m}_r|^2\mathrm{d}r\mathbf{1}_{D_{t}^{N,R}}\mathrm{d}t\Big]\cr
	&&+C_{N}\mathbb{E}\Big[\int_{0}^{T}\Big(\int_{0}^t \frac{\int_{0}^{r}|u^{N,n}_r-u^{N,m}_r-u^{N,n}_q+u^{N,m}_q|{(r-q)^{-1-\alpha}}\mathrm{d}q}{(t-r)^\alpha}\mathrm{d}r \mathbf{1}_{D_{t}^{N,R}} \Big)^2 \mathrm{d}t\Big]\cr
	&&+C_{N}\mathbb{E}\Big[\int_{0}^{T}\Big(\int_{0}^t \frac{|u^{N,n}_r-u^{N,m}_r|}{ (t-r)^{\alpha}}\Big(\int_{0}^{r}\frac{|u^{N,n}_r-u^{N,n}_q|}{(r-q)^{1+\alpha}}\mathrm{d}q\Big)\mathrm{d}r\mathbf{1}_{D_{t}^{N,R}}\Big)^2\mathrm{d}t\Big]\cr
	&&+C_{N}\mathbb{E}\Big[\int_{0}^{T}\Big(\int_{0}^t \frac{ |u^{N,n}_r-u^{N,m}_r|}{ (t-r)^{\alpha}}\Big(\int_{0}^{r}\frac{|u^{N,m}_r-u^{N,m}_q|}{(r-q)^{1+\alpha}}\mathrm{d}q\Big)\mathrm{d}r\mathbf{1}_{D_{t}^{N,R}}\Big)^2\mathrm{d}t\Big]\cr
	&\leq& C_{N}\mathbb{E}\Big[\int_{0}^{T} \int_{0}^t (t-r)^{-\alpha}[ r^{-\alpha}+(t-r)^{-\alpha}]|u^{N,n}_r-u^{N,m}_r|^2\mathrm{d}r\mathbf{1}_{D_{t}^{N,R}}\mathrm{d}t\Big]\cr
	&&+C_N\mathbb{E}\Big[\int_{0}^{T}\int_{0}^t \Big(\int_{0}^{r}\frac{|u^{N,n}_r-u^{N,m}_r-u^{N,n}_q+u^{N,m}_q|}{(r-q)^{1+\alpha}}\mathrm{d}q\Big)^2\mathrm{d}r\mathbf{1}_{D_{t}^{N,R}}\mathrm{d}t\Big]\cr
	&&+C_{N}\mathbb{E}\Big[\int_{0}^{T}\int_{0}^t |u^{N,n}_r-u^{N,m}_r|^2\Big(\int_{0}^{r}\frac{|u^{N,n}_r-u^{N,n}_q|}{(r-q)^{1+\alpha}}\mathrm{d}q\Big)^2\mathrm{d}r\mathbf{1}_{D_{t}^{N,R}}\mathrm{d}t\Big]\cr
	&&+C_{N}\mathbb{E}\Big[\int_{0}^{T}\int_{0}^t |u^{N,n}_r-u^{N,m}_r|^2\Big(\int_{0}^{r}\frac{|u^{N,m}_r-u^{N,m}_q|}{(r-q)^{1+\alpha}}\mathrm{d}q\Big)^2\mathrm{d}r\mathbf{1}_{D_{t}^{N,R}}\mathrm{d}t\Big]\cr
	&\leq& C_{N,R}\int_{0}^{T} \mathbb{E}\big[\|u^{N,n}-u^{N,m}\|^2_{\alpha,t}\mathbf{1}_{D_{t}^{N,R}}\big]\mathrm{d}t.
\end{eqnarray*}

Next, for $\mathcal{N}_{35},$ we have
\begin{eqnarray*}
	\mathcal{N}_{35}
	&\leq&
	C\mathbb{E} \Big[\sup _{t \in[0, T]} \Big(\Lambda^{0,T}_{\alpha,B^{H,N,n,m}}\int_{0}^t(r^{-\alpha}+(t-r)^{-\alpha})(1+|u^{N,n}_r|)\mathrm{d}r\mathbf{1}_{D_{t}^{N,R}}
	\Big)^2\Big]\cr
	&&+
	C\mathbb{E}\Big[\sup _{t \in[0, T]} \Big(\Lambda^{0,T}_{\alpha,B^{H,N,n,m}}\int_{0}^{t}\Big(\int_{0}^{r}\frac{|u^{N,n}_r-u^{N,n}_q|}{(r-q)^{1+\alpha}}\mathrm{d}q\Big)\mathrm{d}r\mathbf{1}_{D_{t}^{N,R}}
	\Big)^2\Big]\cr
	&\leq & C\mathbb{E} \Big[ \int_{0}^{T}\big(\Lambda^{0,T}_{\alpha,B^{H,N,n,m}}\big)^2(1+\|u^{N,n}\|^2_{\alpha, t})\mathbf{1}_{D_{t}^{N,R}}\mathrm{d}t\Big]\cr
	&&+ C \mathbb{E} \Big[\big(\Lambda^{0,T}_{\alpha,B^{H,N,n,m}}\big)^2 \int_{0}^{T}\Big(\int_{0}^{t}\frac{|u^{N,n}_t-u^{N,n}_s|}{(t-s)^{1+\alpha}}\mathrm{d}s\Big)^2\mathbf{1}_{D_{t}^{N,R}}\mathrm{d}t\Big]\cr
	&\leq & C_{R}\mathbb{E} \big[\big(\Lambda^{0,T}_{\alpha,B^{H,N,n,m}}\big)^2\big].
\end{eqnarray*}

In the same way as for the term $\mathcal{M}_{42}$ in the proof of Lemma \ref{lemboundN}, we have
\begin{eqnarray*}
	\mathcal{N}_{36}&\leq&
	C_{N} \mathbb{E}\Big[\sup _{t \in[0, T]} \Big(\int_{0}^t(r^{-\alpha}+(t-r)^{-\alpha})|u^{N,n}_r-u^{N,m}_r|\mathrm{d}r
	\Big)^2\Big]\cr
	&&+C_{N}  \mathbb{E}\Big[\sup_{t\in[0,T]} \Big(\int_{0}^t\Big(\int_{0}^{r}\frac{|u^{N,n}_r-u^{N,m}_r-u^{N,n}_q+u^{N,m}_q|}{(r-q)^{1+\alpha}}\mathrm{d}q\Big)\mathrm{d}r\mathbf{1}_{D_{t}^{N,R}}\Big)^2\Big]\cr
	&&+C_{N}  \mathbb{E}\Big[\sup_{t\in[0,T]}\Big(\int_{0}^t |u^{N,n}_r-u^{N,m}_r|\Big(\int_{0}^{r}\frac{|u^{N,n}_r-u^{N,n}_q|}{(r-q)^{1+\alpha}}\mathrm{d}q\Big)\mathrm{d}r\mathbf{1}_{D_{t}^{N,R}}\Big)^2\Big]\cr
	&&+C_{N}  \mathbb{E}\Big[\sup_{t\in[0,T]}\Big(\int_{0}^t  |u^{N,n}_r-u^{N,m}_r|\Big(\int_{0}^{r}\frac{|u^{N,m}_r-u^{N,m}_q|}{(r-q)^{1+\alpha}}\mathrm{d}q\Big)\mathrm{d}r\mathbf{1}_{D_{t}^{N,R}}\Big)^2\Big]\cr
	&\leq& C_{N,R}  \int_{0}^{T} \mathbb{E}\big[\|u^{N,n}-u^{N,m}\|^2_{\alpha,t}\mathbf{1}_{D_{t}^{N,R}}\big]\mathrm{d}t\cr
	&&+C_{N,R} \mathbb{E}\Big[\int_{0}^T \Big(\int_{0}^{t}\frac{|u^{N,n}_t-u^{N,m}_t-u^{N,n}_s+u^{N,m}_s|}{(t-s)^{1+\alpha}}\mathrm{d}s\mathbf{1}_{D_{t}^{N,R}}\Big)^2\mathrm{d}t\Big]\cr
	&\leq& C_{N,R}  \int_{0}^{T} \mathbb{E}\big[\|u^{N,n}-u^{N,m}\|^2_{\alpha,t}\mathbf{1}_{D_{t}^{N,R}}\big]\mathrm{d}t.
\end{eqnarray*}

Finally, by Gronwall's lemma, we have
\begin{eqnarray*}
	\mathbb{E}\Big[\|u^{N,n}-u^{N,m}\|^2_{\alpha,T}\mathbf{1}_{D_{T}^{N,R}}\Big] \leq C_{N,R} \mathbb{E}\Big[\Big(\sum_{i=1}^{\infty}\sqrt{\lambda_{i}}\|\beta_{i}^{H,N,n}-\beta_{i}^{H,N,m }\|_{\alpha,0,T}\Big)^2\Big].
\end{eqnarray*}
This completes the proof of Lemma \ref{lemboundNm}.\qed

\para{The Proof of Lemma \ref{claim}:}
Let $X^{\varepsilon}_{k\delta},\hat{Y}^{\varepsilon}_{k\delta}, Y^{X^{\varepsilon}_{k\delta},\hat{Y}^{\varepsilon}_{k\delta}}$ and $\bar{W}$ be as Eq. (\ref{construct2}) and let $\mathbb{Q}^{y}$ denote the probability law of the diffusion process $(Y^{x}_t)_{t\geq 0}$ which is governed by the SPDE
\begin{eqnarray}\label{11}
\mathrm{d} Y^{x}_t =(AY^{x}_t+F(x,Y^{x}_t))\mathrm{d}t+G(x,Y^{x}_t) \mathrm{d} \bar{W}_t,
\end{eqnarray}
with initial value $Y^{x}_0=y$ and we denote that solution by $(Y^{x,y}_t)_{t\geq 0}$ (Without ambiguity, one can still use the notation $Y$ to denote the solution of Eq.(\ref{11})). The expectation with respect to $\mathbb{Q}^{y}$ is denoted by $\mathbb{E}^{y}$. Hence, we have $$\mathbb{E}^{y}[\Psi(Y^{x}_t)]=\mathbb{E}[\Psi(Y^{x,y}_t)], \quad  t \geq 0, y \in V,$$
for all bounded function $\Psi$. For more details on $\mathbb{Q}^{y}$, the reader is referred to \cite{oksendal2003stochastic} . Let $\mathscr{F}_{t}^{x}$ be the $\sigma$-field generated by $\{Y^{x,y}_{r},r\leq t\}$
and set
$$\mathcal{J}_{k}(s,\tau,x,y)=\mathbb{E}[\langle S_{(t-k\delta-s\varepsilon)}(b(x,Y^{x,y}
_s)-\bar{b}(x)), S_{(t-k\delta-\tau\varepsilon)}(b(x,Y^{x,y}_{\tau})-\bar{b}(x))\rangle].$$
Then, we have
\begin{eqnarray*}
	\mathcal{J}_{k}(s,\tau,x,y)
	&=&\mathbb{E}^{y}[\langle S_{t-k\delta-s\varepsilon}(b(x,Y^{x}_s)-\bar{b}(x)), S_{t-k\delta-\tau\varepsilon}(b(x,Y^{x}_{\tau})-\bar{b}(x))\rangle]\cr
	&=&\mathbb{E}^{y}[\langle S_{t-k\delta-\tau\varepsilon}(b(x,Y^{x}_{\tau})-\bar{b}(x)), S_{t-k\delta-s\varepsilon}\mathbb{E}^{y}[(b(x,Y^{x}_s)-\bar{b}(x))|\mathscr{F}_{\tau}^{x}]\rangle ].
\end{eqnarray*}

To proceed, by invoking the Markov property of $Y^{x,y}$, we have
\begin{eqnarray*}
	\mathcal{J}_{k}(s,\tau,x,y)=
	\mathbb{E}^{y}[\langle S_{t-k\delta-\tau\varepsilon}(b(x,Y^{x}_{\tau})-\bar{b}(x)), S_{t-k\delta-s\varepsilon}\mathbb{E}^{Y^{x,y}_{\tau}}[b(x,Y^{x}_{s-\tau})-\bar{b}(x)]\rangle],	
\end{eqnarray*}
where $\mathbb{E}^{Y^{x,y}_{\tau}}[b(x,Y^{x}_{s-\tau})-\bar{b}(x)]$ means the function $\mathbb{E}^{y}[b(x,Y^{x}_{s-\tau})-\bar{b}(x)]$ evaluated at $y=Y^{x,y}_{\tau}$.

Using first H\"{o}lder's inequality, then the contraction property of $S_t$ and the boundedness of the function $b$, we obtain
\begin{eqnarray*}
	\mathcal{J}_{k}(s,\tau,x,y)\leq C\big(
	\mathbb{E}^{y}[|b(x,Y^{x}_{\tau})-\bar{b}(x)|^{2}]\big)^{\frac{1}{2}}\big(
	\mathbb{E}^{y}\big[\big|\mathbb{E}^{Y^{x,y}_{\tau}}[b(x,Y^{x}_{s-\tau})-\bar{b}(x)]\big|^{2}\big]\big)^{\frac{1}{2}},
\end{eqnarray*}
where $C>0$ is a constant independent of $k,s,\tau$. Then, in view of Eq. (\ref{measure}) and Eq. (\ref{yxz}), we have
\begin{eqnarray}\label{jk}
\begin{split}
\mathcal{J}_{k}(s,\tau,x,y)\leq C (1+|x|^{2}+|y|^{2})e^{-\frac{\eta}{2}(s-\tau)}.
\end{split}	
\end{eqnarray}

Let $\mathscr{M}_{k\delta}^{\varepsilon}$ be the $\sigma$-field generated by $X^{\varepsilon}_{k\delta}$ and $\hat{Y}^{\varepsilon}_{k\delta}$ that is independent of $(Y^{x,y}_{r})_{r\geq 0}$. By adopting the approach in \cite[Theorem 7.1.2]{oksendal2003stochastic} . We can show
\begin{eqnarray*}
	&&	\mathcal{J}_{k}(s,\tau)\cr
	&&\quad=\mathbb{E}\big[ \mathbb{E}[\langle S_{t-k\delta-s\varepsilon}(b(X^{\varepsilon}_{k\delta},Y^{X^{\varepsilon}_{k\delta},\hat{Y}^{\varepsilon}_{k\delta}}_s)-\bar{b}(X^{\varepsilon}_{k\delta})), S_{t-k\delta-\tau\varepsilon}(b(X^{\varepsilon}_{k\delta},Y^{X^{\varepsilon}_{k\delta},\hat{Y}^{\varepsilon}_{k\delta}}_{\tau})-\bar{b}(X^{\varepsilon}_{k\delta}))\rangle|\mathscr{M}_{k\delta}^{\varepsilon}]]\cr
	&&\quad=\mathbb{E}[\mathcal{J}_{k}(s,\tau,x,y)|_{(x,y)=(X^{\varepsilon}_{k\delta},\hat{Y}^{\varepsilon}_{k\delta})}\big],
\end{eqnarray*}
which, with the aid of Eq. (\ref{jk}), yields
\begin{eqnarray*}
	\mathcal{J}_{k}(s,\tau)
	\leq C\mathbb{E}[
	(1+|X^{\varepsilon}_{k\delta}|^{2}+|\hat{Y}^{\varepsilon}_{k\delta}|^{2})]e^{-\frac{\eta}{2}(s-\tau)}.
\end{eqnarray*}
This completes  the proof of Lemma \ref{claim}. \qed

\section*{Acknowledgement}
B. Pei was partially supported by National Natural Science Foundation (NSF) of China under Grants No.11802216 and No.12172285, NSF of Chongqing under Grant No.cstc2021jcyj-msxmX0296, Fundamental Research Funds for the Central Universities, Young Talent fund of University Association for Science and Technology in Shaanxi, China, and JSPS Grant-in-Aid for JSPS Fellows under Grant No.JP18F18314. Y. Inahama was partially supported by JSPS KAKENHI under Grant No.JP15K04922 and Grant-in-Aid for JSPS Fellows under Grant No.JP18F18314. Y. Xu was partially supported by NSF of China under Grant No.12072264, Key International (Regional) Joint Research
Program of NSF of China under Grant No.12120101002, Research Funds for Interdisciplinary Subject of Northwestern Polytechnical University, and Shaanxi Provincial Key R\&D Program under Grants No.2019TD-010 and No.2020KW-013. B. Pei would like to thank JSPS for Postdoctoral Fellowships for Research in Japan (Standard).



\begin{thebibliography}{99}
	\setlength{\baselineskip}{0.08in}	
	\parskip=0pt
	\bibitem{Bao2016Two}
	J. Bao, G. Yin, C. Yuan,
	\newblock Two-time-scale stochastic partial differential equations driven by alpha-stable noise: averaging principles. \newblock {\em Bernoulli}, 23(1):645--669, 2017.
	
	\bibitem{Biagini2008stochastic}
	F. Biagini, Y. Hu, B. \O ksendal, et al.,  \newblock {\em Stochastic Calculus for Fractional Brownian Motion and Applications}. \newblock Springer Science \& Business Media, 2008.
	
	\bibitem{brehier2012strong}
	C. Br{\'e}hier, \newblock Strong and weak orders in averaging for SPDEs. \newblock {\em Stochastic Processes and Their Applications}, 122(7):2553--2593, 2012.
	
	\bibitem{caraballo2011existence}
	T. Caraballo, M. Garrido-Atienza, T. Taniguchi,
	\newblock The existence and exponential behavior of solutions to stochastic delay evolution equations with a fractional Brownian motion. \newblock {\em Nonlinear Analysis: Theory, Methods \& Applications}, 74(11):3671--3684, 2011.
	
	\bibitem{cerrai2009averaging}
	S. Cerrai and M. Freidlin, \newblock Averaging principle for a class of stochastic reaction--diffusion equations.
	\newblock {\em Probability Theory and Related Fields}, 144(1):137--177, 2009.
	
	\bibitem{chen2014pathwise}
	Y. Chen, H. Gao, M. J. Garrido-Atienza, B. Schmalfuss, \newblock Pathwise solutions of SPDEs
	and random dynamical systems, \newblock {\em Discrete \& Continuous Dynamical Systems}, 34:79--98, 2014.
	
	\bibitem{da2014stochastic}
	G. Da~Prato, J. Zabczyk,
	\newblock {\em Stochastic Equations in Infinite Dimensions}, \newblock Cambridge University Press, Cambridge, 2014.	
	
	\bibitem{decreusefond1998fractional}
	L. Decreusefond, A. {\"U}st{\"u}nel,
	\newblock Fractional Brownian motion: Theory and applications. \newblock {\em ESAIM: Proceedings}, 5:75--86, 1998.
	
	
	\bibitem{freidlin2012random}
	\newblock M. Freidlin, A. Wentzell,
	\newblock \emph{ Random Perturbations of Dynamical Systems}, \newblock Springer, New York, 2012.
	
	\bibitem{fu2011strong}
	H.~Fu, J.~Liu, \newblock Strong convergence in stochastic averaging principle for two time-scales stochastic partial differential equations. \newblock {\em Journal of Mathematical Analysis and Applications}, 384(1):70--86, 2011.
	
	\bibitem{fu2015strong}
	H. Fu, L. Wan, J. Liu,
	\newblock Strong convergence in averaging principle for stochastic hyperbolic--parabolic equations with two time-scales. \newblock {\em Stochastic Processes and their Applications}, 125(8):3255--3279,
	2015.
	
	\bibitem{GarridoAtienza2010dcds}
	M. Garrido-Atienza, K. Lu, B. Schmalfuss,  \newblock Random dynamical systems for stochastic partial differential equations driven by a fractional Brownian motion, \newblock {\em Discrete \& Continuous Dynamical Systems-B,} 14(2): 473--493, 2010.
	
	\bibitem{Garrido-Atienzalu2010}
	M. Garrido-Atienza, K. Lu, B. Schmalfuss,  \newblock Unstable invariant manifolds for stochastic PDEs driven by a fractional Brownian motion, \newblock {\em Journal of Differential Equations,} 248:1637--1667, 2010.
	
	\bibitem{Random2010garr}
	M. Garrido-Atienza, B. Maslowski, B. Schmalfuss, \newblock Random attractors for stochastic
	equations driven by a fractional Brownian motion, \newblock {\em International Journal of Bifurcation and Chaos}, 20: 2761--2782, 2010.
	
	
	\bibitem{Givon2007}	
	D. Givon. \newblock
	\newblock Strong convergence rate for two-time-scale jump-diffusion stochastic differential systems.
	\newblock \emph{ Multiscale Modeling \& Simulation}, 6: 577--594, 2007.
	
	\bibitem{Guerra2008}
	J. Guerra, D. Nualart,  \newblock Stochastic differential equations driven by fractional Brownian motion and standard Brownian motion.  \newblock {\em Stochastic Analysis and Applications,} 26(5): 1053--1075, 2008.
	
	\bibitem{Hairer2019averaging}
	M. Hairer, X.-M. Li,
	\newblock  Averaging dynamics driven by fractional Brownian
	motion. \newblock {\em Annals of Probability},  48:1826--1860, 2020.
	
	
	\bibitem{khas1966limit}
	R. Khasminskii, \newblock On an averaging principle for It\^{o}  stochastic differential equations.
	\newblock {\em Kibernetica},  4: 260--279,1968.
	
	\bibitem{Liu2019}
	W. Liu, M. R\"{o}ckner, X. Sun, Y. Xie, \newblock Averaging principle for slow-fast stochastic differential equations with time dependent locally Lipschitz coefficients. \newblock {\em Journal of Differential Equations}, 268(6):2910--2948, 2020.
	
	\bibitem{liu2010strong}
 D. Liu, \newblock Strong convergence of principle of averaging for multiscale stochastic dynamical systems.  \newblock {\em  Communications in Mathematical Sciences,} 8(4): 999--1020, 2010.
	
	
	\bibitem{Mandelbrot1968}
	B. Mandelbrot, J. Van Ness,  \newblock Fractional Brownian motions, fractional noises and applications,  \newblock {\em SIAM Review,} 10: 422--427, 1968.
	
	\bibitem{maslowski2003evolution}
	B. Maslowski, D. Nualart,
	\newblock Evolution equations driven by a fractional Brownian motion. \newblock {\em Journal of Functional Analysis}, 202(1):277--305, 2003.
	
	\bibitem{mishura2008stochastic}
	Y. Mishura, \newblock {\em Stochastic Calculus for Fractional Brownian Motion and Related Processes}, \newblock Springer, Berlin, 2008.
	
	\bibitem{Mishura2012mixed}
	\newblock Y. Mishura, G. Shevchenko, \newblock Mixed stochastic differential equations with long-range dependence: Existence, uniqueness and convergence of solutions. \newblock {\em Computers \& Mathematics with Applications,} 64 (10): 3217--3227, 2012.
	
	\bibitem{Monin1967Statistical}
	\newblock A. Monin, A. Yaglom, \newblock {\em Statistical Hydromechanics. Mechanics
		of Turbulence. I-II.} \newblock (Russian) Nauka, Moscow, 1967.
	
	\bibitem{oksendal2003stochastic}
	\newblock B.~{\O}ksendal, \newblock {\em Stochastic Differential Equations}. \newblock Springer, Heidelberg, 2003.
	
	
	\bibitem{pei2020averaging}
B. Pei, Y. Inahama, Y. Xu,
	\newblock  Averaging principles for mixed fast-slow systems driven by fractional Brownian motion. \newblock {\em To appear in Kyoto Journal of Mathematics, 23 pages, 2021.}
	
	\bibitem{pei2020convergence}
	B. Pei, Y. Xu, Y. Bai, \newblock Convergence of $p$-th mean in an averaging principle for stochastic partial differential equations driven by fractional Brownian motion. \newblock {\em Discrete \& Continuous Dynamical Systems-Series B}, 25:1141--1158, 2020.
	
	\bibitem{peistochastic2017}
	B. Pei, Y. Xu, G. Yin, \newblock Stochastic averaging for a class of two-time-scale systems of stochastic partial differential equations. \newblock {\em Nonlinear Analysis}, 160: 159--176, 2017.
	
	\bibitem{pei2018sd}
	B. Pei, Y. Xu, G. Yin, \newblock Averaging principles for SPDEs driven by fractional Brownian motions with random delays modulated by two-time-scale Markov switching processes. \newblock {\em Stochastics and Dynamics},18(3):1850023, 2018.
	
	\bibitem{Averaging2018pei}
	B. Pei, Y. Xu, G. Yin, X. Zhang, \newblock Averaging principles for  functional stochastic partial differential equations driven by a fractional Brownian motion modulated by two-time-scale Markovian switching processes,  \newblock {\em Nonlinear Analysis: Hybrid Systems}, 27: 107--124, 2018.
	
	\bibitem{Rascanu2002}
	D. Nualart, A. R\u{a}\c{s}canu,  \newblock Differential equations driven by fractional Brownian motion.  \newblock {\em Collectanea Mathematica,} 53(1): 55--81,2002.
	
	\bibitem{Shevchenko2014mixed}
	G. Shevchenko, \newblock Mixed fractional stochastic differential equations with jumps. \newblock {\em Stochastics}, 86(2): 203--217, 2014.
	
	\bibitem{sun2020baveraging}
	X. Sun, J. Zhai,
	\newblock Averaging principle for stochastic real Ginzburg-Landau equation driven by $\alpha$-stable process. \newblock {\em Communications on Pure \& Applied Analysis}, 19: 1291--1319, 2020.
	
	\bibitem{Taqqu2003fractional}
M. Taqqu, \newblock {\em  Fractional Brownian motion and long-range dependence.}
	\newblock  In: Doukhan, P., Oppenheim, G., Taqqu, M. (eds) Theory
	and Applications of Long-Range Dependence. Birkh\"auser,
	Boston, 5--38, 2003.
	
	\bibitem{tindel2003stochastic}
	S.~Tindel, C.~Tudor,  F.~Viens,
	\newblock Stochastic evolution equations with fractional Brownian motion. \newblock {\em Probability Theory and Related Fields}, 127(2):186--204, 2003.
	
	\bibitem{thompson2015stochastic}
	\newblock W. Thompson, R. Kuske, A. Monahan, \newblock Stochastic averaging of dynamical systems with multiple time scales forced with $\alpha$-stable noise. \newblock {\em Multiscale Modeling \& Simulation}, 13~(4):1194--1223, 2015.
	
	\bibitem{xu2011averaging}
	Y. Xu, J. Duan,  W. Xu,
	\newblock An averaging principle for stochastic dynamical systems with L{\'e}vy noise.
	\newblock {\em Physica D: Nonlinear Phenomena}, 240(17):1395--1401, 2011.
	
	\bibitem{xu2014}
 Y. Xu, R. Guo, D. Liu, H. Zhang, J. Duan,
	\newblock Stochastic averaging principle for dynamical systems with fractional Brownian motion.
	\newblock {\em Discrete \& Continuous Dynamical Systems-Series B}, 19(4):1197--1212, 2014.
	
	\bibitem{xu2015strong}
	J. Xu, Y.~Miao, J. Liu, \newblock Strong averaging principle for slow-fast SPDEs with Poisson random measures. \newblock {\em Discrete \& Continuous Dynamical Systems-Series B}, 20(7): 2233--2256, 2015.
	
	\bibitem{xu2015approximation}
	Y. Xu, B. Pei, Y. Li,
	\newblock Approximation properties for solutions to non-Lipschitz stochastic differential equations with L{\'e}vy noise. \newblock {\em Mathematical Methods in the Applied Sciences}, 38(11): 2120--2131, 2015.
	
	\bibitem{Xu2015Stochastic}
	Y. Xu, B. Pei, R. Guo,	\newblock Stochastic averaging for slow-fast dynamical systems with fractional Brownian motion. \newblock {\em Discrete \& Continuous Dynamical Systems-Series B}, 20(7):2257--2267, 2015.
	
	\bibitem{Xu2017Stochastic}
	Y. Xu, B. Pei,  J. Wu,
	\newblock Stochastic averaging principle for differential equations with non-Lipschitz coefficients driven by fractional Brownian motion. \newblock {\em Stochastics \& Dynamics}, 17(02): 1750013,2017.
	
	\bibitem{Zahle1998}
	M. Z\"{a}hle,  \newblock Integration with respect to fractal functions and stochastic calculus. I.  \newblock {\em Probability Theory and Related Fields}, 111(3): 333--374,1998.
\end{thebibliography}
\end{document}